\documentclass[11pt,a4paper]{article}
\RequirePackage{amsmath}%
\RequirePackage{amsthm}%
\usepackage{amsfonts}%
\usepackage{amssymb}%
\usepackage{graphicx}%
\usepackage{cite}
\usepackage{url}

\newtheorem{example}{Example}

\renewcommand{\i}{\mathrm{i}}
\renewcommand{\Re}{\mathop{\mathrm{Re}}}
\renewcommand{\Im}{\mathop{\mathrm{Im}}}

\newcommand{\CC}{{\mathbb C}}

\newcommand{\bI}{{\bf I}}

\newcommand{\bM}{{\bf M}}
\newcommand{\bN}{{\bf N}}

\begin{document}

\title{Fast computation of the circular map}
\author{Mohamed M.S. Nasser  }

\date{}
\maketitle

\vskip-0.75cm %
\centerline{Department of Mathematics, Faculty of Science, King Khalid University,} 
\centerline{P. O. Box 9004, Abha 61413, Saudi Arabia}
\centerline{E-mail: mms\_nasser@hotmail.com}

\vskip0.7cm

\begin{center}
\begin{quotation}
{\noindent {\bf Abstracts.\;\;}%
This paper presents a new numerical implementation of Koebe's iterative method for computing the circular map of bounded and unbounded multiply connected regions of connectivity $m$. The computational cost of the method is $O(mn\ln n)$ where $n$ is the number of nodes in the discretization of each boundary component. The accuracy and efficiency of the presented method are demonstrated by several numerical examples. These examples include regions with high connectivity, regions whose boundaries are closer together, and regions with piecewise smooth boundaries.
}%
\end{quotation}
\end{center}
\begin{center}
\begin{quotation}
{\noindent {\bf Keywords.\;\;}%
Numerical conformal mapping, Generalized Neumann kernel, Koebe's method.
}%
\end{quotation}
\end{center}
\begin{center}
\begin{quotation}
{\noindent {\bf MSC.\;\;}  30C30; 65R20.}
\end{quotation}
\end{center}
%

\section{Introduction}
\label{sc:int}

Numerous canonical regions have been considered in the literature for conformal mapping of multiply connected regions in the extended complex plane $\overline{\CC}=\CC\cup\{\infty\}$. Thirty-nine canonical slit regions have been catalogued by Koebe~\cite{Koe18}. A novel method for computing the conformal mapping from bounded and unbounded multiply connected regions onto these thirty-nine canonical regions has been presented in~\cite{Nas-cmft09,Nas-siam09,Nas-jmaa11,Nas-jmaa13}. The method has also been used to compute the conformal mapping onto the canonical region obtained by removing rectilinear slits from an infinite strip~\cite{Nas-siam13}. The method is based on a uniquely solvable boundary integral equation with the generalized Neumann kernel. Only the right-hand side of the integral equation is different from one canonical region to another. A fast method for solving the integral equation with the generalized Neumann kernel is given in~\cite{Nas-siam13,Nas-fast}. For multiply connected regions of connectivity $m$, the method requires $O(mn\ln n)$ operations where $n$ is the number of nodes in the discretization of each boundary component. The method presented in~\cite{Nas-siam09,Nas-siam13} can be used to map bounded and unbounded simply connected regions ($m=1$) onto the unit disk and the exterior unit disk, respectively, in $O(n\ln n)$ operations.

An important canonical region which has not been considered in~\cite{Nas-cmft09,Nas-siam09,Nas-jmaa11,Nas-jmaa13,Nas-siam13} is the multiply connected circular region, i.e., a region all of whose boundaries are circles. The canonical multiply connected circular region is important from physically and computationally point of view. For example, recently, analytic formulas for several problems in fluid mechanics are given for multiply connected circular regions. These analytic formulas are described in terms of \emph{the Schottky-Klein prime function} associated with the circular region~\cite{Cro-lift,Cro-cyl,Cro-Geo,Cro-Ass,Cro-Anz,Cro-Kir,Cro-pri}. Circular regions also are an ideal region for using Fourier series and FFT~\cite{Ben07,Del99,Weg01c,Weg01r}. 

For the canonical multiply connected circular region, the known numerical methods are only iterative methods~\cite{Ben07,Del99,Gai64,Gu11,Hal79,Hen3,Klo96,Koe10,Kro14,Luo10,Mar11,Weg01c,Weg05,Zha12}. Koebe's iterative method is the first numerical method for computing the conformal mapping from multiply connected regions on onto the canonical multiply connected circular region~\cite{Koe10}. A convergence proof and the rate of the convergence for Koebe's can be found in~\cite{Gai64,Hen3}. Koebe's iterative method can be used for bounded and unbounded multiply connected regions. For bounded regions of connectivity $m$, each iteration of Koebe's method requires computing the conformal mapping from an unbounded simply connected region onto the exterior unit disk for $m-1$ times and computing the conformal mapping from a bounded simply connected region onto the unit disk for one time. For unbounded regions of connectivity $m$, each iteration of Koebe's method requires computing the conformal mapping from an unbounded simply connected region onto the exterior unit disk for $m$ times. The successive computational region becomes gradually more circular. In the process, we compute approximate values of the centres and radii of the circles. For numerical implementations of Koebe's iterative method, see~\cite{Gu11,Kro14,Luo10,Mar11,Zha12}. Other numerical method for circular regions are Wegmann iterative method~\cite{Weg01c} and Fornberg-like iterative method~\cite{Ben07,Del99} which can be used to compute the inverse conformal mapping form the circular region onto the multiply connected regions. These methods are based on using trigonometric interpolation and FFT. A comparison between these two methods is given in~\cite{Ben07}. For multiply connected regions of connectivity $m$, if $n$-point trigonometric interpolation is used, the computational cost of these methods is $O((mn)^2)$~\cite{Ben07}.

In this paper, based on a boundary integral equation with the generalized Neumann kernel, we present a new numerical implementation of Koebe's iterative method for conformally mapping bounded and unbounded multiply connected regions of connectivity $m$ onto bounded and unbounded canonical multiply connected circular regions of connectivity $m$, respectively. The method provide us with the boundary values of the conformal mapping and its derivative. The interior values of the mapping function are calculated using the Cauchy integral formula. Cauchy integral formula can be also used to calculate the interior values of the derivative of the mapping function as well as the interior values of the inverse mapping function. The computational cost of the method is $O(mn\ln n)$ where $n$ is the number of nodes in the discretization of each boundary component. 

The remainder of this paper is organized as follows: the circle map is defined in Section~2. In Section 3, we present a fast numerical method for computing the conformal mapping of simply connected regions. This fast method with Koebe's iterative methods will be used to compute the circular map of bounded and unbounded multiply connected regions in Sections~4 and~5, respectively. In Section 6, we present eight numerical examples. A short conclusion is given in Section~7.

\section{The circular map}
\label{sc:map}

\begin{figure}[ht]%
\centerline{\scalebox{0.6}{\includegraphics{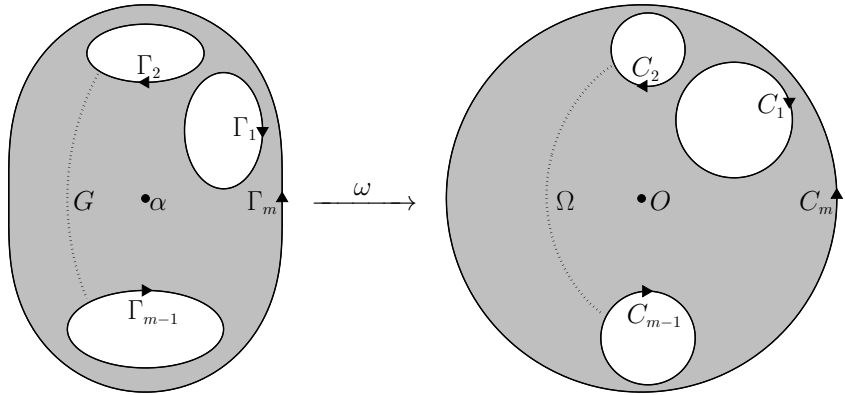}}}
\caption{\rm The original bounded multiply connected region $G$ of connectivity $m$ (left) and the bounded circular region (right).} 
\label{f:Bm}
\end{figure}

\begin{figure}[ht]%
\centerline{\scalebox{0.6}{\includegraphics{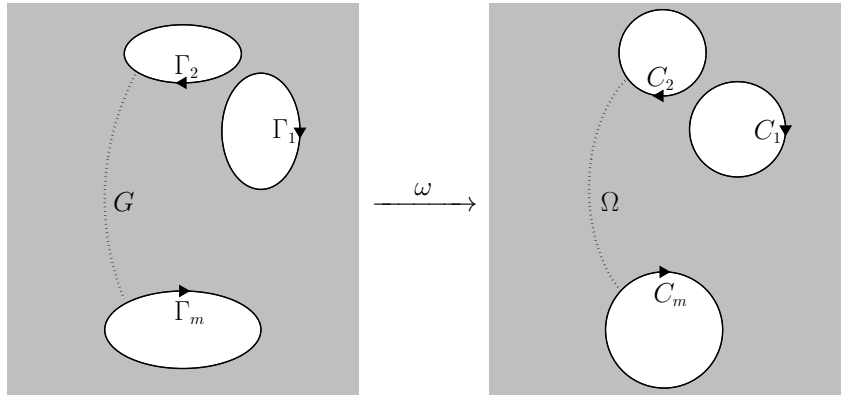}}}
\caption{\rm The original unbounded multiply connected region $G$ of connectivity $m$ (left) and the unbounded circular region (right).} 
\label{f:Um}
\end{figure}

Let $G$ be a multiply connected region of connectivity $m$ in the extended complex plane $\overline{\CC}=\CC\cup\{\infty\}$. The region $G$ can be bounded or unbounded. For bounded $G$, we assume that $\alpha$ is a fixed point in $G$. See Figure~\ref{f:Bm}(left). If $G$ is unbounded, then we assume that $\infty\in G$ and $\alpha$ is a fixed point in the complement of $G$. See Figure~\ref{f:Um}(left). Let $G$ has the boundary
\[
\Gamma=\partial G= \bigcup_{j=1}^{m} \Gamma_j
\]
where $\Gamma_1, \ldots, \Gamma_m$ are closed smooth Jordan curves. The orientation of $\Gamma$ is such that $G$ is always on the left of~$\Gamma$. The curve $\Gamma_j$ is parametrized by a $2\pi$-periodic twice continuously differentiable complex function $\eta_j(t)$ with non-vanishing first derivative $\eta'_j(t)\ne 0$ for $t\in J_j=[0,2\pi]$. The total parameter domain $J$ is the disjoint union of the $m$ intervals $J_1,\ldots,J_m$. We define a parametrization of the whole boundary $\Gamma$ as the complex function $\eta$ defined on $J$ by
\begin{equation}\label{e:eta}
\eta(t)= \left\{ \begin{array}{l@{\hspace{0.5cm}}l}
\eta_1(t),&t\in J_1,\\
\hspace{0.3cm}\vdots\\
\eta_m(t),&t\in J_m. \\
\end{array}
\right.
\end{equation}

For bounded $G$, there exists a conform mapping $w=\omega(z)$ from the bounded region $G$ onto a bounded region $\Omega$ which is bounded by $m$ circles (see Figure~\ref{e:cond-b}(right)). The centre of the exterior circle is $z_m=0$ and the radius is $r_m=1$, i.e., the external circle is the unit circle. For the inner circles $C_j$ for $j=1,2,\ldots,m-1$, the centres $z_j$ and the radii $r_j$ are unknown and should be determined. When the conformal mapping $w=\omega(z)$ is normalized by 
\begin{equation}\label{e:cond-b}
\omega(\alpha)=0, \quad \omega'(\alpha)>0,
\end{equation}
then the conformal mapping $\omega$ as well as the circular region $\Omega$ are uniquely determined by the region $G$~\cite{Hen3,Weg01c}. 

For unbounded $G$, there exists a conform mapping $w=\omega(z)$ from the unbounded region $G$ onto an unbounded circular region $\Omega$ which is bounded by $m$ circles. The centres $z_j$ and the radii $r_j$ of the circles $C_j$ for $j=1,2,\ldots,m$ are unknown and should be determined. When the conformal mapping $w=\omega(z)$ is normalized by the condition near infinity
\begin{equation}\label{e:cond-u}
\omega(z)=z+O\left(\frac{1}{z}\right),
\end{equation}
then the conformal mapping $\omega$ as well as the circular region $\Omega$ are uniquely determined by the region $G$~\cite{Ben07,Del99,Gol69,Hen3,Weg01c}.

For both bounded and unbounded $G$, the region $\Omega$ is called a \emph{circular region} and the mapping function $\omega$ is called the \emph{circular map} of $G$~\cite[p.~488]{Hen3}. The centres and the radii of the circles are called the \emph{parameters of the canonical region} $\Omega$.

\section{The conformal mapping of simply connected regions}
\label{sc:sim}

In this section, we shall present a fast numerical method for computing the conformal mapping from the bounded simply connected region onto the unit disk (see Figure~\ref{f:Bs}). The method can also be used for computing the conformal mapping from the unbounded simply connected region onto the exterior unit disk (see Figure~\ref{f:Us}). This method with Koebe's iterative method will be used in~\S\ref{sc:Koe-bd} and~\S\ref{sc:Koe-ud} to compute the conformal mapping from the bounded and unbounded multiply connected region $G$ onto the bounded and unbounded multiply connected circular region $\Omega$, respectively.

\subsection{The simply connected region}

\begin{figure}[ht]%
\centerline{\scalebox{0.6}{\includegraphics{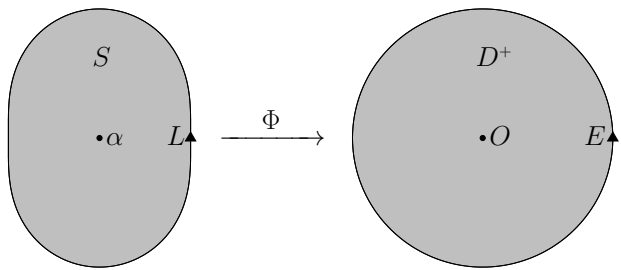}}}
\caption{\rm The original bounded simply connected region $S$ (left) and the bounded circular region $D^+$ (right).} 
\label{f:Bs}
\end{figure}

\begin{figure}[ht]%
\centerline{\scalebox{0.6}{\includegraphics{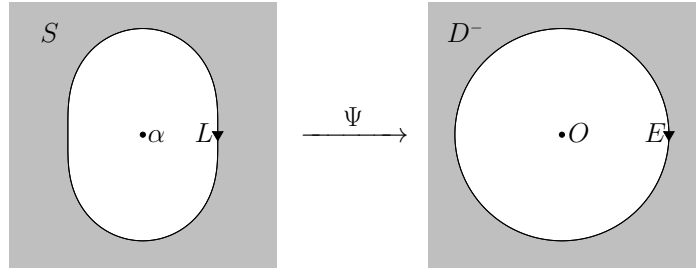}}}
\caption{\rm The original unbounded simply connected region $S$ (left) and the unbounded circular region $D^-$ (right).} 
\label{f:Us}
\end{figure}

Let $S$ be a simply connected region in the extended complex plane $\overline{\CC}=\CC\cup\{\infty\}$. The region $S$ can be bounded or unbounded. For bounded $S$, we assume that $\alpha$ is a fixed point in $S$ (see Figure~\ref{f:Bs}(left)). If $S$ is unbounded, then we assume that $\infty\in S$ and $\alpha$ is a fixed point in the complement of $S$ (see Figure~\ref{f:Us}(left)). The boundary $L=\partial S$ is assumed to be a closed smooth Jordan curves. The orientation of $L$ is such that $S$ is always on the left of~$L$, i.e., $L$ is counterclockwise oriented for bounded $S$ and clockwise oriented for unbounded $S$. 

\subsection{The generalized Neumann kernel}

The curve $L$ is parametrized by a $2\pi$-periodic twice continuously differentiable complex function $\zeta(t)$ with non-vanishing first derivative $\zeta'(t)\ne 0$ for $t\in [0,2\pi]$. We define a complex-valued function $A$ on $L$ by
\begin{equation}\label{e:A}
A(t) = \left\{
\begin{array}{l@{\hspace{0.5cm}}l}
 \zeta(t)-\alpha,     & \text{if $S$ is bounded}, \\
 1,     & \text{if $S$ is unbounded}.
\end{array}%
\right.
\end{equation}
The generalized Neumann kernel formed with $A$ and $\zeta$ is defined by~\cite{Mur-Bul,Weg-Mur-Nas}
\begin{equation}\label{e:N}
 N(s,t) =  \frac{1}{\pi}\Im\left(
 \frac{A(s)}{A(t)}\frac{\zeta'(t)}{\zeta(t)-\zeta(s)}\right).
\end{equation}
We define also a kernel
\begin{equation}\label{e:M}
 M(s,t) =  \frac{1}{\pi}\Re\left(
 \frac{A(s)}{A(t)}\frac{\zeta'(t)}{\zeta(t)-\zeta(s)}\right).
\end{equation}
The kernel $N$ is continuous and the kernel $M$ is singular. Thus, the integral operator
\begin{equation}\label{e:bN}
  \bN \mu(s) = \int_{0}^{2\pi} N(s,t) \mu(t) dt, \quad s\in [0,2\pi],
\end{equation}
is a Fredholm integral operator and the operator
\begin{equation}\label{e:bM}
  \bM\mu(s) = \int_{0}^{2\pi}  M(s,t) \mu(t) dt, \quad s\in [0,2\pi],
\end{equation}
is a singular integral operator. 

For more details on generalized Neumann kernel, see~\cite{Mur-Bul,Nas-cmft09,Nas-siam09,Weg-Mur-Nas,Weg-Nas}.

\subsection{The bounded simply connected region}

Let $\Phi$ be the conformal mapping from the bounded simply connected region $S$ onto the unit disk $D^+$. The boundary values of $\Phi$ are given by
\begin{equation}\label{e:Phi-thet}
\Phi(\zeta(t))=e^{\i\theta(t)}
\end{equation}
where $\theta(t)$ is the \emph{boundary correspondence function} of the mapping function $\Phi$. The function $\theta(t)-t$ is a $2\pi$-periodic function and $\theta'(t)>0$ for all $t\in[0,2\pi]$. By differentiating both sides of~(\ref{e:Phi-thet}) with respect to the parameter $t$, we obtain
\begin{equation}\label{e:Phi'-thet1}
\zeta'(t)\Phi'(\zeta(t))=\i\theta'(t)e^{\i\theta(t)}.
\end{equation}
Hence
\begin{equation}\label{e:Phi'-thet2}
\zeta'(t)\frac{\Phi'(\zeta(t))}{\Phi(\zeta(t))}=\i\theta'(t).
\end{equation}
It is clear from~(\ref{e:Phi-thet}) and (\ref{e:Phi'-thet1}) that determining the functions $\theta$ and $\theta'$ provides us with the boundary values of the mapping function $\Phi$ and its derivative $\Phi'$.

With the normalization
\begin{equation}\label{e:Phi-cond}
\Phi(\alpha)=0, \quad \Phi'(\alpha)>0,
\end{equation}
the mapping function $\Phi$ is unique. The function $w=\Phi(z)$ can be written as
\begin{equation}\label{e:Phi-f}
\Phi(z)=c(z-\alpha)e^{(z-\alpha)f(z)},
\end{equation}
where $f$ is analytic function on $S$ and $c=\Phi'(\alpha)>0$. Hence
\begin{equation}\label{e:Phi-log}
\log\Phi(z)=\ln c+\log(z-\alpha)+(z-\alpha)f(z)
\end{equation}
Then, in view of~(\ref{e:Phi-thet}), we obtain
\begin{equation}\label{e:thet-f}
\i\theta(t)=\ln c+\log(\zeta(t)-\alpha)+A(t)f(\zeta(t)).
\end{equation}
Hence, the boundary values of the analytic function $f$ are given by
\begin{equation}\label{e:Af}
A(t)f(\zeta(t))=\gamma(t)+h+\i(\mu(t)+\theta(t))
\end{equation}
where $h=-\ln c$ and
\begin{equation}\label{e:gam-mu}
\gamma(t)+\i\mu(t)=-\log(\zeta(t)-\alpha).
\end{equation}

By differentiating both sides of~(\ref{e:Phi-log}), we obtain
\begin{equation}\label{e:Phi'-f'}
\frac{\Phi'(z)}{\Phi(z)}=\frac{1}{z-\alpha}+f(z)+(z-\alpha)f'(z).
\end{equation}
Then function
\begin{equation}\label{e:F-Phi'}
F(z)=(z-\alpha)\frac{\Phi'(z)}{\Phi(z)}
\end{equation}
is analytic is $S$ and its boundary values are given by
\begin{equation}\label{e:F-thet'}
\tilde{A}(t)F(\zeta(t))=\i\theta'(t).
\end{equation}

\subsection{The unbounded simply connected region}

Let $\Psi$ be the mapping function from the unbounded simply connected region $S$ onto the exterior unit disk $D^-$. The boundary values of $\Psi$ satisfies
\begin{equation}\label{e:Psi-thet}
\Phi(\zeta(t))=e^{-\i\theta(t)}
\end{equation}
where $\theta(t)$ is the \emph{boundary correspondence function} of the mapping function $\Psi$. The function $\theta(t)-t$ is a $2\pi$-periodic function and $\theta'(t)>0$ for all $t\in[0,2\pi]$. By differentiating both sides of~(\ref{e:Psi-thet}) with respect to $t$, we obtain
\begin{equation}\label{e:Psi'-thet1}
\zeta'(t)\Psi'(\zeta(t))=-\i\theta'(t)e^{-\i\theta(t)}.
\end{equation}
Hence
\begin{equation}\label{e:Psi'-thet2}
\zeta'(t)\frac{\Psi'(\zeta(t))}{\Psi(\zeta(t))}=-\i\theta'(t).
\end{equation}

With the normalization
\begin{equation}\label{e:Psi-cond}
\Psi(\infty)=0, \quad \Psi'(\infty)>0,
\end{equation}
the mapping function $\Psi$ is unique. The function $w=\Psi(z)$ can be written as
\begin{equation}\label{e:Psi-f}
\Psi(z)=c(z-\alpha)e^{-f(z)},
\end{equation}
where $f$ is analytic function on $S$ and $c=\Psi'(\infty)>0$. Then
\begin{equation}\label{e:Psi-log}
\log\Psi(z)=\ln c+\log(z-\alpha)-f(z).
\end{equation}
Then, in view of~(\ref{e:Psi-thet}), we obtain
\begin{equation}\label{e:thet-f-u}
-\i\theta(t)=\ln c+\log(\zeta(t)-\alpha)-A(t)f(\zeta(t))
\end{equation}
Hence, the boundary values of the function $f$ are given by
\begin{equation}\label{e:Af-u}
A(t)f(\zeta(t))=\gamma(t)+h+\i(\mu(t)+\theta(t))
\end{equation}
where $h=\ln c$ and
\begin{equation}\label{e:gam-mu-u}
\gamma(t)+\i\mu(t)=\log(\zeta(t)-\alpha).
\end{equation}

By differentiating both sides of~(\ref{e:Psi-log}), we obtain
\begin{equation}\label{e:Psi'-f'}
\frac{\Psi'(z)}{\Psi(z)}=\frac{1}{z-\alpha}-f'(z).
\end{equation}
Then function
\begin{equation}\label{e:F-Psi'}
F(z)=\frac{\Psi'(z)}{\Psi(z)}
\end{equation}
is analytic is $S$ with $F(\infty)=0$ and its boundary values are given by
\begin{equation}\label{e:F-thet'-u}
\tilde{A}(t)F(\zeta(t))=\i\theta'(t).
\end{equation}

\subsection{The fast numerical method}
\label{sc:m1}

For both bounded and unbounded $S$, in view of~(\ref{e:Af}) and~(\ref{e:Af-u}), it follows from~\cite{Nas-cmft09,Nas-siam09} the the function $\phi=\mu+\theta$ is the unique solution of the integral equation
\begin{equation}\label{e:ie}
(\bI-\bN)\phi=-\bM\gamma
\end{equation}
and the constant $h$ is given by
\begin{equation}\label{e:h}
h=[(\bI-\bN)\gamma-\bM\phi]/2,
\end{equation}
where the functions $\gamma$ and $\mu$ are given by~(\ref{e:gam-mu}) for bounded $S$ and by~(\ref{e:gam-mu-u}) for unbounded $S$. Thus the boundary correspondence function $\theta$ is given by
\[
\theta=\phi-\mu.
\]

In this paper, we solve the equation~(\ref{e:ie}) by the fast method presented in~\cite{Nas-siam13,Nas-fast}. The numerical method is based on strict discretization of the integrals in~(\ref{e:ie}) and~(\ref{e:h}) by the trapezoidal rule with the $n$ equidistant collocation points 
\begin{equation}\label{e:pt-tj}
t_{i}=(i-1)\frac{2\pi}{n}, \quad i=1,2,\ldots,n,
\end{equation}
for a given even positive integer $n$. For $2\pi$-periodic function $\gamma(t)$, the trapezoidal rule approximate the integral $I=\int_{0}^{2\pi}\gamma(t)dt$ by $I_n=\frac{2\pi}{n}\sum_{i=1}^{n}\gamma(t_i)$. If the function $\gamma(t)$ is $k$ times continuously differentiable, then the rate of convergence of the trapezoidal rule is $O(1/n^k)$. For analytic $\gamma(t)$, the rate of convergence is better than $O(1/n^k)$ for any positive integer $k$~\cite[p.~83]{Kro-Ueb}. 

Discretizing the integral equation~(\ref{e:ie}) yields an $n\times n$ linear system which is solved in $O(n\ln n)$ operations by a combination of the GMRES method and the FMM~\cite{Nas-siam13,Nas-fast}. The GMRES method converges significantly faster since the eigenvalues of the discretizing matrix are clustered around $1$. In fact, since $S$ is a simply connected region and the function $A$ is defined by~(\ref{e:A}), the generalized Neumann kernel $N$ has only real eigenvalues in the interval $[-1,1)$ where $-1$ is a simple eigenvalue~\cite{Nas-amc11}. Thus, the discretizing matrix of the integral equation~(\ref{e:ie}) has only real eigenvalues on the interval $(0,2]$ with $2$ as a simple eigenvalue and the other eigenvalues are clustered around $1$~\cite{Nas-amc11,Nas-Mur12}. For nearly circular region $S$, the eigenvalues become much clustered around $1$. When $L$ is the unit circle, the generalized Neumann kernel becomes
\[
N(s,t) = -\frac{1}{2\pi}
\]
which has only the eigenvalues $0$ and $-1$. Thus the discretizing matrix of the integral equation~(\ref{e:ie}) has only two eigenvalues $2$ and $1$ where $2$ is a simple eigenvalue and $1$ has the algebraic multiplicity $n-1$.

Solving the integral equation with generalized Neumann kernel~(\ref{e:ie}) using the fast method presented in~\cite{Nas-fast} does not require the second derivative $\eta''(t)$ of the parametrization of the boundary. Thus, in view of~(\ref{e:Phi-thet}), (\ref{e:Phi'-thet2}), (\ref{e:Psi-thet}), and (\ref{e:Psi'-thet2}), successive application of the fast method, requires calculating the values of the boundary correspondence function $\theta(t)$ and its first derivative $\theta'(t)$ at the points~(\ref{e:pt-tj}). Thus, by calculating the values of the boundary correspondence function $\theta(t)$ at the points~(\ref{e:pt-tj}), we shall represent the $2\pi$-periodic function $\theta(t)-t$ on $[0,2\pi]$ by the interpolating trigonometric polynomial of degree $n/2$
	\begin{equation}\label{e:thet-F}
	\theta(t)-t=a_0+\sum_{j=1}^{n/2}a_j\cos jt +\sum_{j=1}^{n/2-1}b_j\sin jt, 
	\end{equation}
that interpolate $\theta(t)-t$ at the $n$ equidistant points~(\ref{e:pt-tj}) (see~\cite[p.~364]{Weg05}). The coefficients $a_0,a_1,\ldots,a_n,b_1,\ldots,b_n$ are calculated by the FFT in $O(n\ln n)$ operations. Then the function $\theta'(t)-1$ is approximated by the trigonometric polynomials of degree $n/2$
\begin{equation}\label{e:thet'-F}
\theta'(t)-1=a'_0+\sum_{j=1}^{n/2}a'_j\cos jt +\sum_{j=1}^{n/2-1}b'_j\sin jt, 
\end{equation}
where
\[
a'_0=a'_n=0, \quad a'_j=jb_j, \quad b'_j=-ja_j, \quad j=1,2,\ldots,n-1.
\]
The values of the function $\theta'(t)$ at the $n$ points~(\ref{e:pt-tj}) are calculated by the FFT. 

For bounded region $S$, by obtaining the functions $\theta$ and $\theta'$, we obtain the boundary values of the mapping function $\Phi$ and its derivative $\Phi'$. The functions $\Phi$ and $\Phi'$ are analytic in $S$. Hence the values of the functions $\Phi(z)$ and $\Phi'(z)$ for interior points $z\in S$ can be computed by the Cauchy integral formula
\begin{eqnarray}
\label{e:Phi(z)}
\Phi(z)&=&\frac{1}{2\pi\i}\int_{L}\frac{\Phi(\zeta)}{\zeta-z}d\zeta, \\
\label{e:Phi'(z)}
\Phi'(z)&=&\frac{1}{2\pi\i}\int_{L}\frac{\Phi'(\zeta)}{\zeta-z}d\zeta.
\end{eqnarray}

For unbounded region $S$, obtaining the functions $\theta$ and $\theta'$ yields the boundary values of the mapping function $\Psi$ and its derivative $\Psi'$. The functions $\frac{\Psi(z)}{z-\alpha}$ and $\Psi'(z)$ are analytic in $S$ and have the same value $c=e^h$ at $\infty$ where the constant $h$ is given by~(\ref{e:h}). By the Cauchy integral formula, the values of the functions $\Psi(z)$ and $\Psi'(z)$ for interior points $z\in S$ can be computed by~\cite[p.~2]{Gak66}
\begin{eqnarray}
\label{e:Psi(z)}
\Psi(z)&=&-c(z-\alpha)+\frac{z-\alpha}{2\pi\i}\int_{L}\frac{\Psi(\zeta)}{\zeta-\alpha}\frac{1}{\zeta-z}d\zeta, \\
\label{e:Psi'(z)}
\Psi'(z)&=&-c+\frac{1}{2\pi\i}\int_{L}\frac{\Psi'(\zeta)}{\zeta-z}d\zeta.
\end{eqnarray}

For the computational cost of the numerical method, computing the right-hand side of the integral equation~(\ref{e:ie}), requires two FMMs and three FFTs, each iteration of the GMRES method requires one FMM, computing the function $h$ in~(\ref{e:h}) requires two FMMs and three FFTs, and computing the function $\theta'$ requires two FFTs. Thus, the method requires four FMMs, eight FFTs, and one FMM for each iteration of the GMRES method. Since one application of the FMM requires $O(n)$ operations and one application of the FFT requires $O(n\ln n)$ operations, the complexity of the method is $O(n\ln n)$ operations.

\section{Koebe's iterative method for bounded multiply connected regions}
\label{sc:Koe-bd}

In this section, based on the results of the previous section, a numerical implementation of Koebe's iterative method for computing the conformal mapping $w=\omega(z)$ from the bounded multiply connected region $G$ onto the bounded circular region $\Omega$ (see Figure~\ref{f:Bm}) will be describe. We present a method for computing the boundary values of the mapping function
\[
\xi(t)=\omega(\eta(t)), \quad t\in J,
\]
 the boundary values of its derivative 
\[
\xi'(t)=\eta'(t)\omega'(\eta(t)), \quad t\in J,
\]
and the parameters $z_j,r_j$, $j=1,2,\ldots,m-1$, of the canonical region $\Omega$ where $\eta(t)$, given by~(\ref{e:eta}), is the parametrization of the boundary $\Gamma$ of $G$.

\subsection{Initializations}

At the beginning, we set
\[
C_i^{0,0}=\Gamma_i \quad{\rm for}\quad i=1,2,\ldots,m.
\]
The curve $C_i^{0,0}$ is parametrized by $\xi_i^{0,0}(t)$ which is defined by
\[
\xi_i^{0,0}(t)=\eta_i(t), \quad t\in J_i, \quad i=1,2,\ldots,m.
\]
Hence
\[
\frac{d}{dt}\xi_i^{0,0}(t)=\eta'_i(t), \quad t\in J_i, \quad i=1,2,\ldots,m.
\]
We Assume also that $z_i^{0,0}$ is a fix point inside $C_i^{0,0}$ for $i=1,\ldots,m-1$ and $z_{m}^{0,0}=\alpha$. 
Hence, initial values of the boundary values of the mapping function are given by
\[
\omega^{0}(\eta(t))= \left\{ \begin{array}{l@{\hspace{0.5cm}}l}
\xi_1^{0,0}(t),&t\in J_1,\\
\hspace{0.3cm}\vdots\\
\xi_m^{0,0}(t),&t\in J_m. \\
\end{array}
\right.
\]

\subsection{Iterations}

For $k=1,2,3,\ldots$, where $k$ denotes the iteration number, we repeat the following three steps:

\subsubsection{Step I: The internal curves}
\label{ssc:int}

For $j=1,2,\ldots,m-1$, let $\Psi_{k,j}$ be the conformal mapping from the exterior region of the curve $C^{k-1,j-1}_{j}$ onto the exterior unit disk. Then
\[
C^{k-1,j}_{j} = \Psi_{k,j}(C^{k-1,j-1}_{j})
\]
is the unit circle. The curve $C_j^{k-1,j}$ is parametrized by $\xi_j^{k-1,j}(t)$ which is defined as the boundary values of the conformal mapping $\Psi_{k,j}$, i.e.,
\[
\xi_j^{k-1,j}(t)=\Psi_{k,j}(\xi_j^{k-1,j-1}(t)), \quad t\in J_j.
\]
The derivative of the function $\xi_j^{k-1,j}(t)$ can be computed from the boundary values of the derivative of the mapping function $\Psi_{k,j}$ as follows
\[
\frac{d}{dt}\xi_j^{k-1,j}(t)=\Psi'_{k,j}(\xi_j^{k-1,j-1}(t))\frac{d}{dt}\xi_j^{k-1,j-1}(t),  \quad t\in J_j.
\]
The boundary values of the mapping function $\Psi_{k,j}$ and its derivative $\Psi'_{k,j}$ can be computed using the fast method presented in~\S\ref{sc:m1}. 

The smooth Jordan curves $C^{k-1,j-1}_{i}$, $i=1,2,\ldots,m$, $i\ne j$, are external to the curve $C^{k-1,j-1}_{j}$. Thus, the function $\Psi_{k,j}$ maps the curves $C^{k-1,j-1}_{i}$ onto smooth Jordan curves 
\[
C^{k-1,j}_{i}= \Psi_{k,j}(C^{k-1,j-1}_{i}),\quad i=1,2,\ldots,m, \quad i\ne j,
\]
external to $C^{k-1,j}_{j}$. For $i=1,2,\ldots,m$ such that $i\ne j$, the curve $C_i^{k-1,j}$ is parametrized by 
\[
\xi_i^{k-1,j}(t)=\Psi_{k,j}(\xi_i^{k-1,j-1}(t)), \quad t\in J_i.
\]
The derivative of the function $\xi_i^{k-1,j}(t)$ can be computed from the boundary values of the derivative of the mapping function $\Psi_{k,j}$. Since $\xi_i^{k-1,j-1}(t)$, $t\in J_i$, are in the exterior region of $C^{k-1,j-1}_{j}$, then the values of the function $\Psi_{k,j}$ and its first derivative at the points $\xi_i^{k-1,j-1}(t)$, i.e., $\Psi_{k,j}(\xi_i^{k-1,j-1}(t))$ and $\Psi'_{k,j}(\xi_i^{k-1,j-1}(t))$, can be computed using the Cauchy integral formula as explained in~(\ref{e:Psi(z)}) and~(\ref{e:Psi'(z)}).

Finally, we set
\[
z_j^{k-1,j}=0.
\]
For $i=1,2,\ldots,m$, $i\ne j$, the point $z^{k-1,j-1}_{i}$ inside the curve $C^{k-1,j-1}_{i}$ in the exterior region of the curve $C^{k-1,j-1}_{j}$ will be mapped by the function $\Psi_{k,j}$ into a point $\Psi_{k,j}\left(z^{k-1,j-1}_{i}\right)$ inside the curve $C^{k-1,j}_{i}$ in the exterior of the circle $C^{k-1,j}_{j}$. The values of the function $\Psi_{k,j}$ at the points $z^{k-1,j-1}_{i}$ can be computed by the Cauchy integral formula. We define
\[
z^{k-1,j}_{i}=\Psi_{k,j}\left(z^{k-1,j-1}_{i}\right), \quad i=1,2,\ldots,m, \quad i\ne j.
\]

\subsubsection{Step II: The external curve}
\label{ssc:ext}

Let $\Phi_{k}$ be the conformal mapping from the interior region of the curve $C^{k-1,m-1}_{m}$ onto the unit disk. Then
\[
C^{k-1,m}_{m} = \Phi_{k}(C^{k-1,m-1}_{m})
\]
is the unit circle. The curve $C_{m}^{k-1,{m}}$ is parametrized by $\xi_{m}^{k-1,{m}}(t)$ which is defined as the boundary values of the conformal mapping $\Phi_{k}$, i.e.,
\[
\xi_{m}^{k-1,{m}}(t)=\Phi_{k}(\xi_{m}^{k-1,m-1}(t)), \quad t\in J_{m}.
\]
The derivative of the function $\xi_{m}^{k-1,{m}}(t)$ can be computed from the boundary values of the derivative of the mapping function $\Phi_{k}$. The boundary values of the mapping function $\Phi_{k}$ and its derivative $\Phi'_{k}$ can be computed using the fast method presented in~\S\ref{sc:m1}. 

The smooth Jordan curves $C^{k-1,m-1}_{i}$, $i=1,2,\ldots,m-1$, are internal to the curve $C^{k-1,m-1}_{m}$. Thus, the function $\Phi_{k}$ maps the curves $C^{k-1,m-1}_{i}$ onto smooth Jordan curves 
\[
C^{k-1,m}_{i}= \Phi_{k}(C^{k-1,m-1}_{i}),\quad i=1,2,\ldots,m-1,
\]
internal to $C^{k-1,m}_{m}$. For $i=1,2,\ldots,m-1$, the curve $C_i^{k-1,m}$ is parametrized by 
\[
\xi_i^{k-1,m}(t)=\Phi_{k}(\xi_i^{k-1,m-1}(t)).
\]

The derivative of the function $\xi_{m+1}^{k-1,m}(t)$ can be computed from the boundary values of the derivative of the mapping function $\Phi_{k}$. Since $\xi_i^{k-1,m-1}(t)$, $t\in J_i$, $i=1,2,\ldots,m-1$, are in the interior region of $C^{k-1,m-1}_{m}$, then the values of the function $\Phi_{k}$ and its derivative at the points $\xi_i^{k-1,m-1}(t)$, i.e., $\Phi_{k}(\xi_i^{k-1,m-1}(t))$ and $\Phi'_{k}(\xi_i^{k-1,m-1}(t))$, can be computed using the Cauchy integral formula as explained in~(\ref{e:Phi(z)}) and~(\ref{e:Phi'(z)}).
  
Then, we set
\[
z_{m}^{k-1,{m}}=0.
\]
For $i=1,2,\ldots,m-1$, the point $z^{k-1,m-1}_{i}$ inside the curve $C^{k-1,m-1}_{i}$ in the interior region of the curve $C^{k-1,m-1}_{m}$ are mapped by the function $\Phi_{k}$ into a point $\Phi_{k}\left(z^{k-1,m-1}_{i}\right)$ inside the curve $C^{k-1,m}_{i}$ in the interior of the circle $C^{k-1,m}_{j}$. We define
\[
z^{k-1,m}_{i}=\Phi_{k}\left(z^{k-1,m-1}_{i}\right), \quad i=1,2,\ldots,m-1,
\]
where $\Phi_{k}\left(z^{k-1,m-1}_{i}\right)$ are computed by the Cauchy integral formula.

\subsubsection{Step III: Update and conditions of convergence}
\label{ssc:con}

Let $w=\omega^k(z)$ be the approximate mapping function obtained in the $k^{\rm th}$ iteration. Then the boundary values of $\omega^k$ are given by
\begin{equation}\label{e:om-k}
\omega^{k}(\eta(t))= \left\{ \begin{array}{l@{\hspace{0.5cm}}l}
\xi_1^{k-1,m}(t),&t\in J_1,\\
\hspace{0.3cm}\vdots\\
\xi_m^{k-1,m}(t),&t\in J_m.
\end{array}
\right.
\end{equation}
Then the boundary values of the derivative $\frac{d}{dz}\omega^k(z)$ can be computed by differentiate both sides of~(\ref{e:om-k}).

By obtaining the boundary values of $\omega^k$, we test the convergence of the method. We stop the iteration if
\begin{equation}\label{e:conv}
\|\omega^{k}-\omega^{k-1}\|_\infty<\varepsilon\quad {\rm or}\quad k>Max,
\end{equation}
where $\varepsilon$ is a given tolerance and $Max$ is the maximum number of iterations allowed.

If the condition~(\ref{e:conv}) is not satisfied, we set	
\[
C_i^{k,0}=C_i^{k-1,m}, \quad i=1,\ldots,m.
\]
For $i=1,2,\ldots,m$, the curve $C_i^{k,0}$ is parametrized by
\[
\xi_i^{k,0}(t)=\xi_i^{k-1,m}(t).
\] 
Then, we set $k=k+1$ and repeat Steps I--III.

\subsection{The interior values}

If the method converges, then we consider also the boundary values of the approximate mapping function $\omega^k$ in~(\ref{e:om-k}) as an approximation of the boundary values of the mapping function $\omega$, i.e., we set 
\begin{equation}\label{e:om}
\xi(t)=\omega(\eta(t))=\omega^k(\eta(t)), \quad t\in J.
\end{equation}
We consider the bounded multiply connected region bounded by the circles $C_1^{k-1,m}$, \ldots, $C_m^{k-1,m}$, as the canonical region $\Omega$. The boundaries of $\Omega$ are then given by
\[
C_i=C_i^{k-1,m}, \quad i=1,\ldots,m,m.
\] 
The centre $z_i$ and the radius $R_i$ of the circle $C_i$ are approximated by
\[
z_i=z^{k-1,m}_{i}, \quad r_i=\frac{\sum_{j=1}^{n}\left|\xi_i^{k-1,m}(t_j)-z^{k-1,m}_{i}\right|}{n}, \quad i=1,\ldots,m.
\]
The boundary $C=\cup_{i=1}^{m}C_i$ is parametrized by the function $\xi(t)$.

Since $\omega$ is analytic in the region $G$, thus once we obtain its boundary values from~(\ref{e:om-k}), we can compute the values of $w=\omega(z)$ at interior points $z\in G$ using the Cauchy integral formula,
\[
w=\omega(z)=\frac{1}{2\pi\i}\int_J \frac{\omega(\eta(t))}{\eta(t)-z}\eta'(t)dt.
\]
Since the derivative of $\xi(t)$ is known, we can find the boundary values of the derivative of the function $\omega$ by differentiation both sides of~(\ref{e:om}) with respect to $t$. Hence, we can also use the Cauchy integral formula to find the values of $\omega'(z)$ for interior points $z\in G$.

\subsection{The inverse conformal mapping}

The inverse mapping function $\omega^{-1}$ is analytic in the circular region $\Omega$. Since the boundary $C$ is paramterized by $\xi(t)$, $t\in J$, the values of $z=\omega^{-1}(w)$ at interior points $w\in\Omega$ can be computed using the Cauchy integral formula
\begin{equation}\label{e:cau-i}
z=\omega^{-1}(w)
=\frac{1}{2\pi\i}\int_C \frac{\omega^{-1}(\xi)}{\xi-w}d\xi
=\frac{1}{2\pi\i}\int_J \frac{\omega^{-1}(\xi(t))}{\xi(t)-w}\xi'(t)dt,
\end{equation}
where $\omega^{-1}(\xi(t))=\eta(t)$ and $\xi'(t)=\omega'(\eta(t))\eta'(t)$. 

\section{Koebe's iterative method for unbounded multiply connected regions}
\label{sc:Koe-ud}

Based on the results of \S\ref{sc:sim}, this section presents a numerical implementation of Koebe's iterative method for computing the conformal mapping $w=\omega(z)$ from the unbounded multiply connected region $G$ onto the unbounded circular region $\Omega$ (see Figure~\ref{f:Um}). We present a method for computing the boundary values of the mapping function
\[
\xi(t)=\omega(\eta(t)), \quad t\in J,
\]
the boundary values of its derivative 
\[
\xi'(t)=\eta'(t)\omega'(\eta(t)), \quad t\in J,
\]
and the parameters $z_j,r_j$, $j=1,2,\ldots,m$, of the canonical region $\Omega$. 

The details are similar to the bounded case presented in the previous section.

\subsection{Initializations}

At the beginning, we set
\[
C_i^{0,0}=\Gamma_i \quad{\rm for}\quad i=1,2,\ldots,m.
\]
The curve $C_i^{0,0}$ is parametrized by $\xi_i^{0,0}(t)$ where
\[
\xi_i^{0,0}(t)=\eta_i(t), \quad 
\frac{d}{dt}\xi_i^{0,0}(t)=\eta'_i(t), \quad t\in J_i, \quad i=1,2,\ldots,m.
\]
We Assume also that $z_i^{0,0}$ is a fix point inside $C_i^{0,0}$ for $i=1,\ldots,m$. Thus, initial values of the boundary values of the mapping function are given by
\[
\omega^{0}(\eta(t))= \left\{ \begin{array}{l@{\hspace{0.5cm}}l}
\xi_1^{0,0}(t),&t\in J_1,\\
\hspace{0.3cm}\vdots\\
\xi_m^{0,0}(t),&t\in J_m. \\
\end{array}
\right.
\]

\subsection{Iterations}

For $k=1,2,3,\ldots$, where $k$ denotes the iteration number, we shall repeat the following three steps:

\subsubsection{Step I: The curves}
\label{sscu:int}

For $j=1,2,\ldots,m$, let $\Psi_{k,j}$ be the conformal mapping from the exterior region of the curve $C^{k-1,j-1}_{j}$ onto the exterior unit disk. Then
\[
C^{k-1,j}_{j} = \Psi_{k,j}(C^{k-1,j-1}_{j})
\]
is the unit circle. The curve $C_j^{k-1,j}$ is parametrized by 
\[
\xi_j^{k-1,j}(t)=\Psi_{k,j}(\xi_j^{k-1,j-1}(t)), \quad t\in J_j.
\]

The smooth Jordan curves $C^{k-1,j-1}_{i}$, $i=1,2,\ldots,m$, $i\ne j$, are external to the curve $C^{k-1,j-1}_{j}$. Thus, the function $\Psi_{k,j}$ maps the curves $C^{k-1,j-1}_{i}$ onto smooth Jordan curves 
\[
C^{k-1,j}_{i}= \Psi_{k,j}(C^{k-1,j-1}_{i}),\quad i=1,2,\ldots,m, \quad i\ne j,
\]
external to $C^{k-1,j}_{j}$. The curve $C_i^{k-1,j}$ is parametrized by 
\[
\xi_i^{k-1,j}(t)=\Psi_{k,j}(\xi_i^{k-1,j-1}(t)), \quad t\in J_i,\quad i=1,2,\ldots,m, \quad i\ne j.
\]
The derivative of the function $\xi_i^{k-1,j}(t)$, $i=1,2,\ldots,m$, can be computed from the boundary values of the derivative of the mapping function $\Psi_{k,j}$. 

Finally, we set
\[
z_j^{k-1,j}=0.
\]
For $i=1,2,\ldots,m$, $i\ne j$, the point $z^{k-1,j-1}_{i}$ inside the curve $C^{k-1,j-1}_{i}$ in the exterior region of the curve $C^{k-1,j-1}_{j}$ will be mapped by the function $\Psi_{k,j}$ into a point 
\[
z^{k-1,j}_{i}=\Psi_{k,j}\left(z^{k-1,j-1}_{i}\right), \quad i=1,2,\ldots,m, \quad i\ne j,
\]
inside the curve $C^{k-1,j}_{i}$ in the exterior of the circle $C^{k-1,j}_{j}$.

\subsubsection{Step II: Normalization}
\label{sscu:nor}

After computing $\xi_i^{k-1,m}(t)$, $t\in J_i$, for $i=1,\ldots,m$, then the function $\hat\omega^k$ with the boundary values
\begin{equation}\label{e:hw-k}
\hat\omega^k(\eta(t))= \left\{ \begin{array}{l@{\hspace{0.5cm}}l}
\xi_1^{k-1,m}(t),&t\in J_1,\\
\hspace{0.3cm}\vdots\\
\xi_m^{k-1,m}(t),&t\in J_m,
\end{array}
\right.
\end{equation}
is the conformal mapping from the region $G$ onto the exterior region of the curves $C^{k-1,m}_{i}$, $i=1,\ldots,m$. However, the function $\hat\omega^k$ does not satisfies the normalization~(\ref{e:cond-u}). The function $\hat\omega^k$ has the expansion near $\infty$,
\[
\hat\omega^k(z) = bz+c_0+c_1z^{-1}+c_2z^{-2}+\cdots
\]
with positive real constant $b$. Since $\alpha$ is in the exterior of $G$, then the constants $b$ and $c_0$ can be computed by~\cite[p.~2]{Gak66}
\begin{eqnarray*}
b&=&-\frac{1}{2\pi\i}\int_\Gamma\frac{\hat\omega_k(\eta)}{\eta-\alpha} \frac{d\eta}{\eta-\alpha}, \\
c_0&=&-\frac{1}{2\pi\i}\int_\Gamma[\hat\omega_k(\eta)-a\eta] \frac{d\eta}{\eta-\alpha}.
\end{eqnarray*}

Define a function $\psi^k$ by
\[
\psi^k(z)=\frac{z-c_0}{b}.
\]
Then the function $\omega^k$ defined by
\begin{equation}\label{e:om-ku}
\omega^{k}(\eta(t))=\psi^k\circ\hat\omega^k
\end{equation}
is the conformal mapping from the region $G$ onto the exterior region of the curves 
\[
C^{k-1,m+1}_{i}=\psi^k(C^{k-1,m}_{i}),\quad i=1,\ldots,m, 
\]
and satisfies the normalization~(\ref{e:cond-u}). The curve $C^{k-1,m+1}_{i}$ is parametrized by
\[
\xi_i^{k-1,m+1}(t)=\psi^k(\xi_i^{k-1,m}(t)),\quad i=1,\ldots,m.
\] 
Then, we set
\[
z_i^{k-1,m+1}=\psi^k(z_i^{k-1,m}),\quad i=1,\ldots,m.
\]

\subsubsection{Step III: Update and conditions of convergence}
\label{sscu:con}

The boundary values of the approximate mapping function $w=\omega^k(z)$ are given by~(\ref{e:om-ku}) and the boundary values of the derivative $\frac{d}{dz}\omega^k(z)$ can be computed from~(\ref{e:om-ku}) and~(\ref{e:hw-k}). We stop the iteration if
\begin{equation}\label{e:convu}
\|\omega^{k}-\omega^{k-1}\|_\infty<\varepsilon\quad {\rm or}\quad k>Max,
\end{equation}
where $\varepsilon$ is a given tolerance and $Max$ is the maximum number of iterations allowed.

If the condition~(\ref{e:conv}) is not satisfied, we set	
\[
C_i^{k,0}=C_i^{k-1,m+1}, \quad i=1,\ldots,m.
\]
For $i=1,2,\ldots,m$, the curve $C_i^{k,0}$ is parametrized by
\[
\xi_i^{k,0}(t)=\xi_i^{k-1,m+1}(t).
\] 
Then, we set $k=k+1$ and repeat Steps I--III.

\subsection{The interior values}

If the method converges, then we consider the bounded multiply connected region bounded by the circles $C_1^{k-1,m+1}$, \ldots, $C_m^{k-1,m+1}$, as the canonical region $\Omega$. The boundaries of $\Omega$ are then given by
\[
C_i=C_i^{k-1,m+1}, \quad i=1,\ldots,m.
\] 
The centre $z_i$ and the radius $R_i$ of the circle $C_i$ are approximated by
\[
z_i=z^{k-1,m+1}_{i}, \quad R_i=\frac{\sum_{j=1}^{n}\left|\xi_i^{k-1,m+1}(t_j)-z^{k-1,m+1}_{i}\right|}{n}, \quad i=1,\ldots,m.
\]

We consider also the boundary values of the approximate mapping function $\omega^k$ in~(\ref{e:om-ku}) as an approximation of the boundary values of the mapping function $\omega$, i.e., we have 
\begin{equation}\label{e:omu}
\xi(t)=\omega(\eta(t))=\omega^k(\eta(t)), \quad t\in J.
\end{equation}
The function $\xi(t)$ is the parametrization of the boundary $C$.

Since $\omega$ is analytic in the region $G$ with the normalization~(\ref{e:cond-u}) and $\alpha$ is in the exterior of $G$, the function $\frac{\omega(z)-z}{z-\alpha}$ is analytic in $G$ and its value at $\infty$ equals to $1$. Thus once we obtain its boundary values from~(\ref{e:om-ku}) and~(\ref{e:omu}), we can compute the values of $w=\omega(z)$ at interior points $z\in G$ using the Cauchy integral formula~\cite[p.~2]{Gak66}
\begin{equation}\label{e:cau-u}
\omega(z)=-(z-\alpha)+(z-\alpha)\frac{1}{2\pi\i}\int_J \frac{\omega(\eta(t))}{\eta(t)-\alpha} 
\frac{1}{\eta(t)-z}\eta'(t)dt.
\end{equation}

Since $\omega$ is analytic in the region $G$ with the normalization~(\ref{e:cond-u}), the function $\omega'$ is analytic in $G$ with $\omega'(\infty)=1$. In view of~(\ref{e:omu}), the boundary values of the derivative of the function $\omega$ can be obtained by differentiation both sides of~(\ref{e:om-ku}) and~(\ref{e:hw-k}). Thus, we can compute the values of $\omega'(z)$ at interior points $z\in G$ using the Cauchy integral formula~\cite[p.~2]{Gak66}
\begin{equation}\label{e:cau-u'}
\omega'(z)=-1+\frac{1}{2\pi\i}\int_J \frac{\omega'(\eta(t))}{\eta(t)-z}\eta'(t)dt.
\end{equation}

\subsection{The inverse conformal mapping}

The inverse mapping function $\omega^{-1}$ is analytic in the circular region $\Omega$ with a simple pole at $\infty$. The normalization~(\ref{e:cond-u}) implies that the inverse function $\omega^{-1}$ satisfies~\cite{Wen92}
\begin{equation}\label{e:cond-uiu}
\lim_{w\to\infty}\frac{\omega^{-1}(w)}{w}=1.
\end{equation} 
Let $\hat\alpha$ is a fixed point inside $C_1$ in the exterior of $\Omega$, then function $\frac{\omega^{-1}(w)-w}{w-\hat\alpha}$ is analytic in $\Omega$ and its value at $\infty$ equals to $1$. 
Thus, by the Cauchy integral formula, we can compute the values of $z=\omega^{-1}(w)$ at interior points $w\in\Omega$ using~\cite[p.~2]{Gak66} 
\begin{equation}\label{e:cau-iu}
z=\omega^{-1}(w)
=-(w-\hat\alpha)+(w-\hat\alpha)\frac{1}{2\pi\i}\int_J \frac{\omega^{-1}(\xi(t))}{\xi(t)-\hat\alpha}
\frac{1}{\xi(t)-w}\xi'(t)dt,
\end{equation}
where $\xi'(t)=\omega'(\eta(t))\eta'(t)$ and, in view of~(\ref{e:omu}), 
\[
\omega^{-1}(\xi(t))=\eta(t), \quad t\in J.
\]

\section{Numerical examples}

We consider 8 numerical examples. In the first two examples, we consider examples with known solutions. In Examples~3 and~4, we consider examples from~\cite{Weg01c,Weg05}. A region which boundaries are closer together is given in Example 5. In Examples 6 and 7, we consider regions with high connectivity. Finally, we consider a region with piecewise smooth boundaries in Example~8. We presents the discrete errors only for the Examples~1 and~2 since the exact maps are known. For these two examples, as in~\cite{Ben07}, we compute the maximum numerical errors $E_{\omega,n}$, $E_{z,n}$, and $E_{R,n}$ for the boundary values of the mapping function, the centres, and the radii, respectively. For the remaining examples, the exact maps are unknowns.

For function \verb|zfmm2dpart|, we assume that ${\tt iprec}=5$ which means that the tolerance of the FMM is $0.5\times 10^{-15}$. For the function \verb|gmres|, we choose the parameters ${\tt restart}=10$, ${\tt gmrestol}=0.5\times10^{-14}$, and ${\tt maxit}=10$, which means that the GMRES method is restarted every $10$ inner iterations, the tolerance of the GMRES method is $0.5\times10^{-14}$, and the maximum number of outer iterations of GMRES method is $10$. See~\cite{Nas-siam13,Nas-fast} for more details. For Koebe's iterations, we iterate until 
\[
\|\omega^{k}-\omega^{k-1}\|_\infty<0.5\times10^{-13}.
\]
The method converges after few iterations when the boundaries are well separated. The number of iterations increases if the boundaries are closer together. 

For the direct mapping $w=\omega(z)$, we plot images of horizontal and vertical lines from the $z$-plane. For the inverse mapping $z=\omega^{-1}(w)$, we plot images of radial lines and circles from the $w$-plane. The values of the mapping function $\omega(z)$ for interior points $z\in G$ and the values of the inverse mapping function $\omega^{-1}(w)$ for interior points $w\in\Omega$ are computed using the Cauchy integral formula. A fast and accurate method for computing the Cauchy integral formula for interior points is given in~\cite{Nas-siam13,Nas-fast}. 

\begin{figure}[p]%
\centerline{
\scalebox{0.25}{\includegraphics{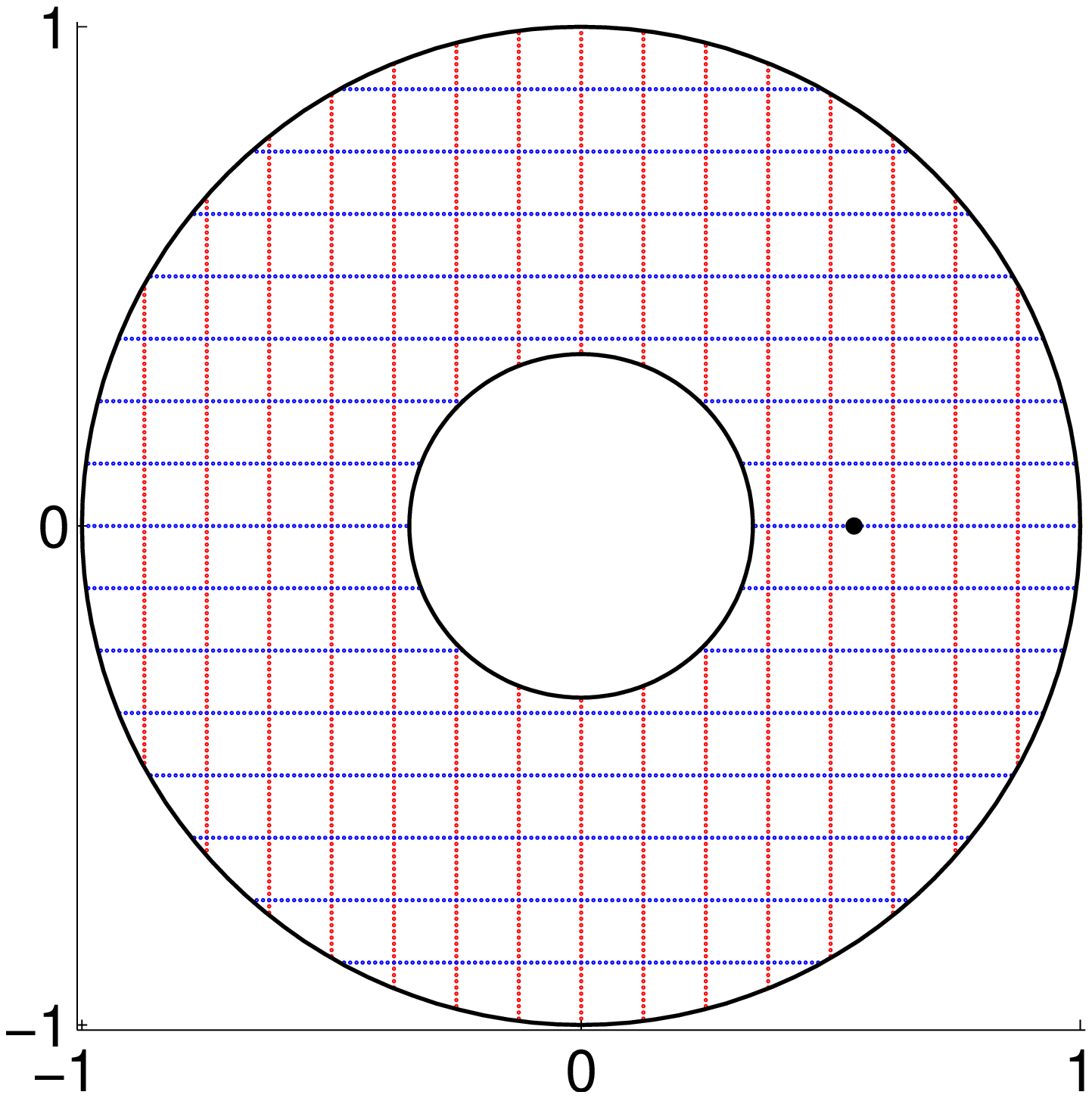}}
\scalebox{0.25}{\includegraphics{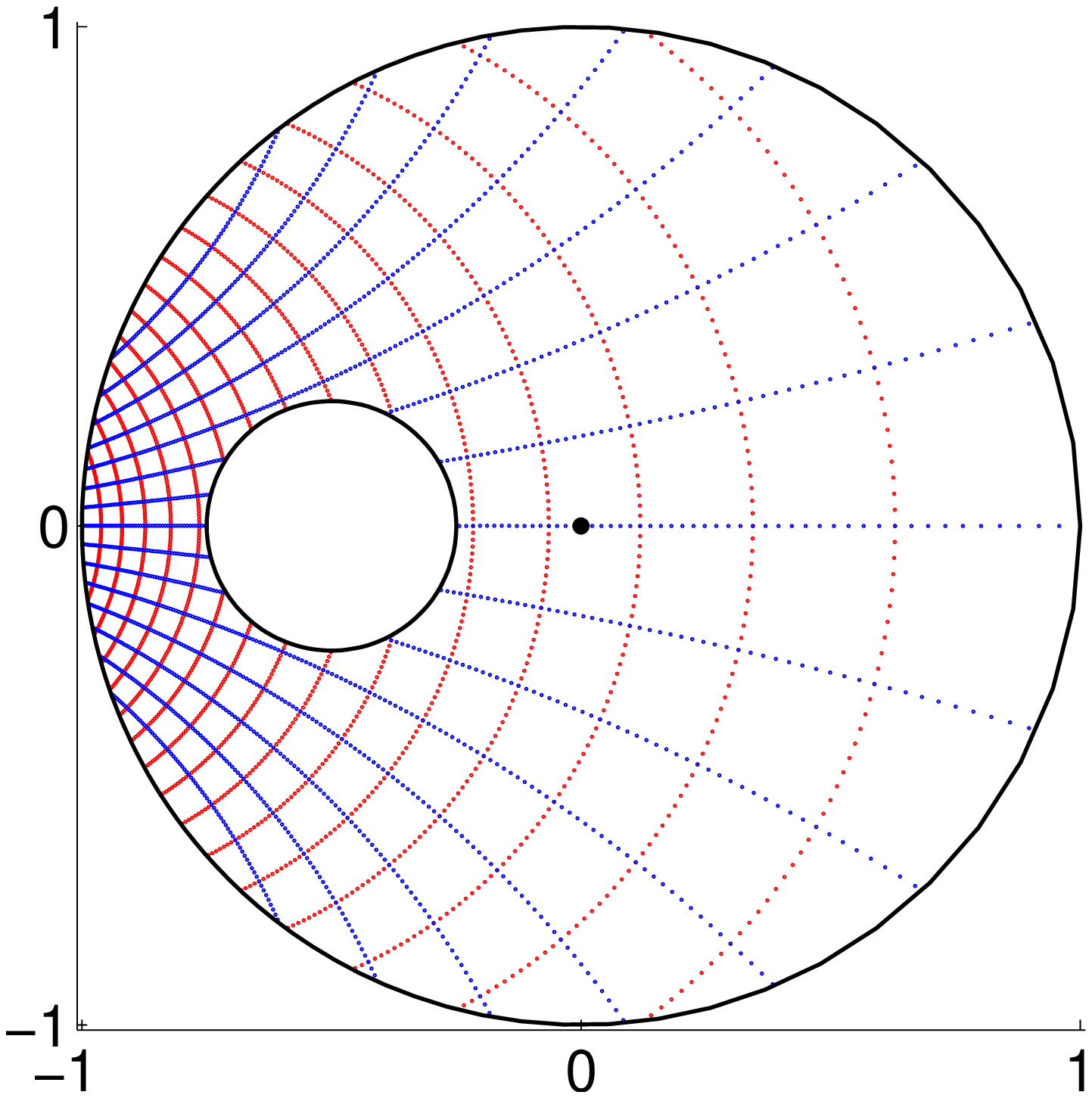}}
}
\caption{\rm The original region $G$ for Example 1 (left) and its image obtained with $n=128$ (right).} 
\label{f:ex1-im}
\end{figure}

\begin{example}\label{ex:1}{\rm
In this example, we consider an example with known exact mapping function (see Figure~\ref{f:ex1-im}). The mapping function (see~\cite[p.~1279]{Gre98})
\begin{equation}\label{e:ex1-om-1}
\omega(z)=\frac{z-a}{1-az}
\end{equation}
maps the region $G$ between two concentric circles with centres $0$, and radii $R$, $1$, in the $z$-plane onto a bounded doubly connected region $\Omega$ between two circles with centres $-0.5$, $0$, and radii $0.25$, $1$, in the $w$-plane where
\[
a=\frac{16}{19+\sqrt{105}},\quad R=\frac{ 8}{13+\sqrt{105}}.
\]
}\end{example}

\begin{figure}[p]%
\centerline{
\scalebox{0.235}{\includegraphics{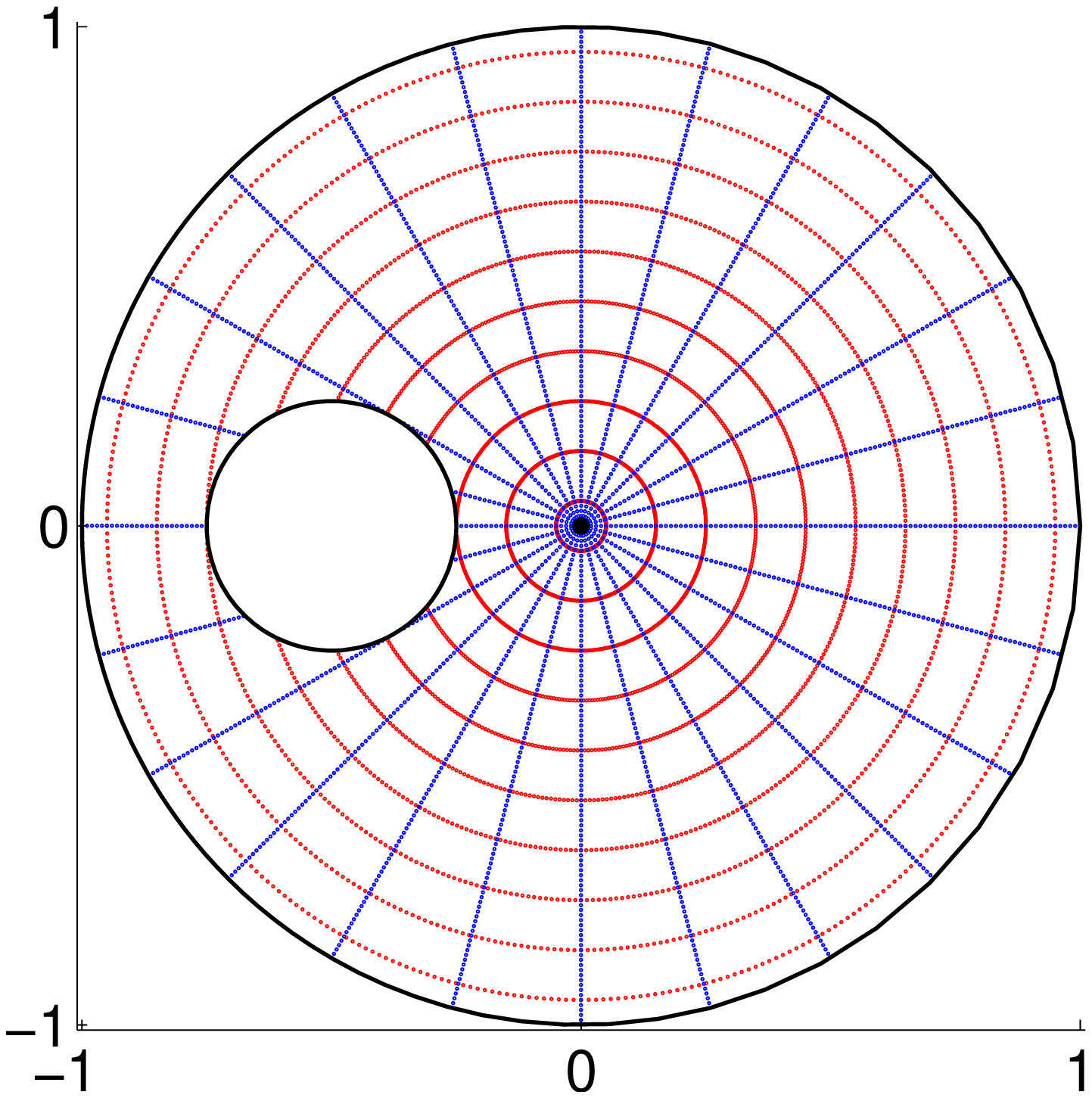}}
\scalebox{0.235}{\includegraphics{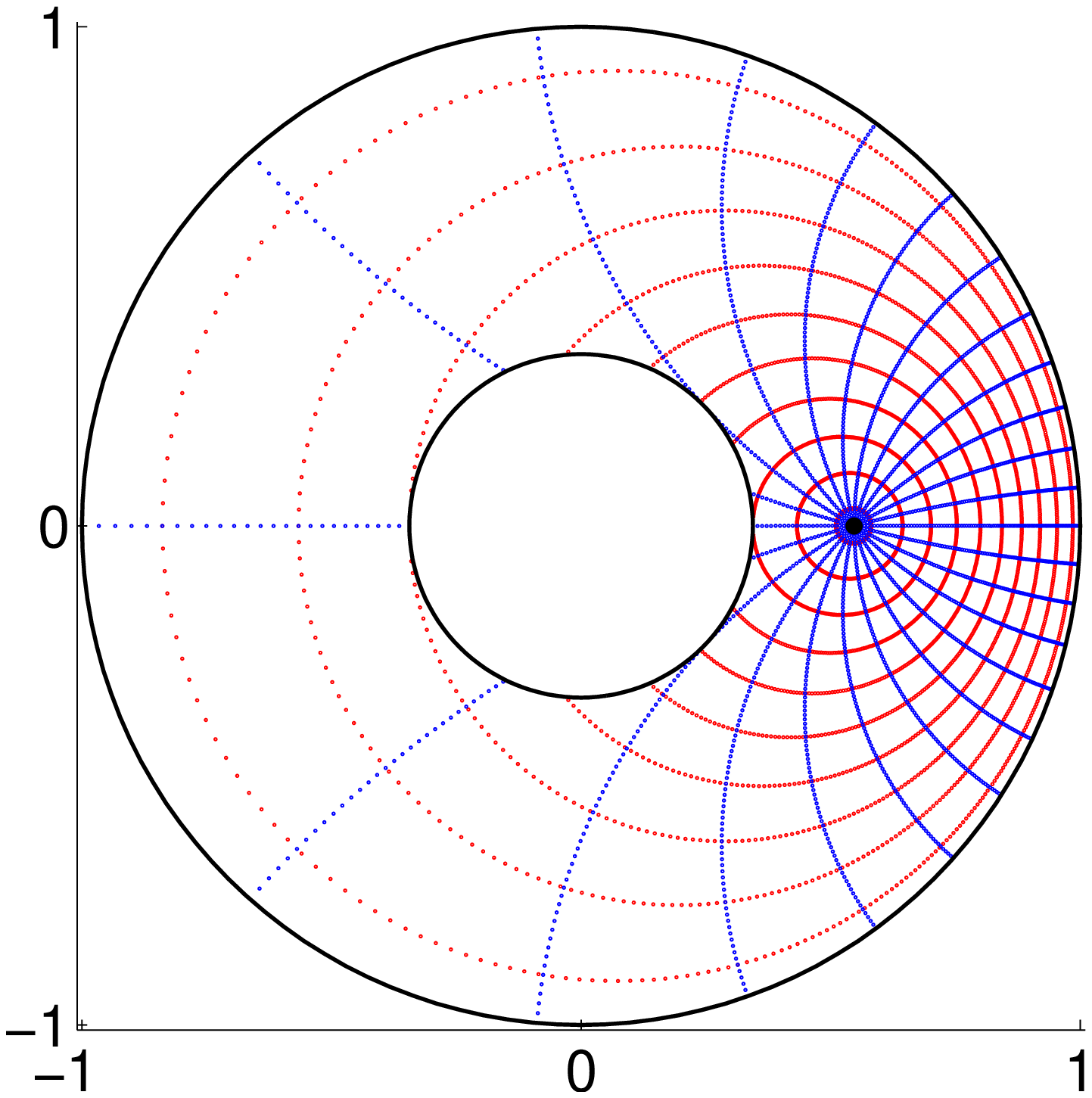}}
}
\caption{\rm The circular region $\Omega$ for Example 1 (left) and its inverse image obtained with $n=128$ (right).} 
\label{f:ex1-inv}
\end{figure}

\begin{figure}[p]%
\centerline{
\scalebox{0.235}{\includegraphics{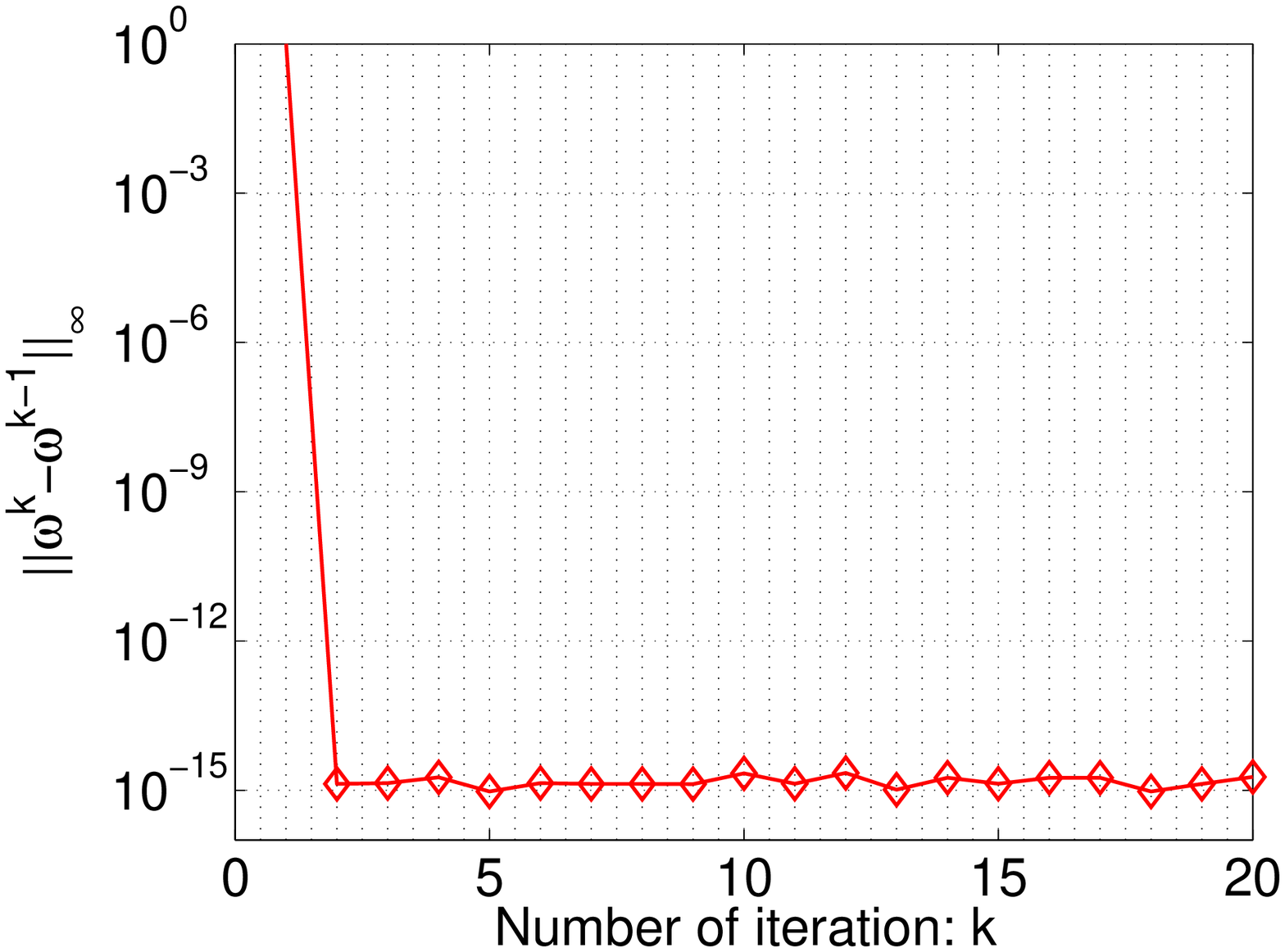}}
\scalebox{0.235}{\includegraphics{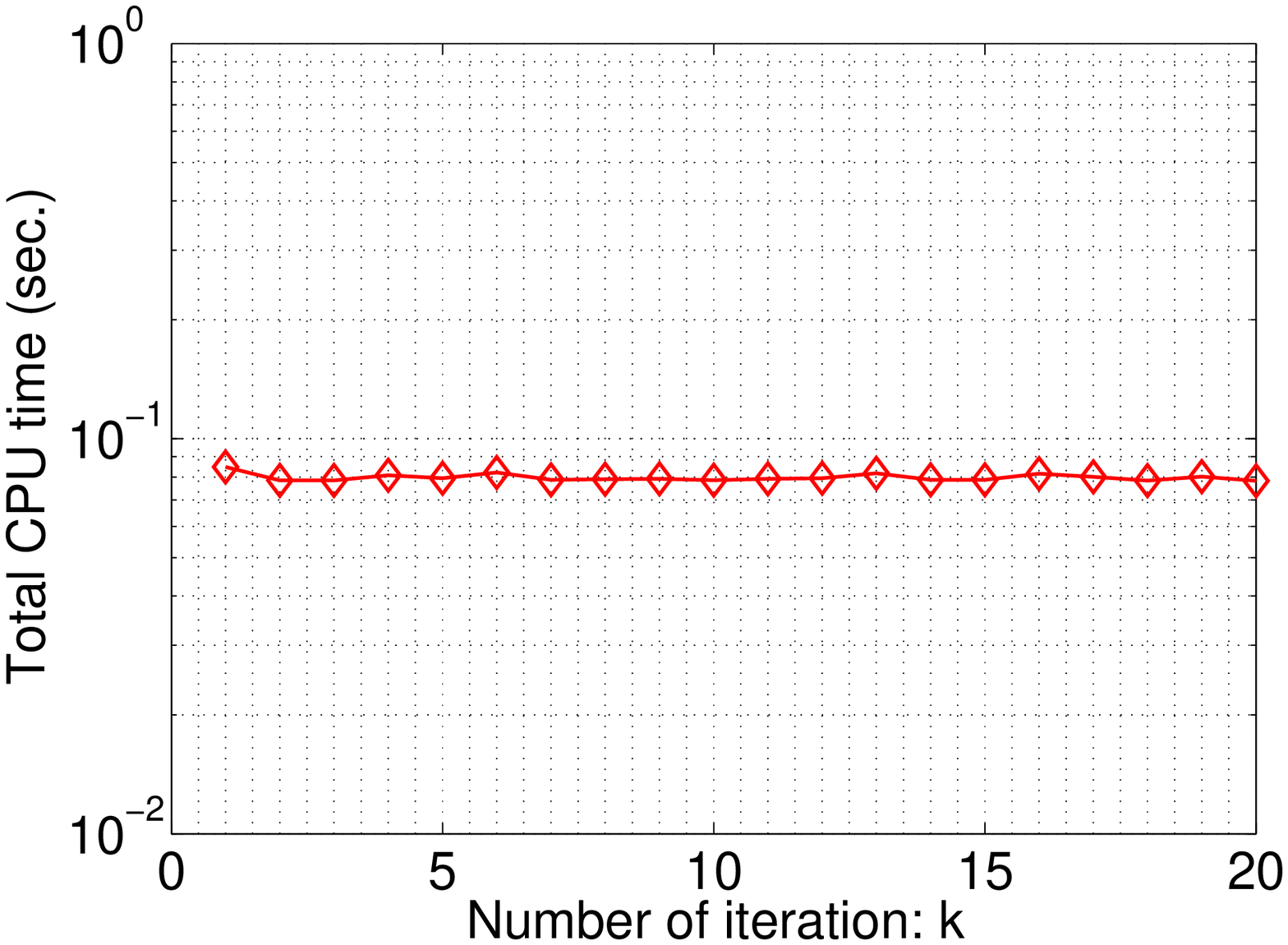}}
\scalebox{0.225}{\includegraphics{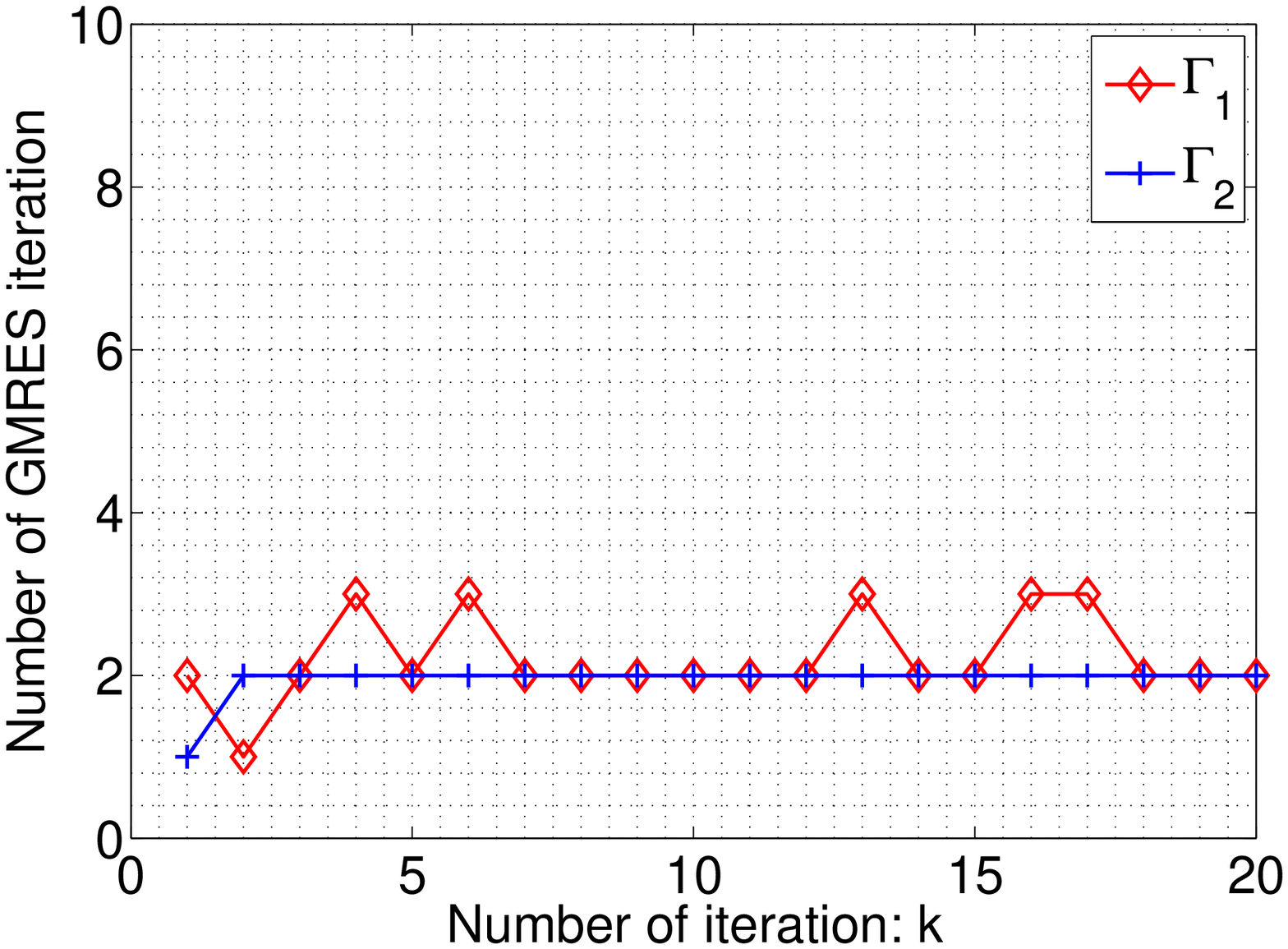}}
}
\centerline{
\scalebox{0.225}{\includegraphics{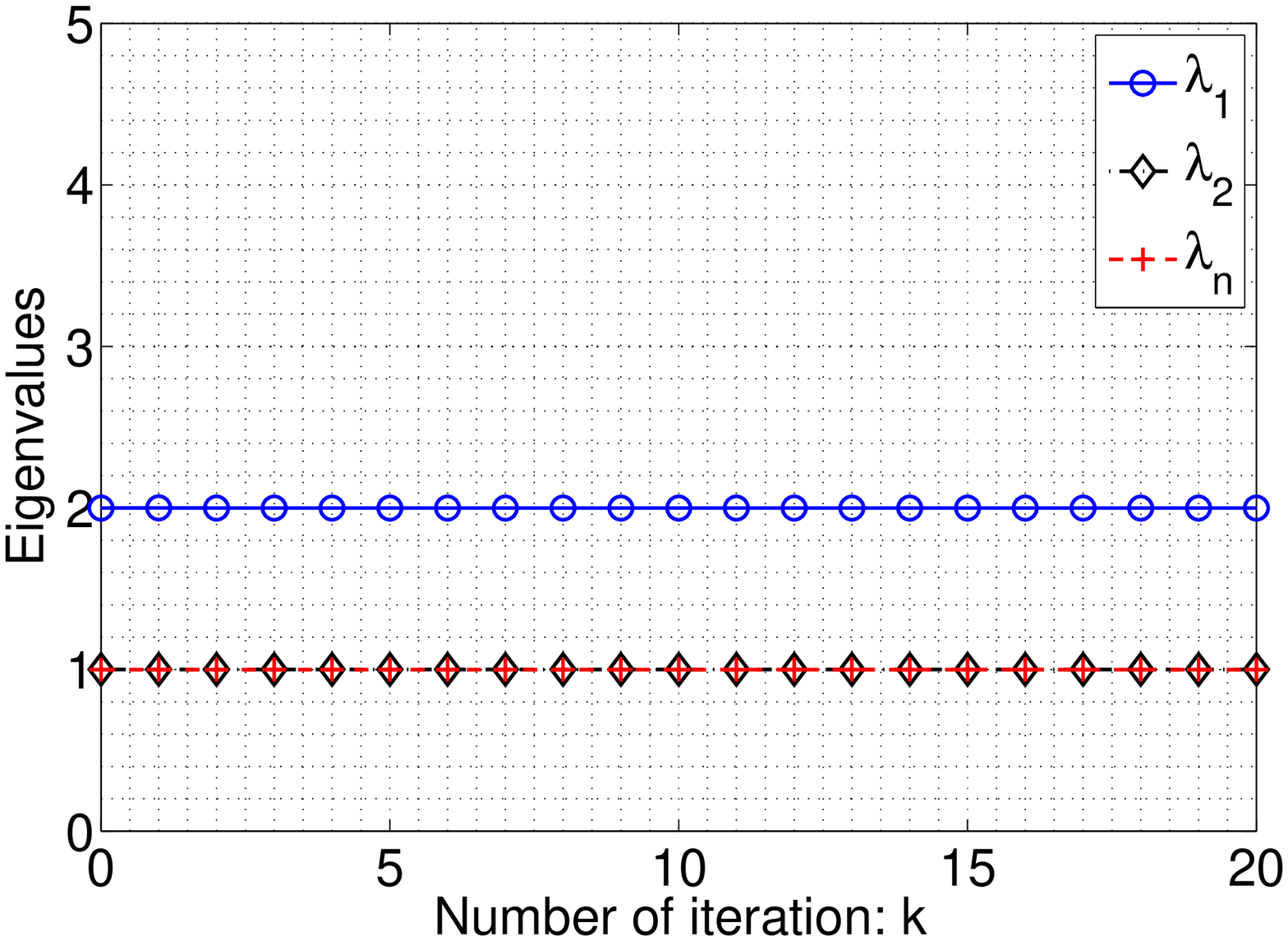}}
\scalebox{0.225}{\includegraphics{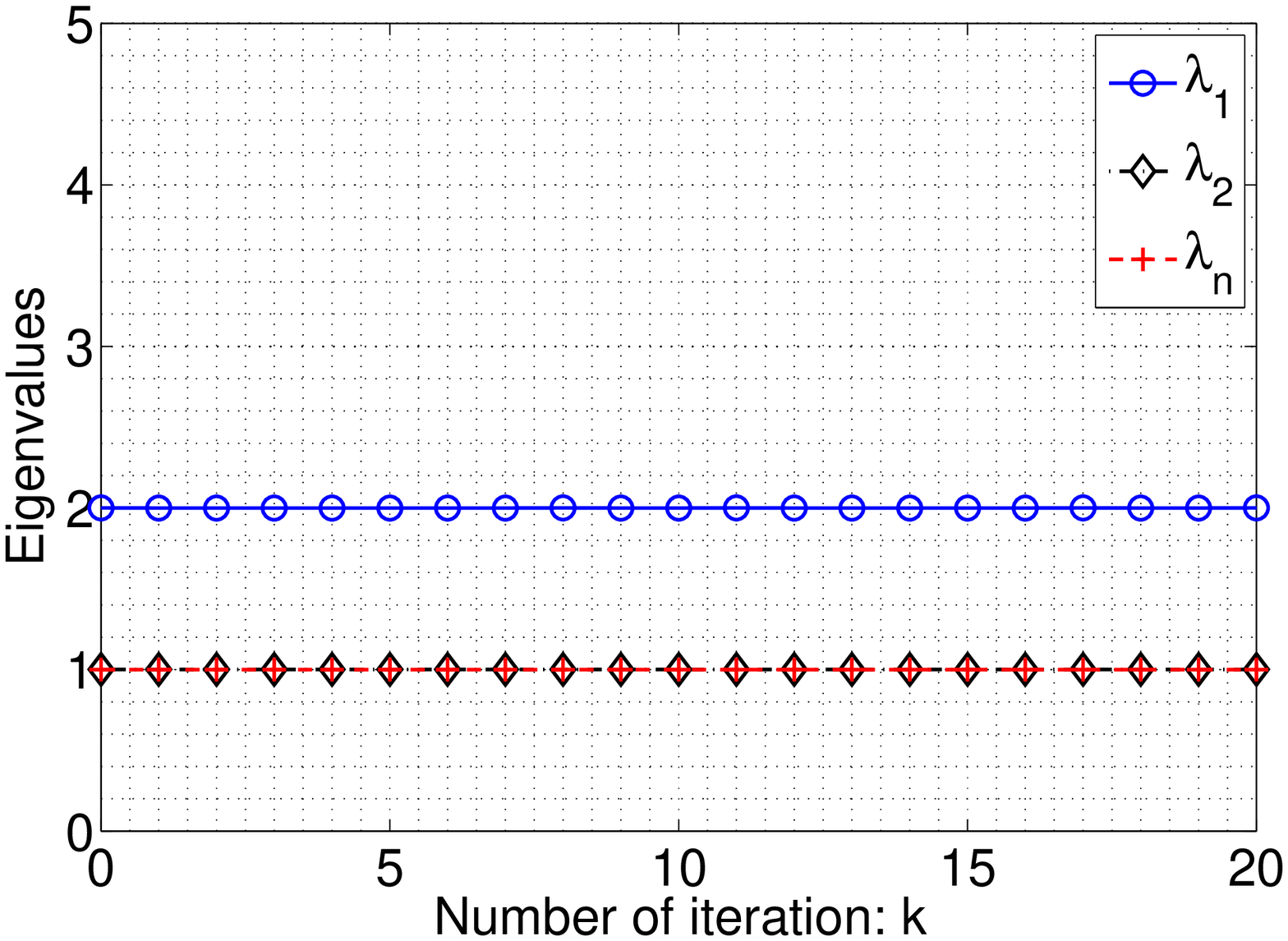}}
\scalebox{0.225}{\includegraphics{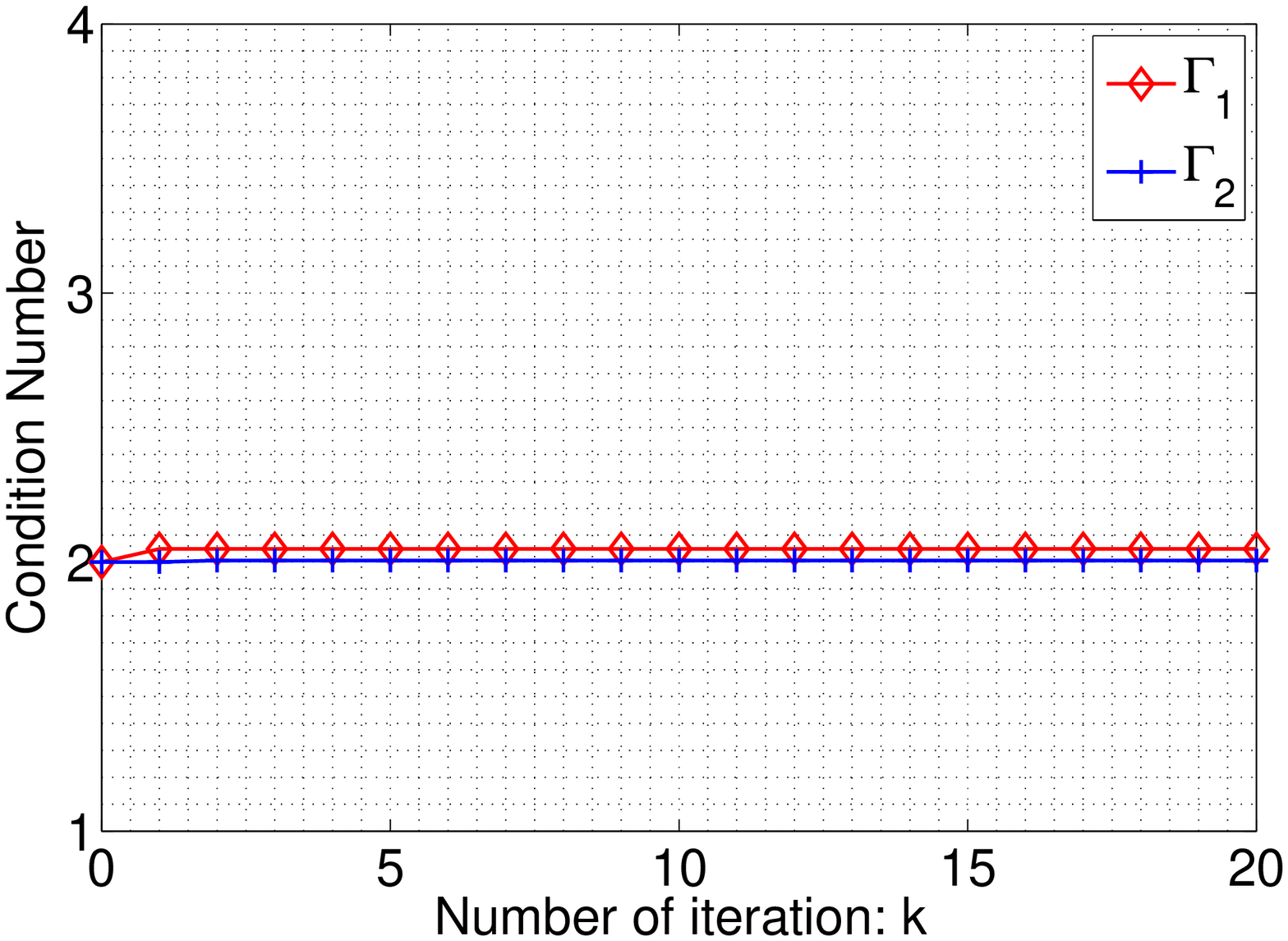}}
}
\caption{\rm Numerical results for Example 1 obtained with $n=128$. 
First row: the successive iteration errors (left), the CPU time in seconds (middle), the number of GMRES iterations (right).
Second row: the eigenvalues $\lambda_1,\lambda_2,\lambda_n$ for the coefficient matrix of the linear system on $\Gamma_1$ (left)  and $\Gamma_2$ (middle), the condition number of the coefficient matrices (right).} 
\label{f:ex1-err}
\end{figure}

Table~\ref{t:ex1} shows the error $E_{\omega,n}$, $E_{Z,n}$, and $E_{R,n}$ for various values of $n$. The original region $G$ and its image are shown in Figure~\ref{f:ex1-im}. In Figure~\ref{f:ex1-inv}, we show the canonical region $\Omega$ and its inverse image. The successive error $\|\omega^{k}-\omega^{k-1}\|_\infty$ vs. the number of iteration $k$ for $k=1,2,\ldots,20$ is shown in Figure~\ref{f:ex1-err}. The successive error is less than $10^{-14}$ after only $2$ iterations. For the first $20$ iterations, the CPU time (seconds) and the number of GMRES iterations for each of the boundaries $\Gamma_1$ and $\Gamma_2$ vs. the number of iteration $k$ are shown in Figure~\ref{f:ex1-err}. The largest eigenvalue $\lambda_1$, the second largest eigenvalue $\lambda_2$, and the smallest eigenvalue $\lambda_n$ for each of the boundaries $\Gamma_1$ and $\Gamma_2$ are shown in Figure~\ref{f:ex1-err}. The other $n-3$ eigenvalues are in the interval $[\lambda_n,\lambda_2]$. When $k$ increases, $\lambda_1\approx2$ and $\lambda_n\approx\lambda_2\approx1$. Thus, the other eigenvalues are also approximately equal to $1$. Figure~\ref{f:ex1-err} shows also the condition number of the coefficient matrices of the linear systems for each of the boundaries $\Gamma_1$ and $\Gamma_2$. The condition number becomes constant when $k$ increases. 

\begin{table}[ht]
\caption{Discretization errors for Example 1.}
\label{t:ex1}%
\vskip-0.75cm
\[
\begin{array}{l@{\hspace{1.5cm}}c@{\hspace{1.0cm}}c@{\hspace{1.0cm}}c} \hline %
n&E_{\omega,n} &E_{Z,n}  &E_{R,n} \\  %
\hline %
16  &1.5(-02) &2.8(-03) &1.9(-03)   \\
32  &7.5(-06) &3.4(-07) &2.3(-07)   \\
64  &1.9(-10) &7.4(-16) &4.2(-16)   \\
128 &3.3(-15) &3.6(-16) &5.3(-16)   \\
256 &3.7(-15) &5.0(-16) &1.4(-16)   \\
\hline %
\end{array}
\]
\end{table}
       
\begin{figure}%
\centerline{
\scalebox{0.235}{\includegraphics{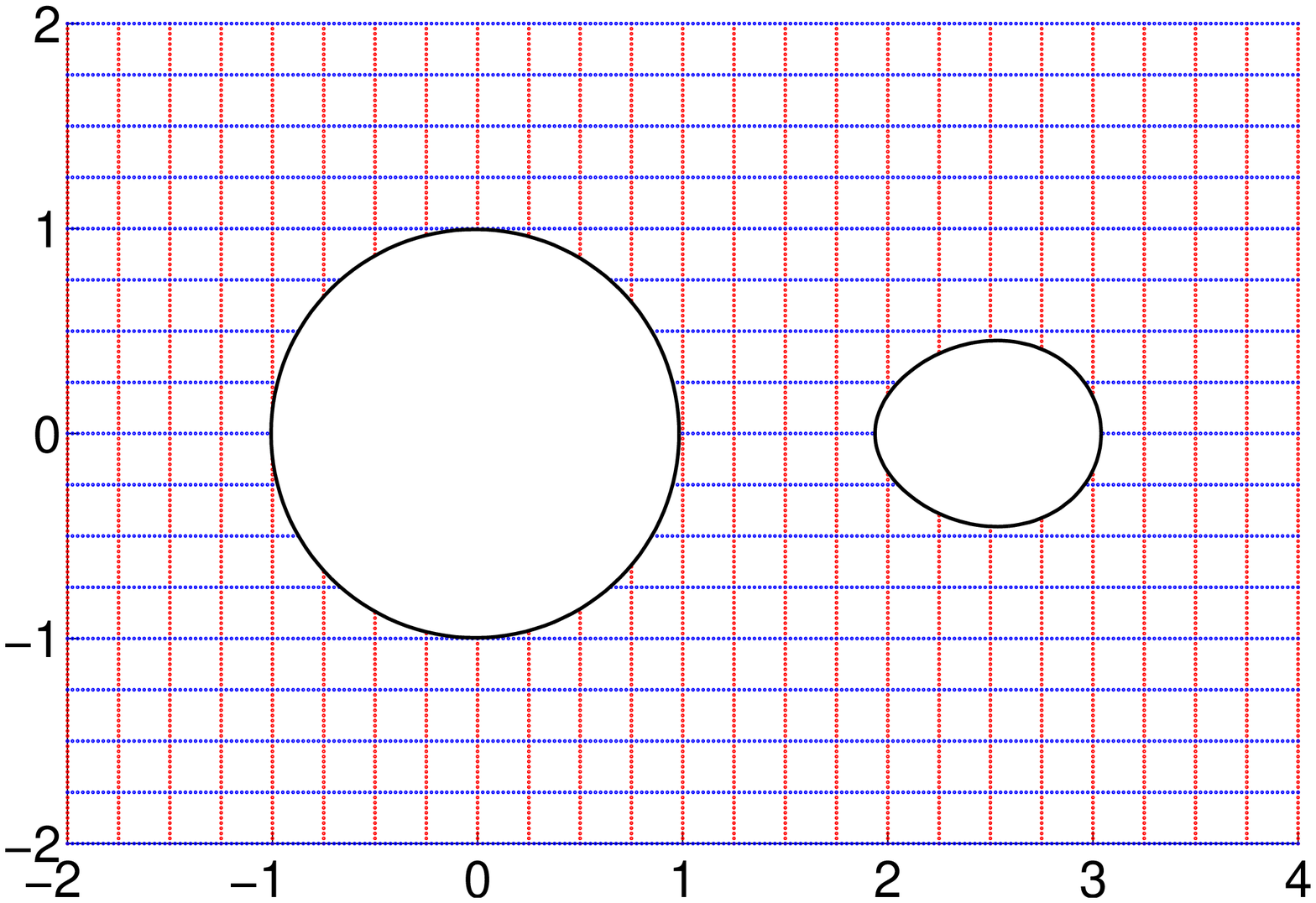}}
\scalebox{0.235}{\includegraphics{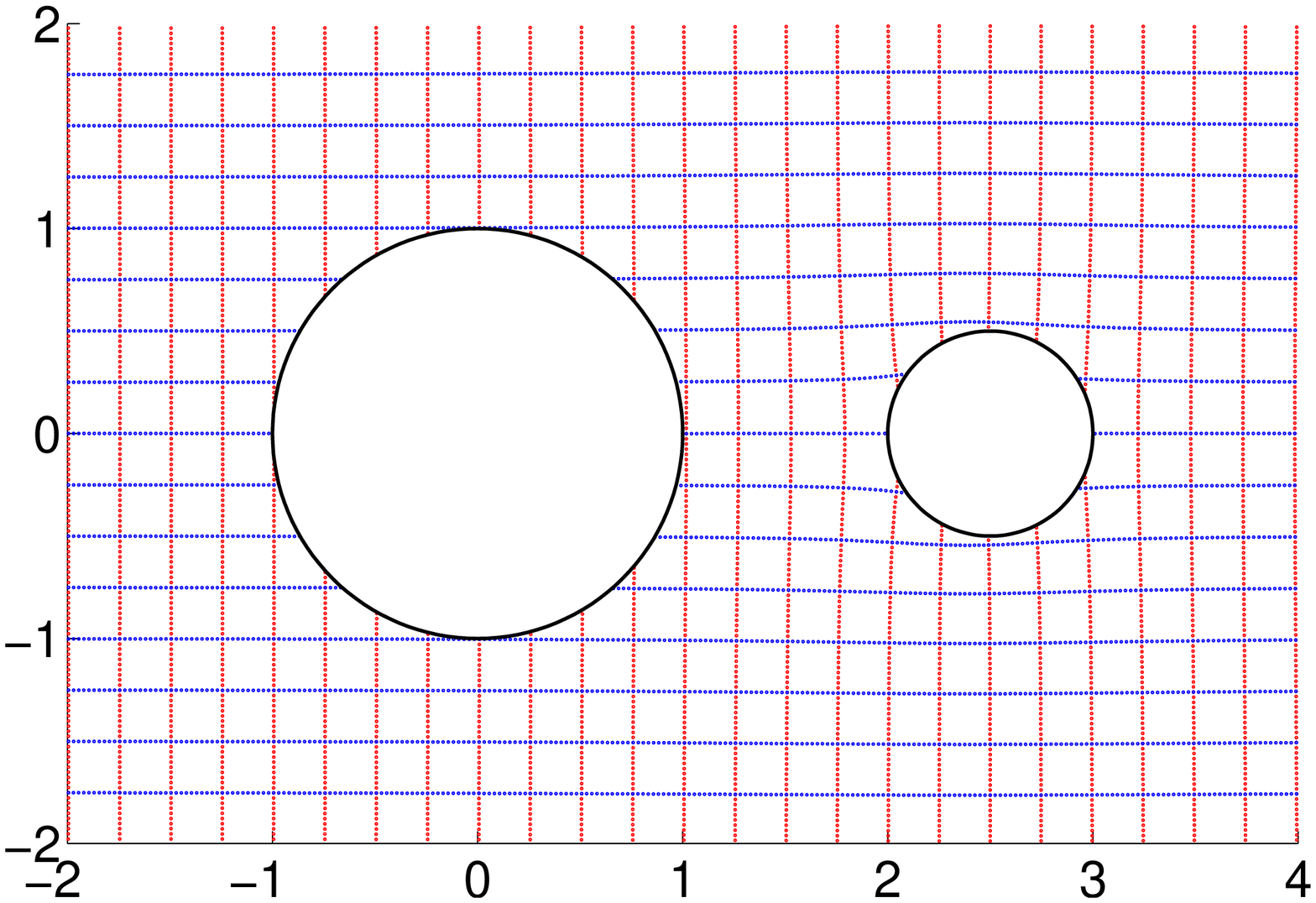}}
}
\caption{\rm The original region $G$ for Example 2 (left) and its image obtained with $n=128$ (right).} 
\label{f:ex2-im}
\end{figure}

\begin{figure}%
\centerline{
\scalebox{0.235}{\includegraphics{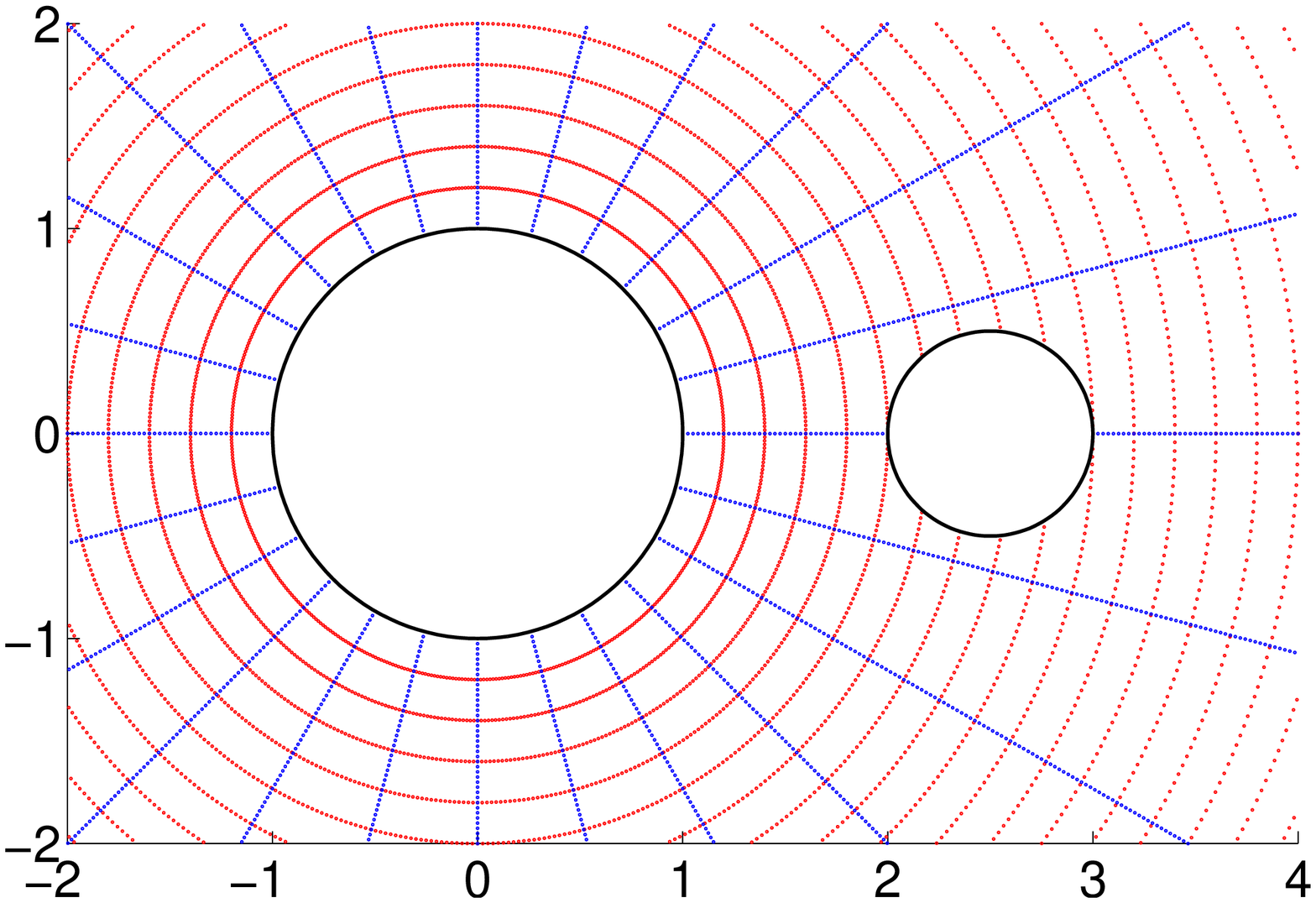}}
\scalebox{0.235}{\includegraphics{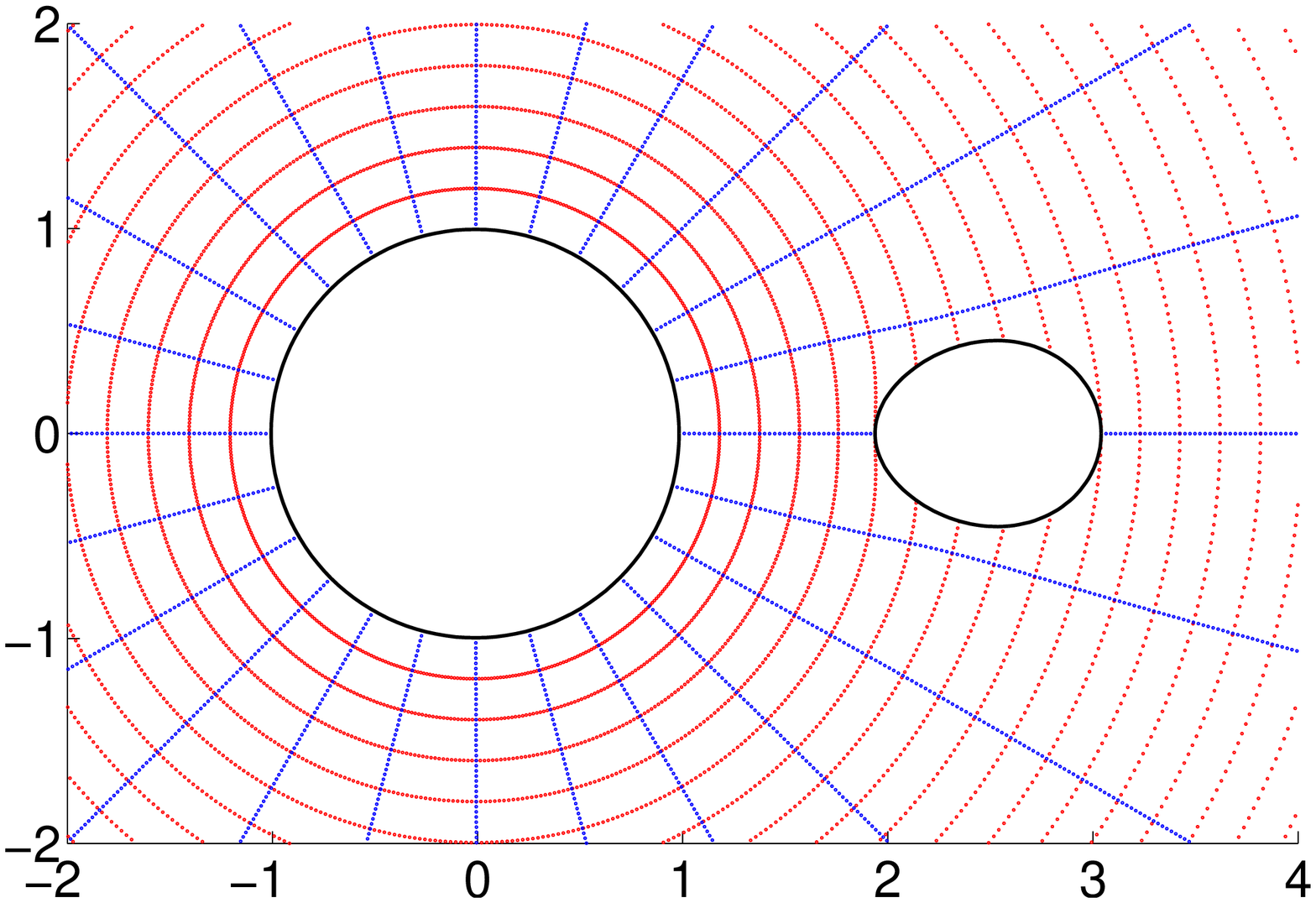}}
}
\caption{\rm The circular region $\Omega$ for Example 2 (left) and its inverse image obtained with $n=128$ (right).} 
\label{f:ex2-inv}
\end{figure}

\begin{figure}%
\centerline{
\scalebox{0.235}{\includegraphics{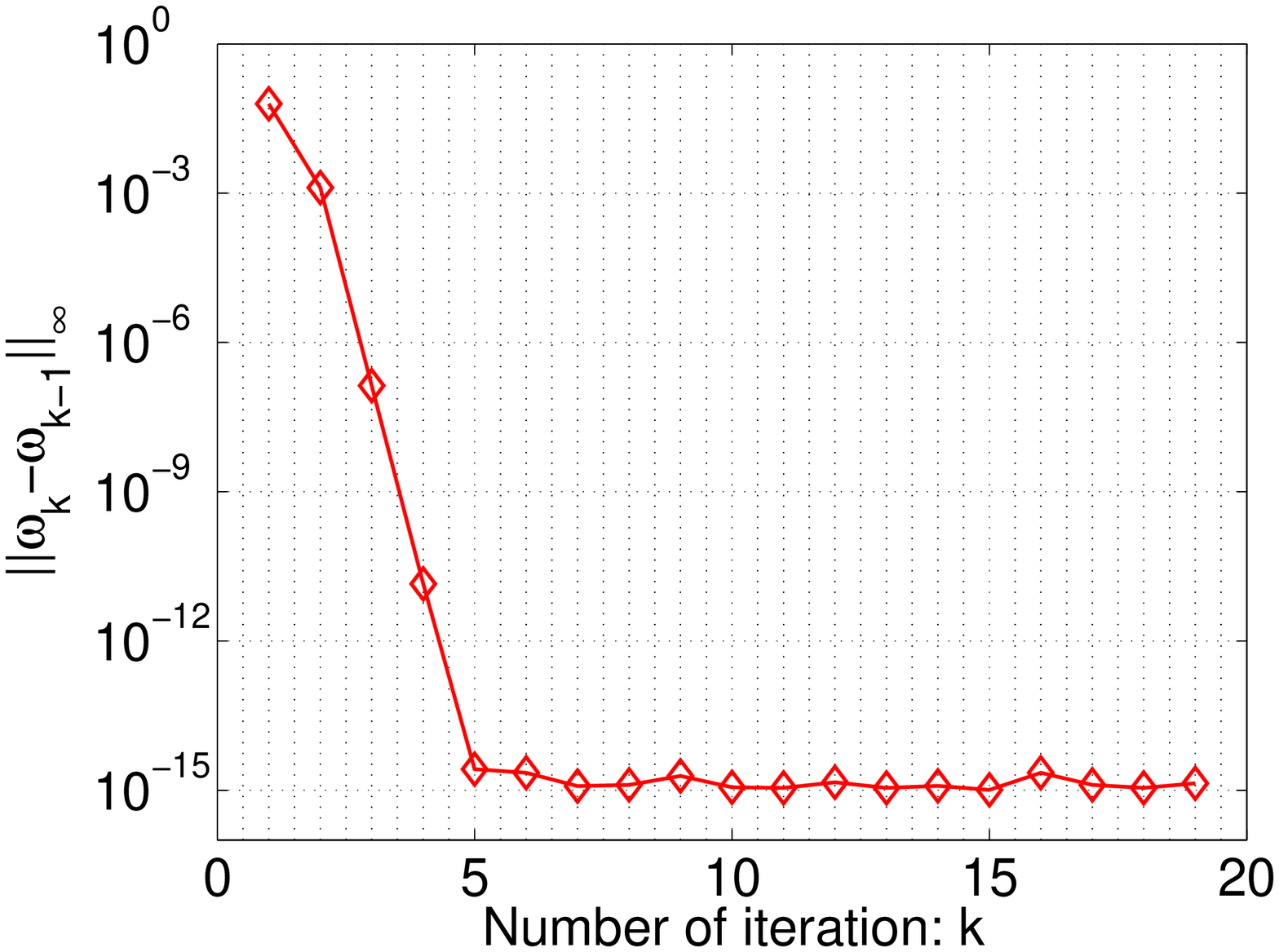}}
\scalebox{0.235}{\includegraphics{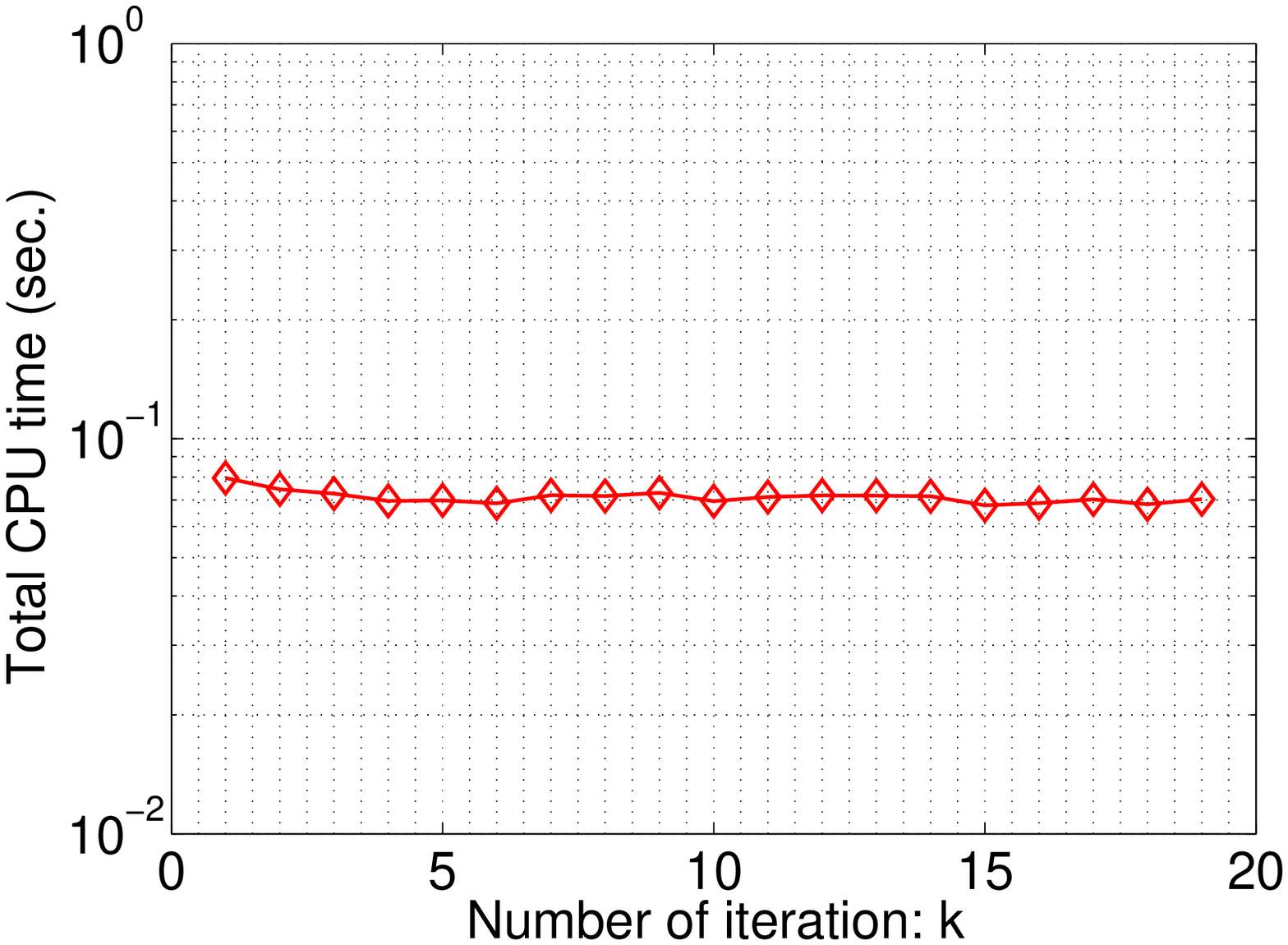}}
\scalebox{0.225}{\includegraphics{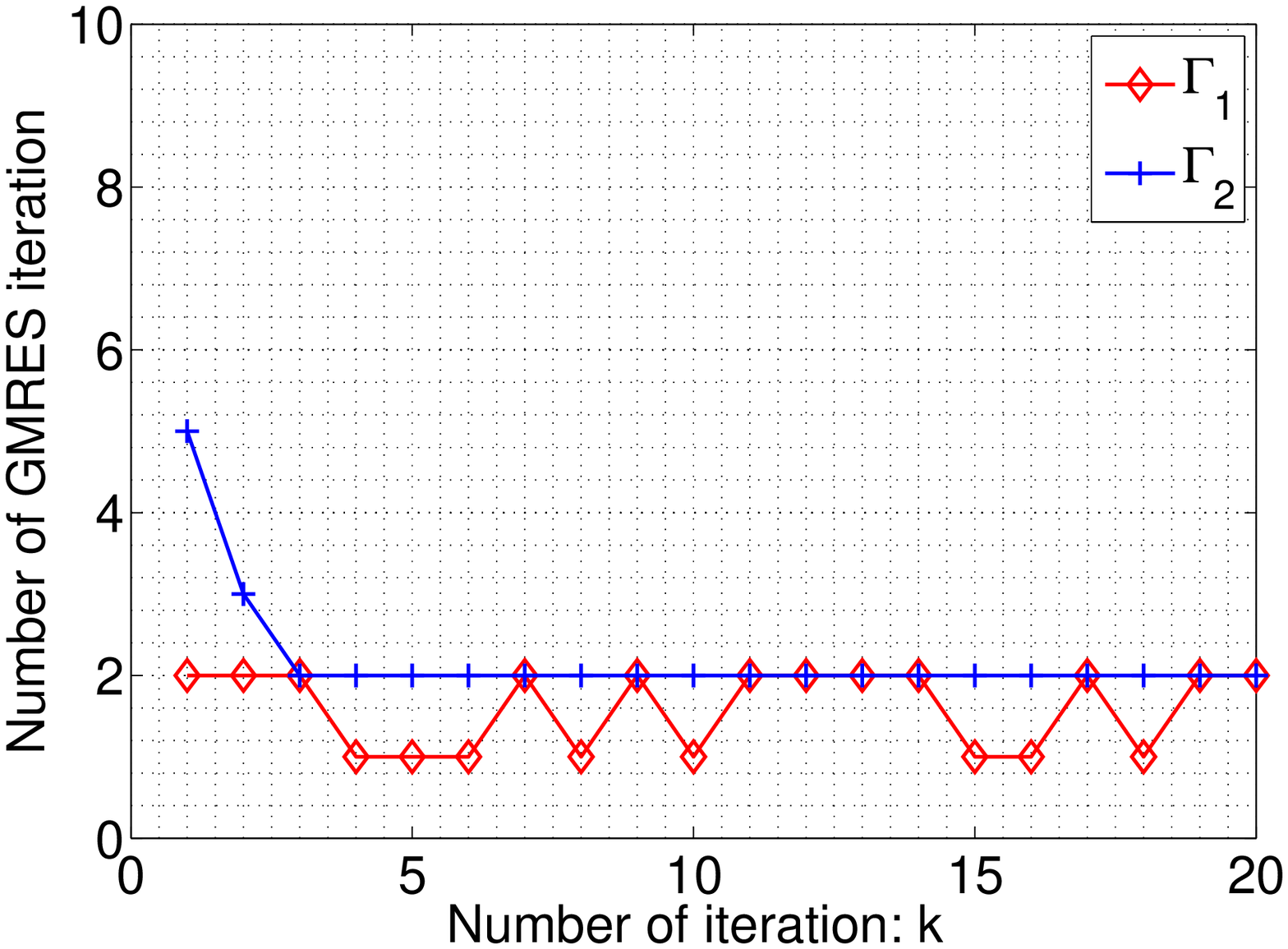}}
}
\centerline{
\scalebox{0.225}{\includegraphics{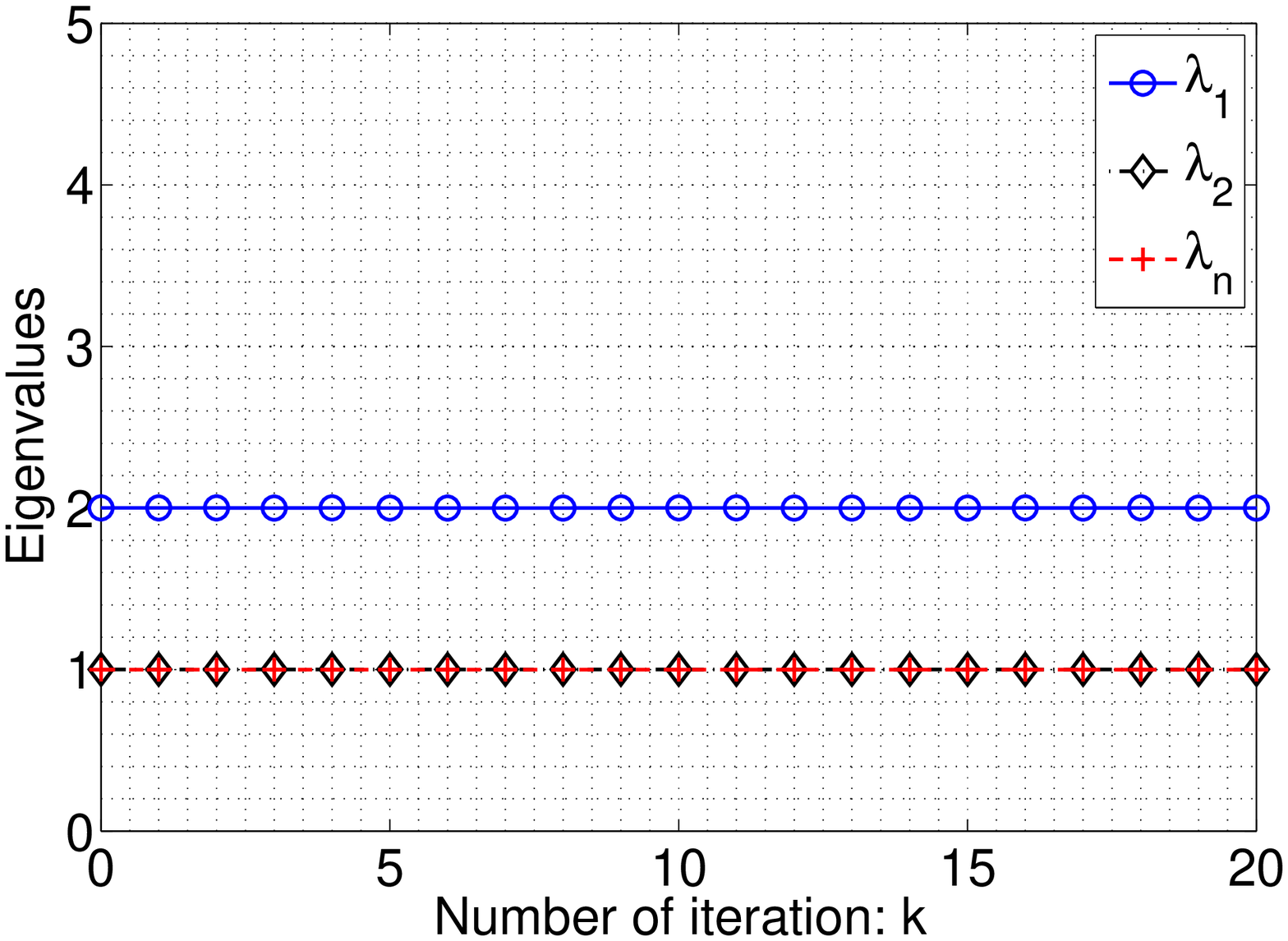}}
\scalebox{0.225}{\includegraphics{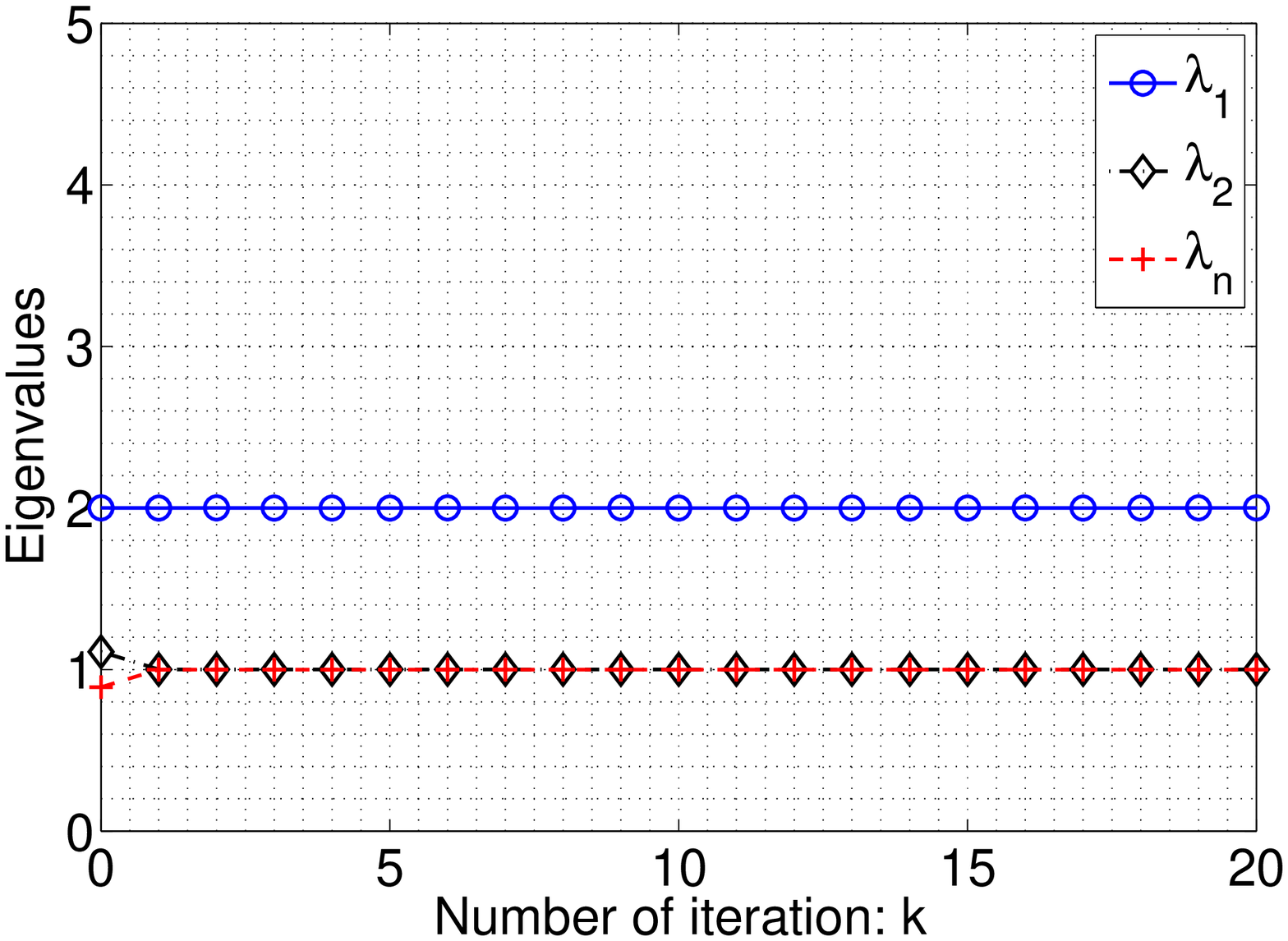}}
\scalebox{0.225}{\includegraphics{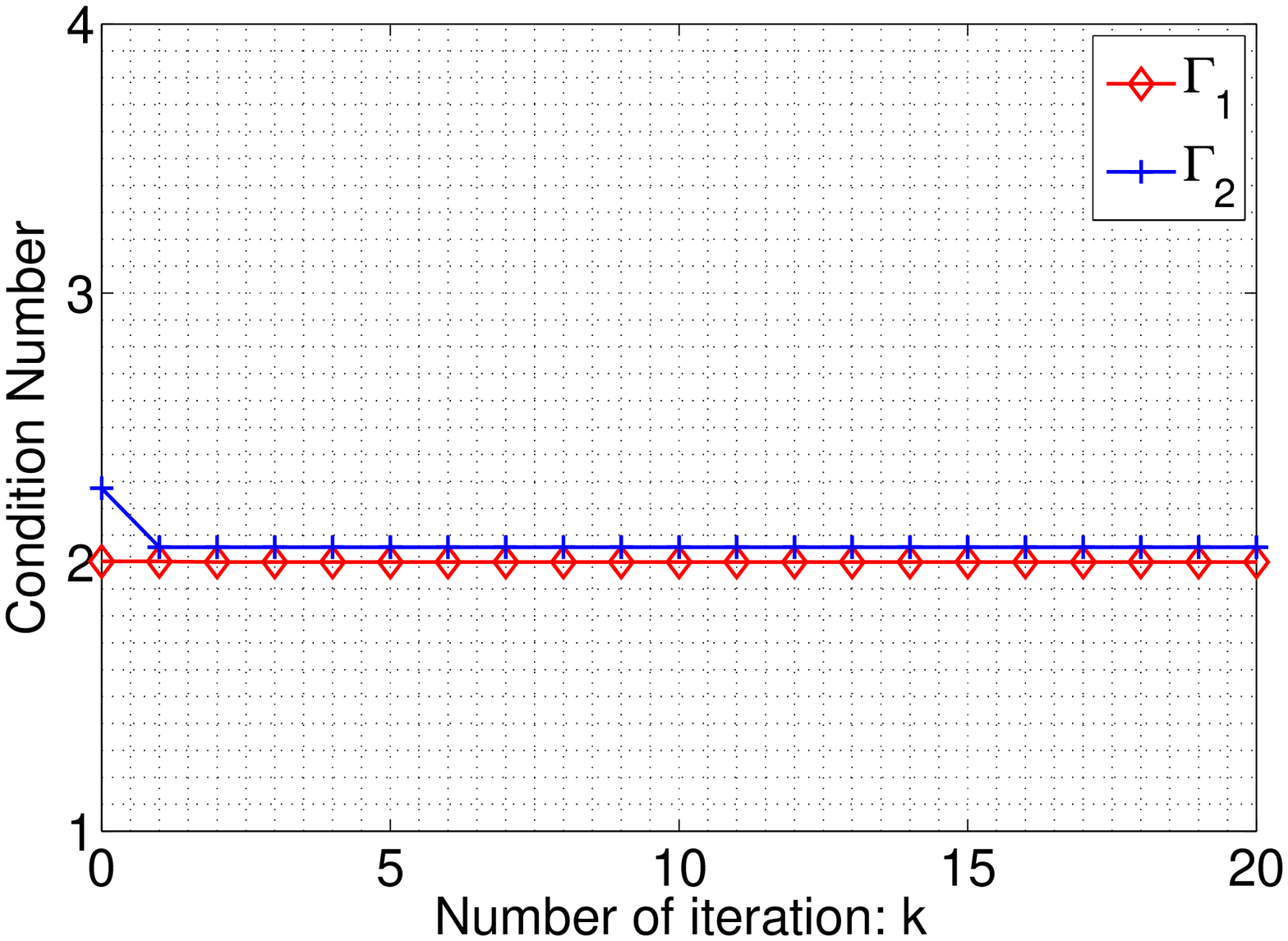}}
}
\caption{\rm The same as in the Figure~\ref{f:ex1-err}, but for Example 2.} 
\label{f:ex2-err}
\end{figure}

\begin{example}\label{ex:2}{\rm
In this example, we consider an example with known exact mapping function from~\cite{Ben07}. The inverse mapping function $z=\omega^{-1}(w)$ maps the region $\Omega$ exterior to two circles with centres $0$, $2.5$ and radii $1$, $0.5$ in the $w$-plane onto an unbounded doubly connected region $G$ exterior to curves $\Gamma_1$, $\Gamma_2$ in the $z$-plane is given by
\begin{equation}\label{e:ex2-om-1}
\omega^{-1}(w)=f_5(f_4(f_3(f_2(f_1(w))))),
\end{equation}
where
\[
f_1(w)=\frac{w-a}{aw-1}, \quad
f_2(w)=\beta w, \quad
f_3(w)=w+\frac{1}{w}, \quad
f_4(w)=\frac{\beta}{\beta^2+1}w,
\]
and
\[
a=\frac{7+2\sqrt{6}}{5}, \quad \beta=30.
\]
The function $f_5$ is given by
\[
f_5(w)=\frac{C_1w+C_2}{w+C_3}, 
\]
where
\[
C_1=\frac{a^4-\beta^2}{a(a^2-\beta^2)}, \quad
C_2=-\frac{3a^2\beta^2+a^2-\beta^4-3\beta^2}{a^2\beta^2+a^2-\beta^2-\beta^4}, \quad
C_3=\frac{a^2+\beta^2}{a(\beta^2+1)}.
\]
The exact mapping function $w=\omega(z)$, which is the inverse of the function $\omega^{-1}$ in~(\ref{e:ex2-om-1}), satisfies the normalization~(\ref{e:cond-u}). The numerical results are shown in Table~\ref{t:ex2} and Figures~\ref{f:ex2-im}--\ref{f:ex2-err}.
}\end{example}

Numerical computing of the inverse mapping function $\omega^{-1}$ for this example using Wegmann's and Fornberg's methods has been given~\cite{Ben07}. It is clear from Tabel~\ref{t:ex2} and from~\cite[Tables~1--4]{Ben07}, the accuracy of our method is almost the same as the accuracy of Wegmann's and Fornberg's methods although our method is used for computing $\omega$ and Wegmann's and Fornberg's methods are used for computing $\omega^{-1}$.

\begin{table}[ht]
\caption{Discretization errors for Example 2.}
\label{t:ex2}%
\vskip-0.5cm
\[
\begin{array}{l@{\hspace{1.5cm}}c@{\hspace{1.0cm}}c@{\hspace{1.0cm}}c} \hline %
n&E_{\omega,n} &E_{Z,n}  &E_{R,n} \\  %
\hline %
16  &5.5(-07) &3.6(-07) &5.7(-09)   \\
32  &3.0(-11) &8.0(-13) &4.7(-14)   \\
64  &1.2(-13) &7.1(-14) &5.0(-14)   \\
128 &1.2(-13) &7.3(-14) &4.6(-14)   \\
256 &1.2(-13) &6.7(-16) &5.2(-16)   \\
\hline %
\end{array}
\]
\end{table}

\begin{figure}%
\centerline{
\scalebox{0.235}{\includegraphics{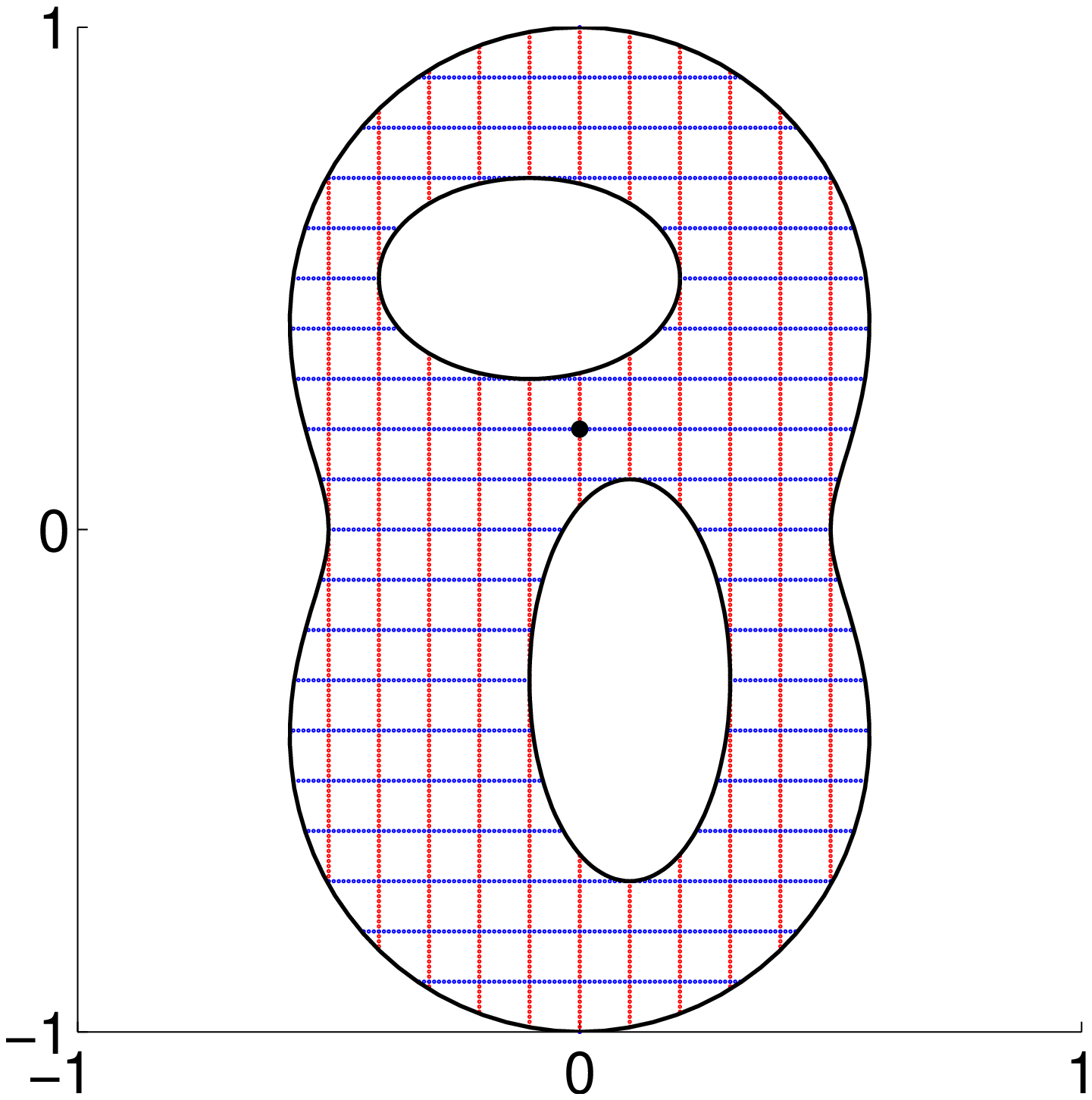}}
\scalebox{0.235}{\includegraphics{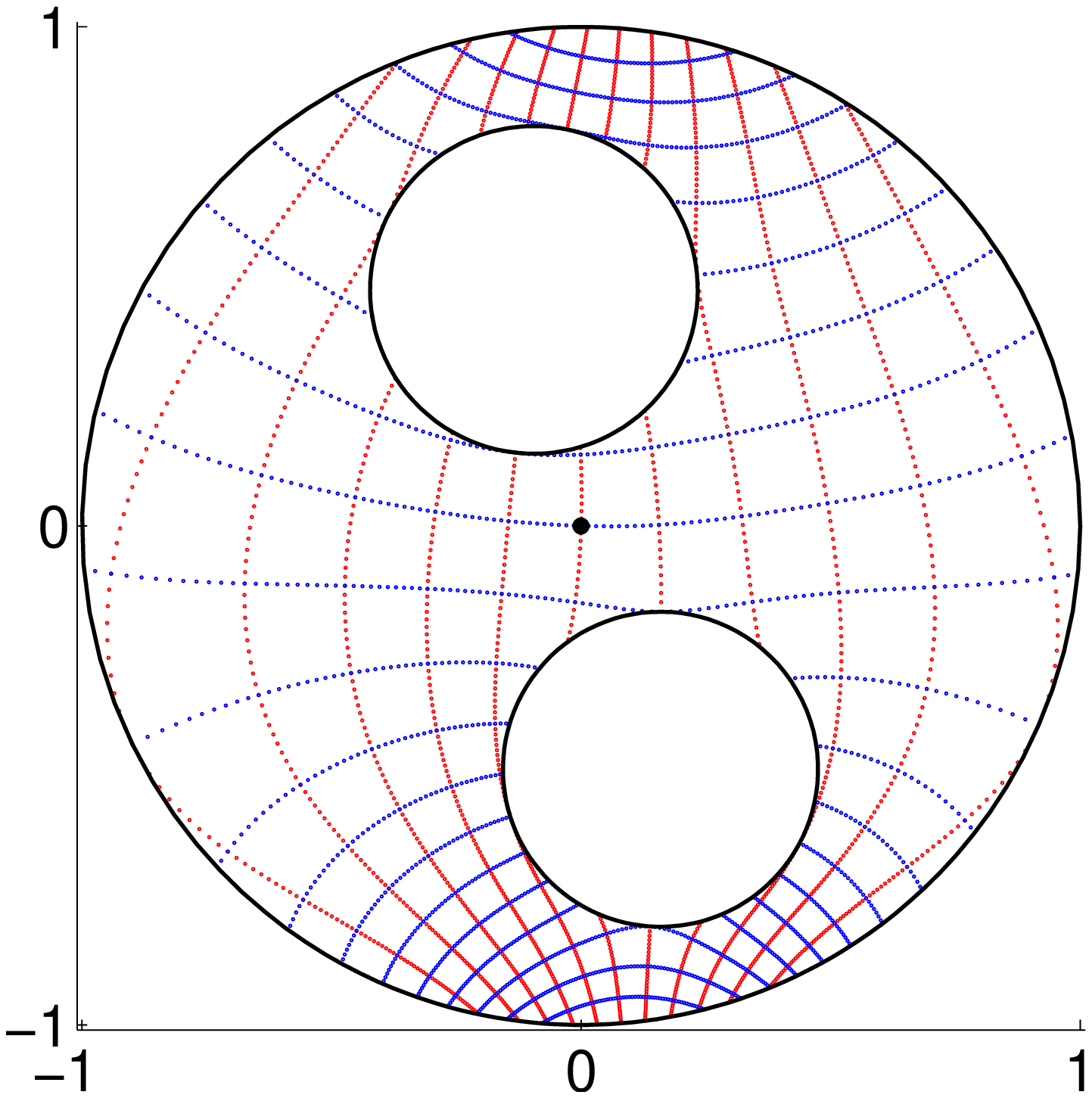}}
}
\caption{\rm The original region $G$ for Example 3 (left) and its image obtained with $n=128$ (right).} 
\label{f:ex3-im}
\end{figure}

\begin{figure}%
\centerline{
\scalebox{0.235}{\includegraphics{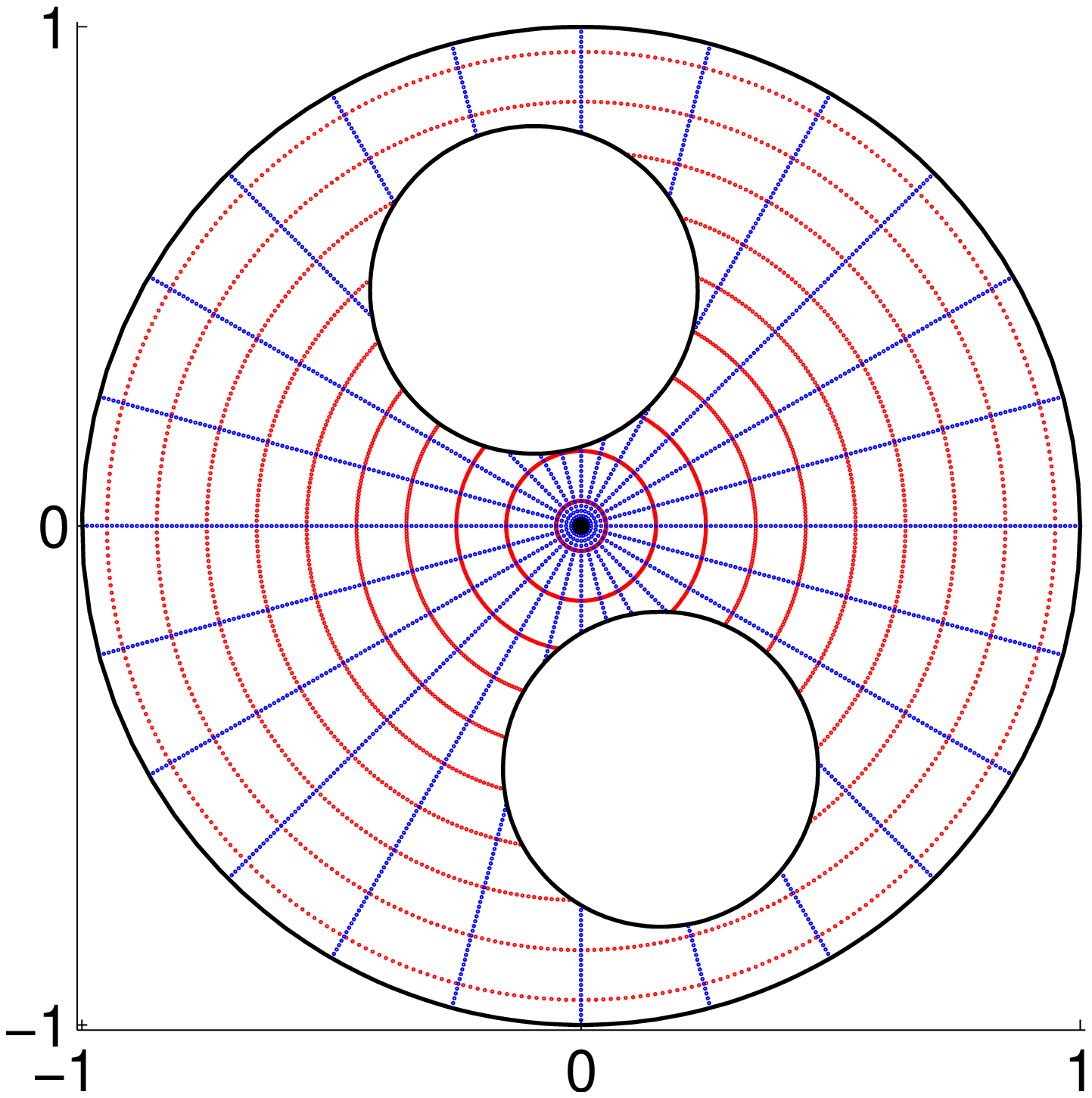}}
\scalebox{0.235}{\includegraphics{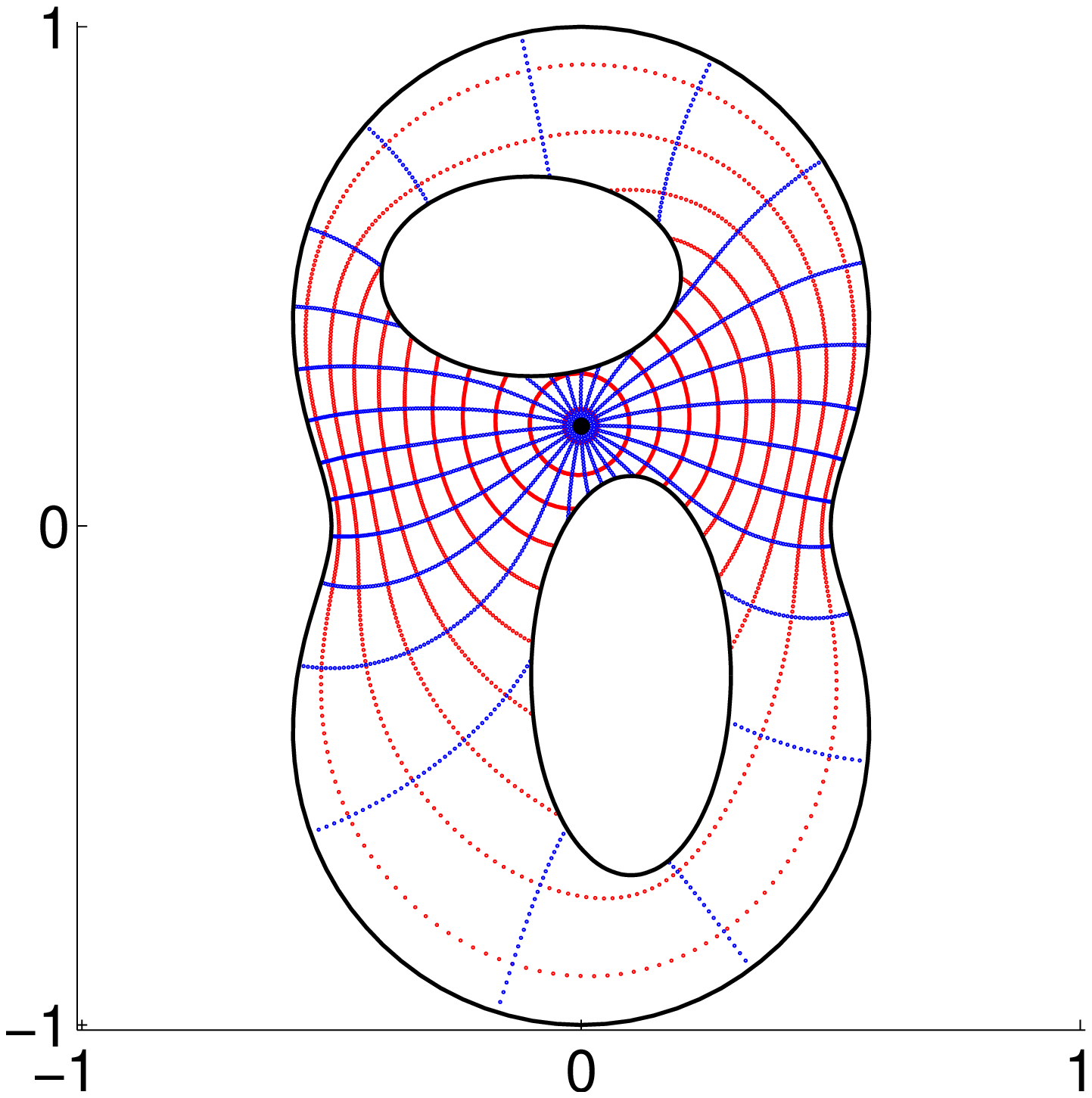}}
}
\caption{\rm The circular region $\Omega$ for Example 3 (left) and its inverse image obtained with $n=128$ (right).} 
\label{f:ex3-inv}
\end{figure}

\begin{figure}%
\centerline{
\scalebox{0.235}{\includegraphics{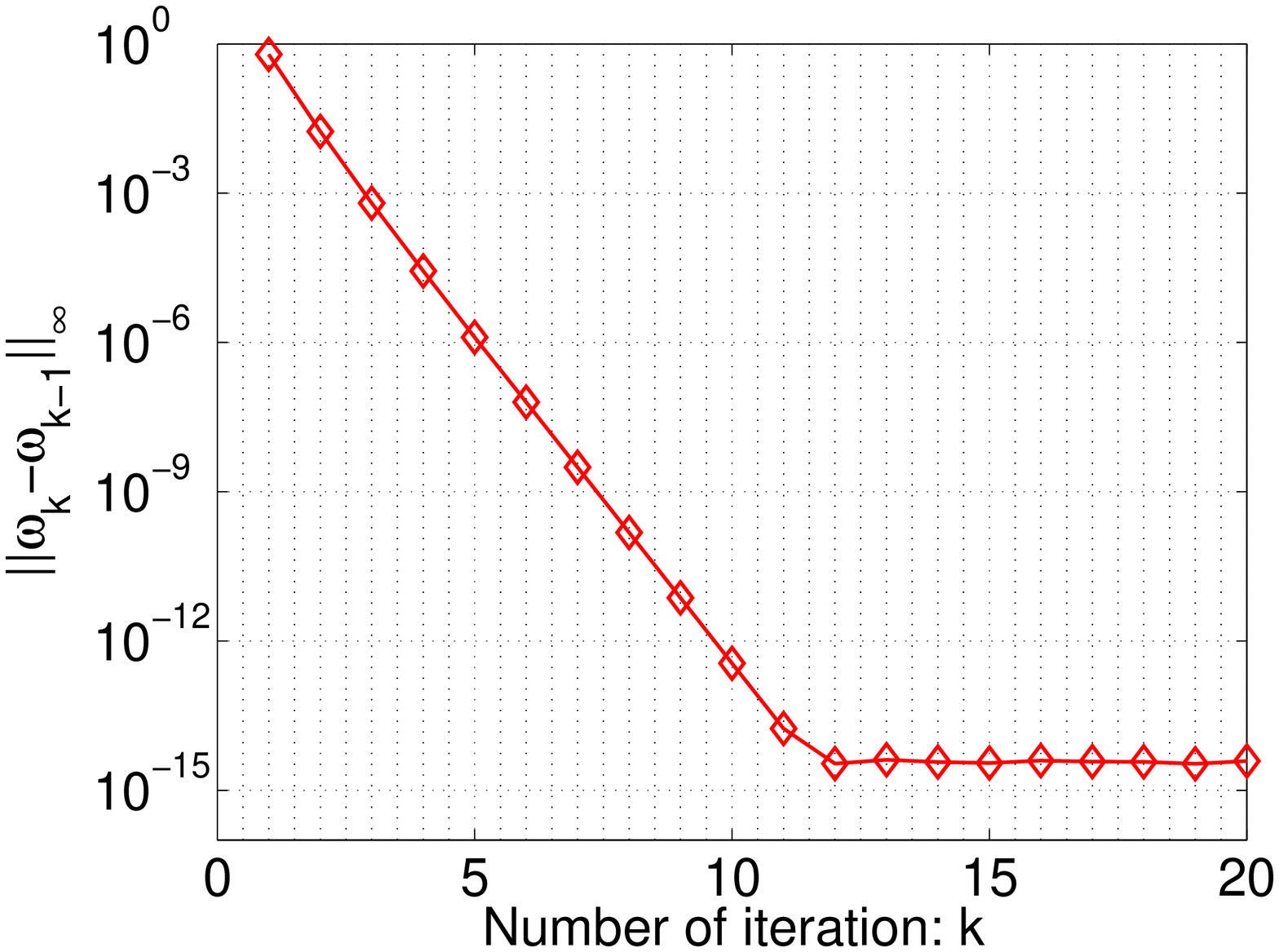}}
\scalebox{0.235}{\includegraphics{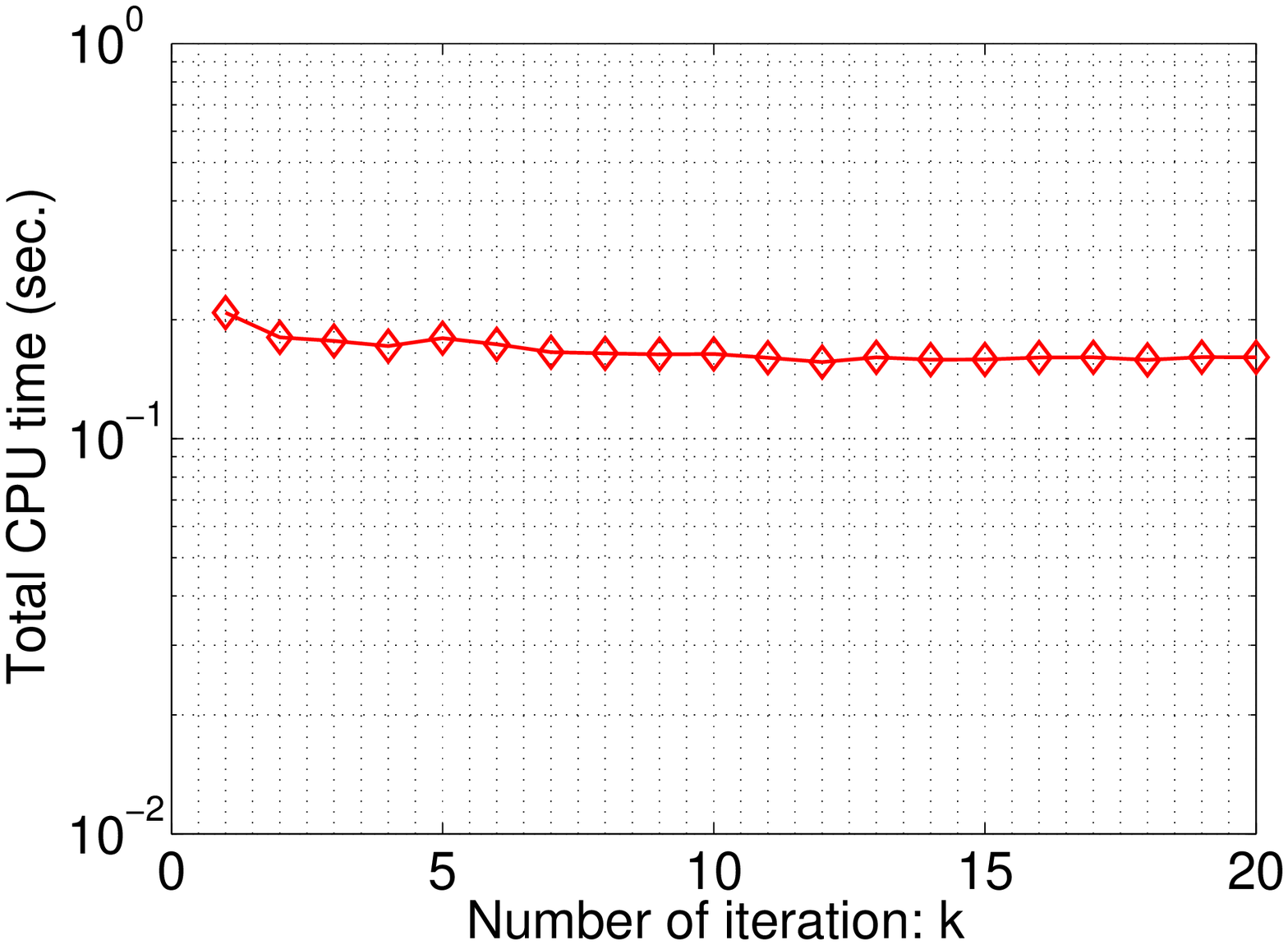}}
\scalebox{0.225}{\includegraphics{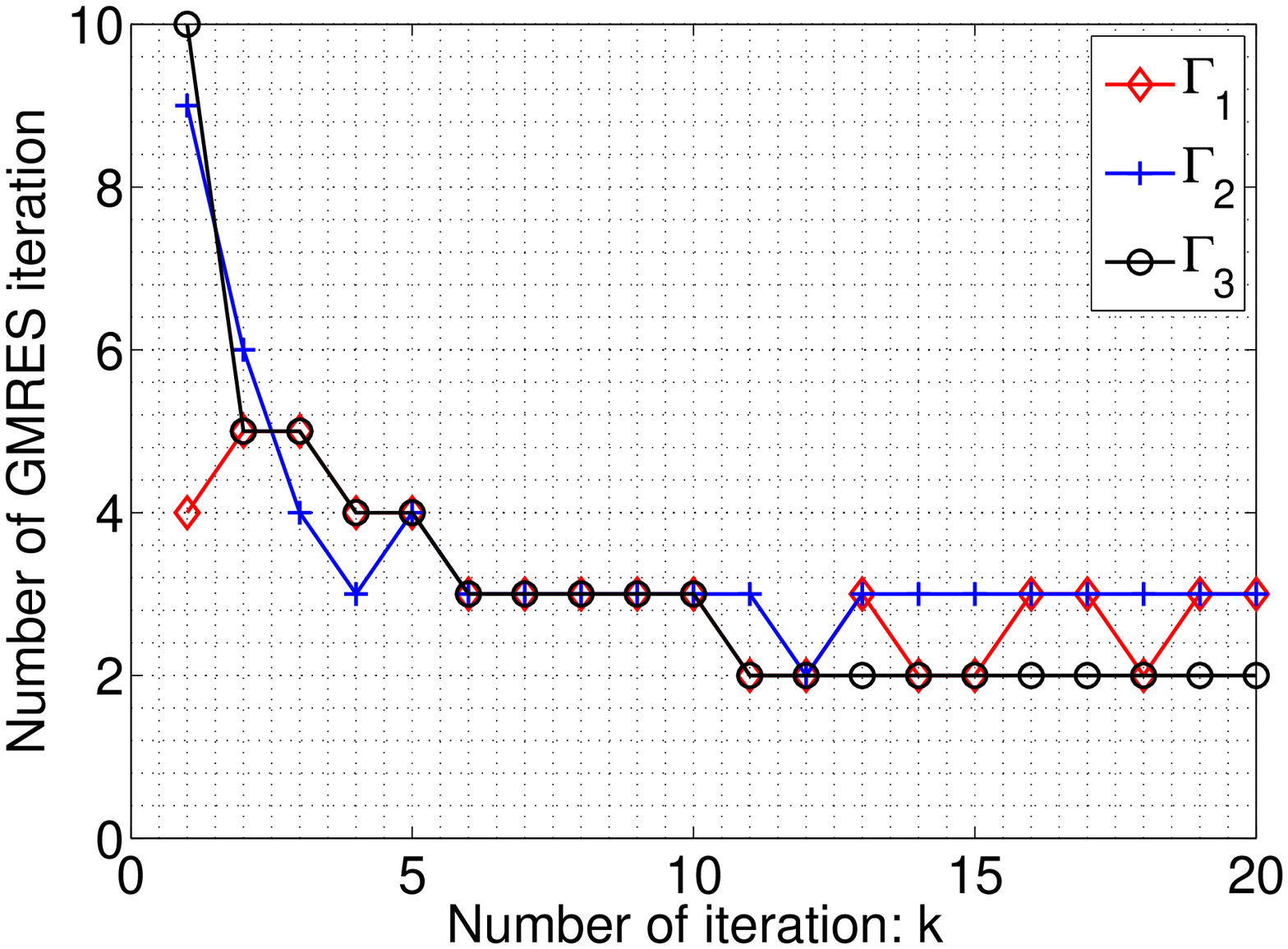}}
}
\centerline{
\scalebox{0.225}{\includegraphics{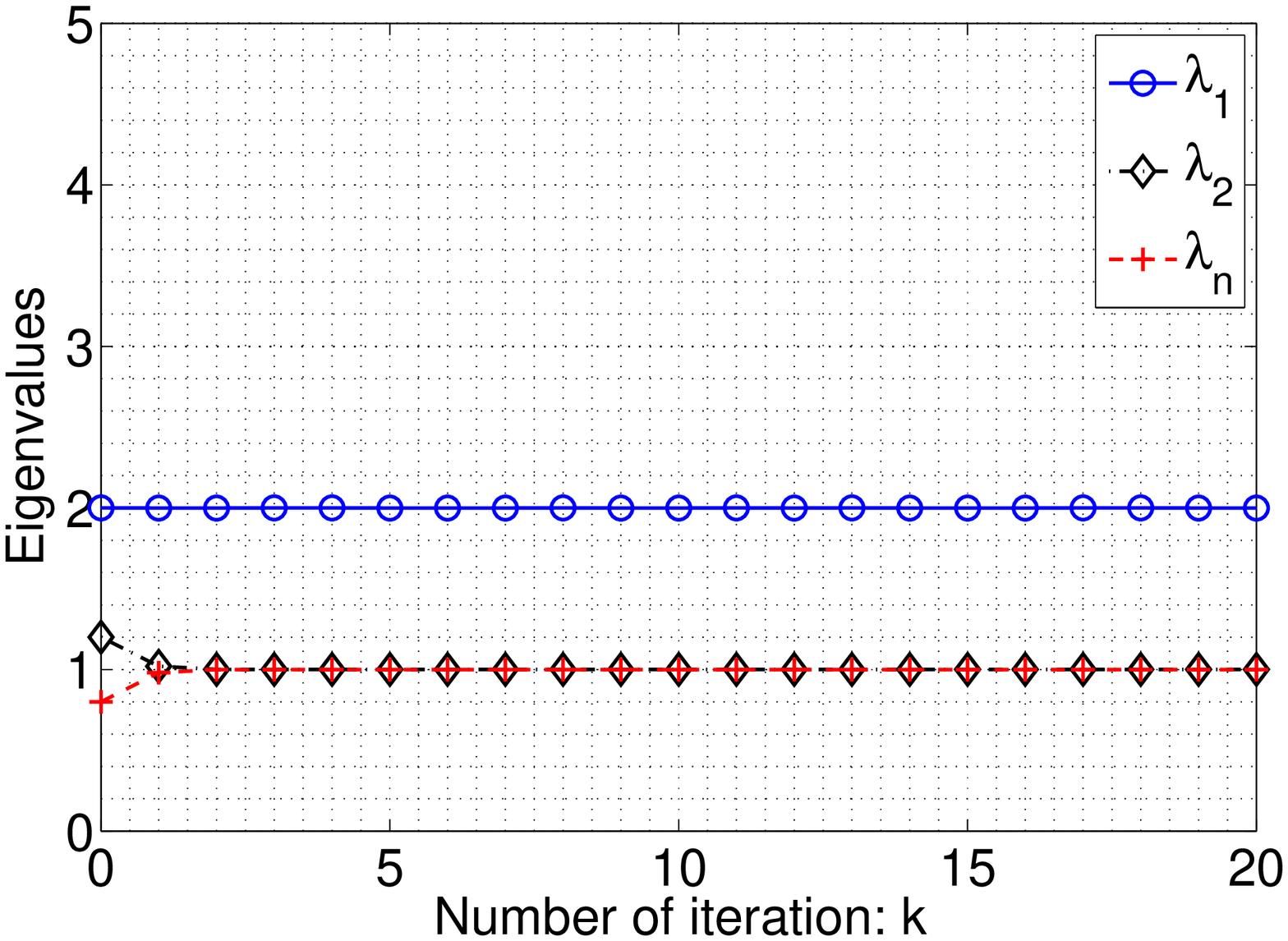}}
\scalebox{0.225}{\includegraphics{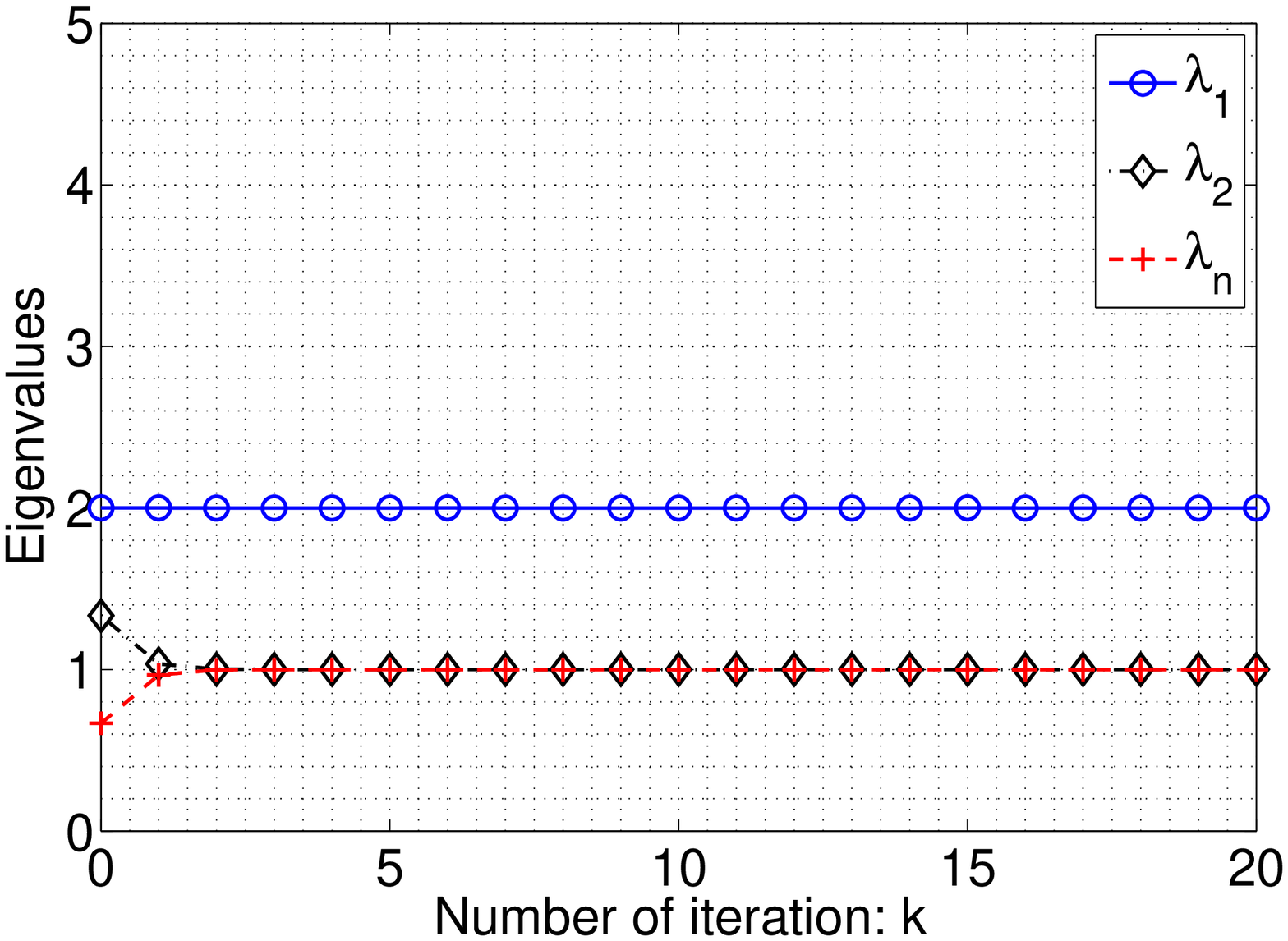}}
\scalebox{0.225}{\includegraphics{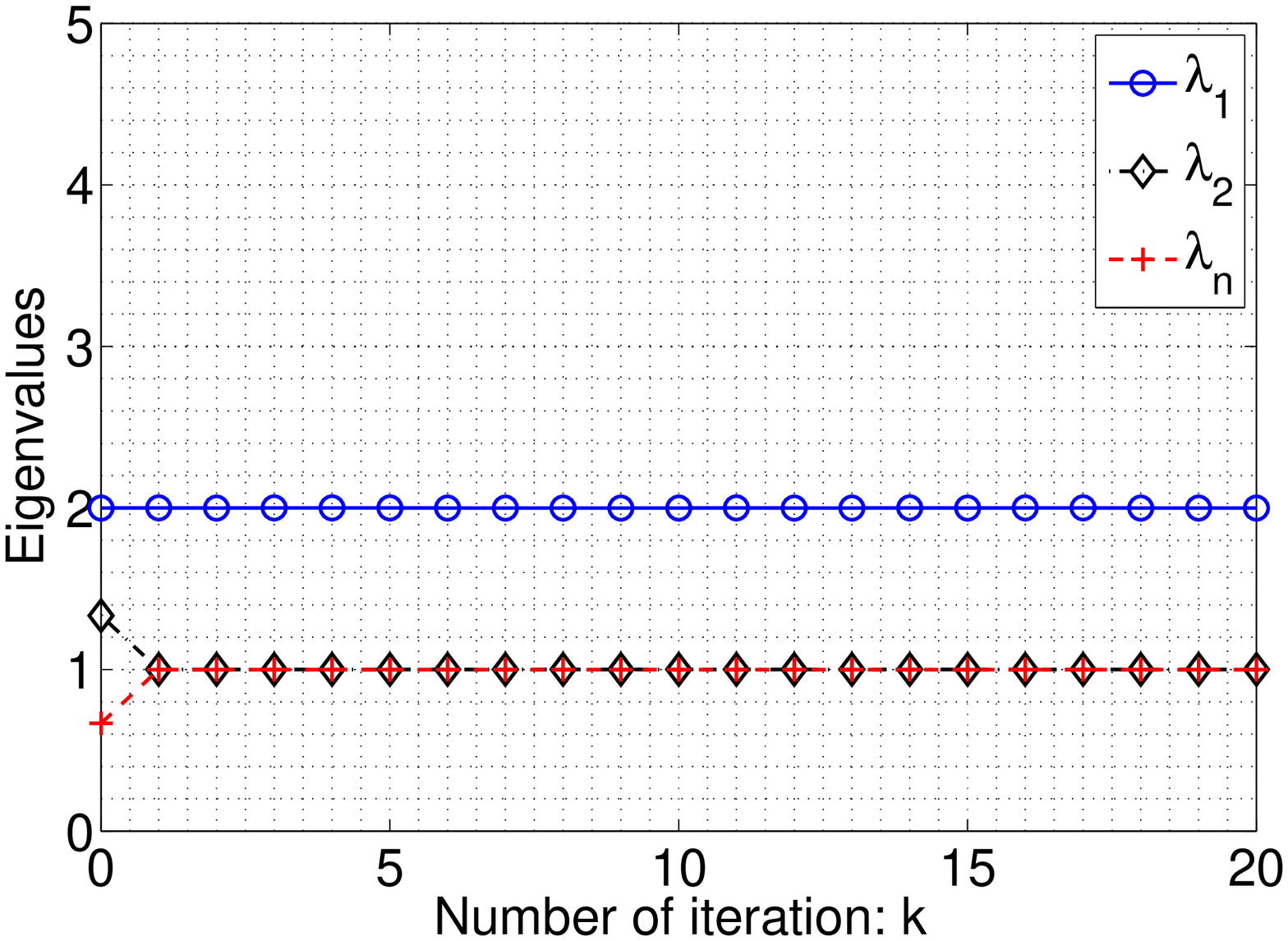}}
}
\centerline{
\scalebox{0.225}{\includegraphics{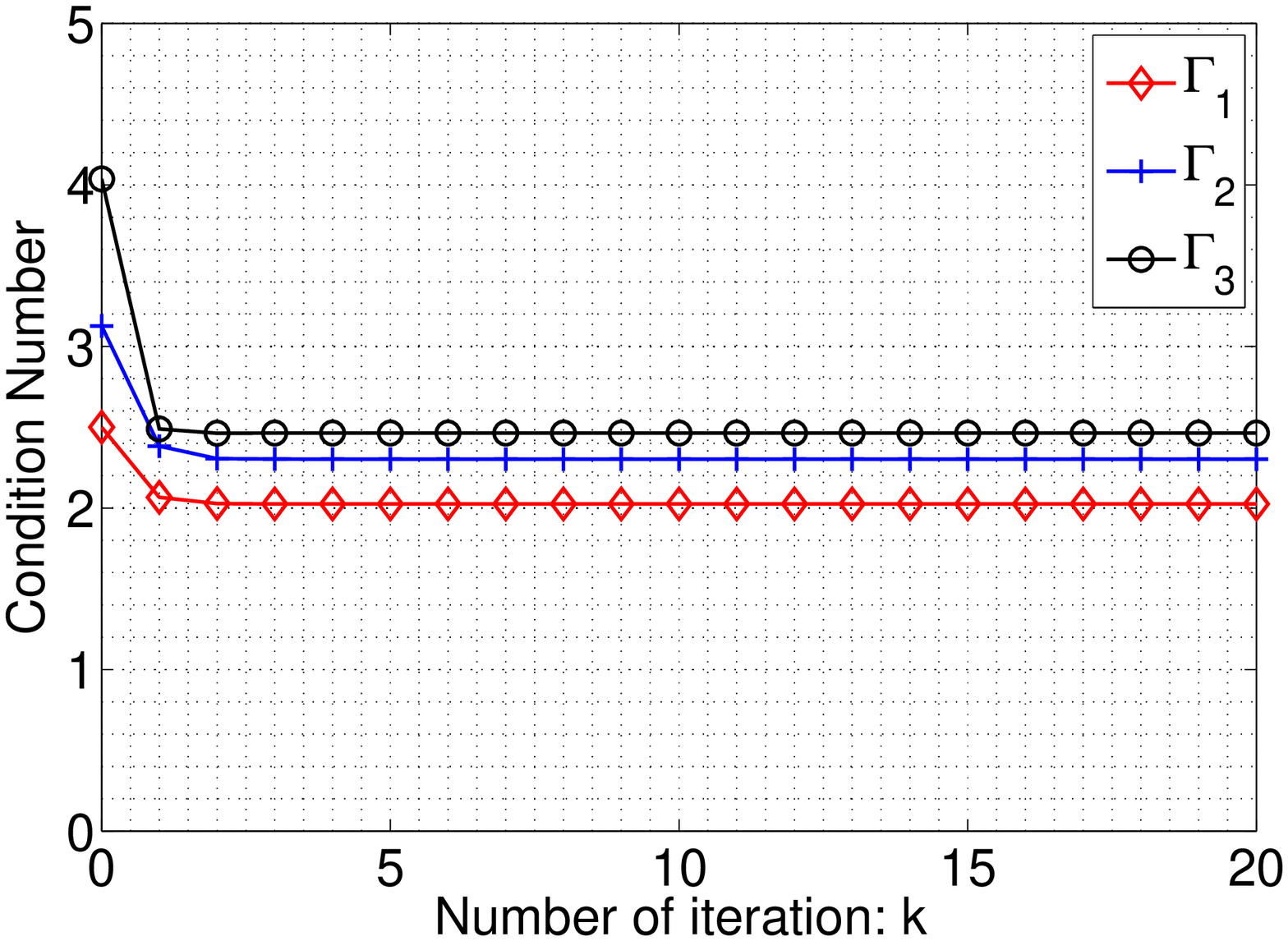}}
}
\caption{\rm Numerical results for Example 3 obtained with $n=128$. 
First row: the successive iteration errors (left), the CPU time in seconds (middle), the number of GMRES iterations (right). Second row: the eigenvalues $\lambda_1,\lambda_2,\lambda_n$ for the coefficient matrix of the linear system on $\Gamma_1$ (left), $\Gamma_2$ (middle), and $\Gamma_3$ (right). Third row: The condition number of the coefficient matrices.} 
\label{f:ex3-err}
\end{figure}

\begin{example}\label{ex:3}{\rm
In this example, we calculate the mapping function $w=\omega(z)$ that maps the region $G$ bounded by three ellipse in $z$-plane onto a bounded multiply connected circular region $\Omega$ in the $w$-plane. The same example has been considered in~\cite[Example~19]{Weg05} for computing $\omega^{-1}$ but with different normalization. The inner ellipses $\Gamma_1$ and $\Gamma_2$ are parametrized by
\begin{eqnarray*}
\Gamma_1 &:& \eta_1(t)=-0.1+0.5\i+0.3\cos t-0.2\i\sin t,\\
\Gamma_2 &:& \eta_2(t)=+0.1-0.3\i+0.2\cos t-0.4\i\sin t, 
\end{eqnarray*}
for $0\le t\le 2\pi$. The external boundary $\Gamma_3$ is the inverted ellipse parametrized by 
\[
\Gamma_3\;\;:\;\; \eta_3(t)=\sqrt{1-(1-p^2)\cos^2t}\;e^{\i t}, \quad p=0.5.
\]
The numerical results are shown in Figures~\ref{f:ex3-im}--\ref{f:ex3-err}.
}\end{example}

\begin{figure}%
\centerline{
\scalebox{0.235}{\includegraphics{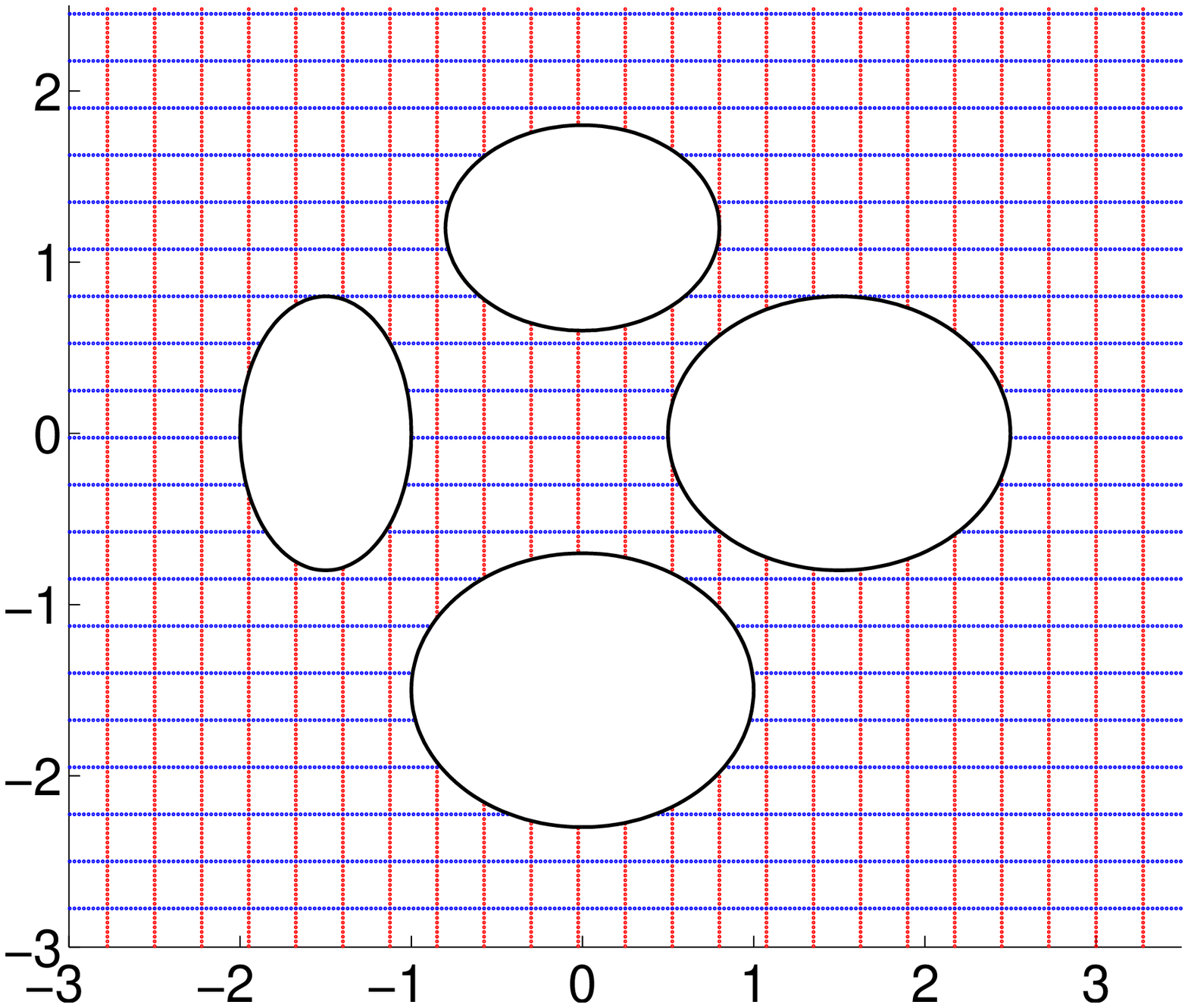}}
\scalebox{0.235}{\includegraphics{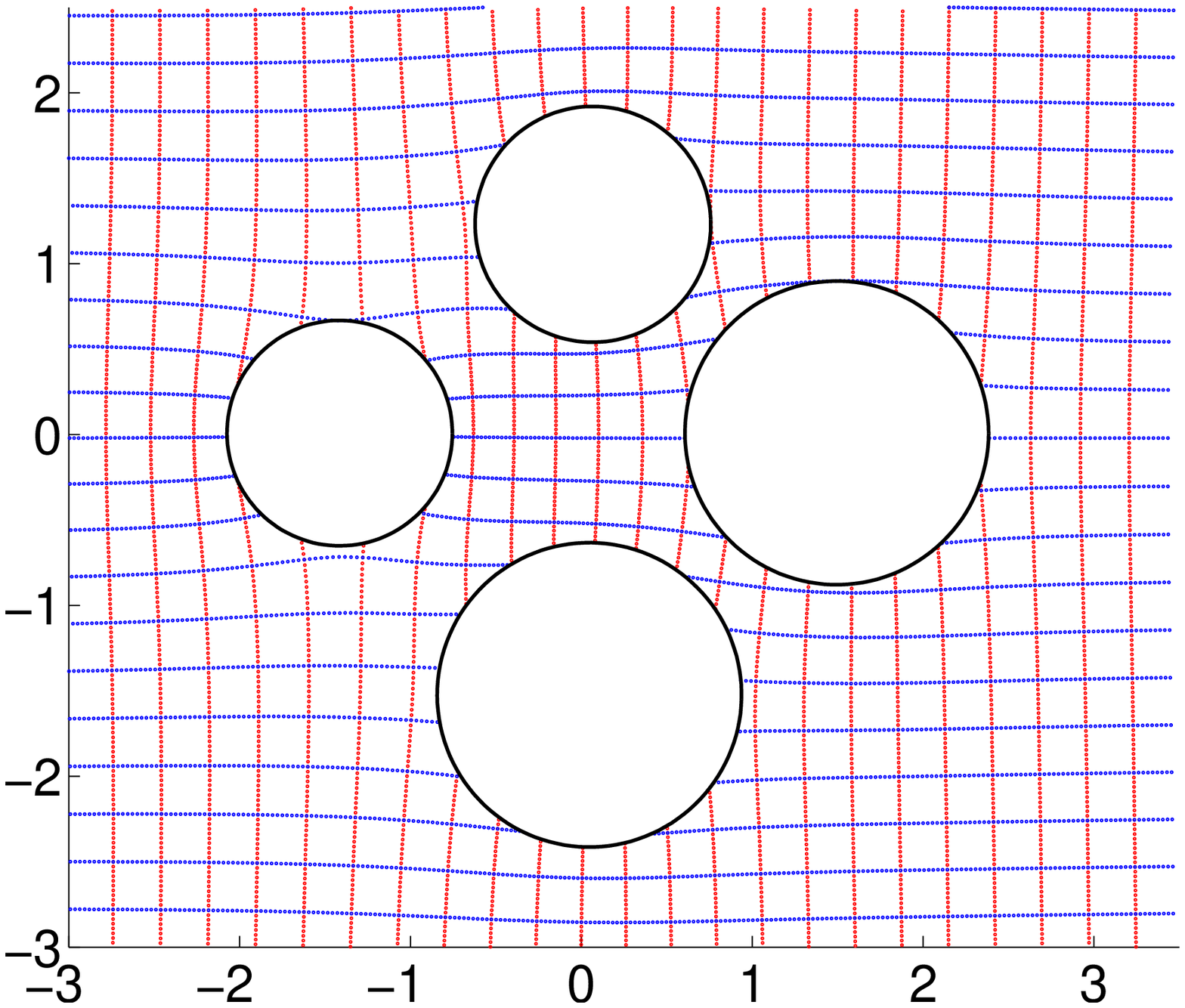}}
}
\caption{\rm The original region $G$ for Example 4 (left) and its image obtained with $n=128$ (right).} 
\label{f:ex4-im}
\end{figure}

\begin{figure}%
\centerline{
\scalebox{0.235}{\includegraphics{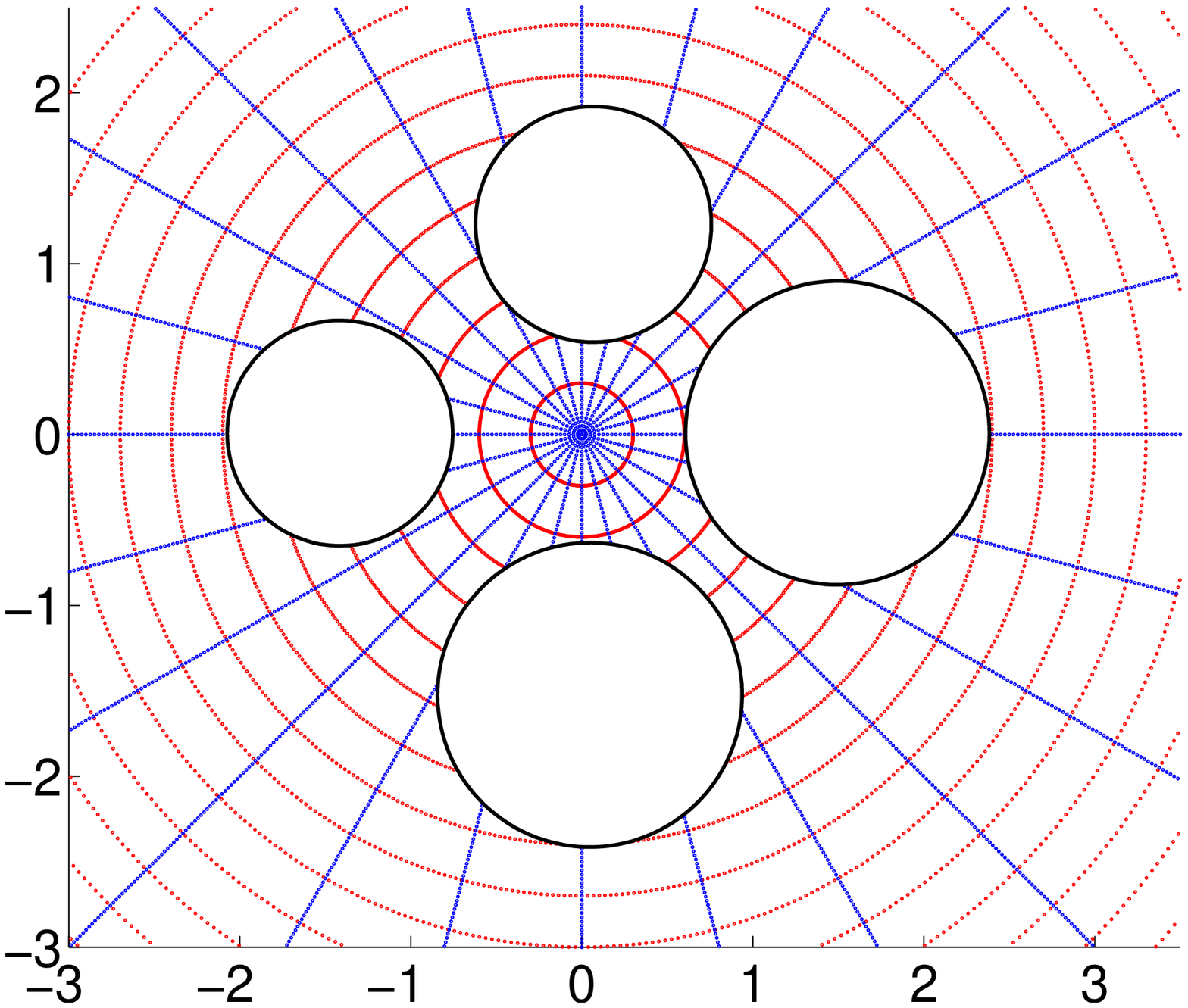}}
\scalebox{0.235}{\includegraphics{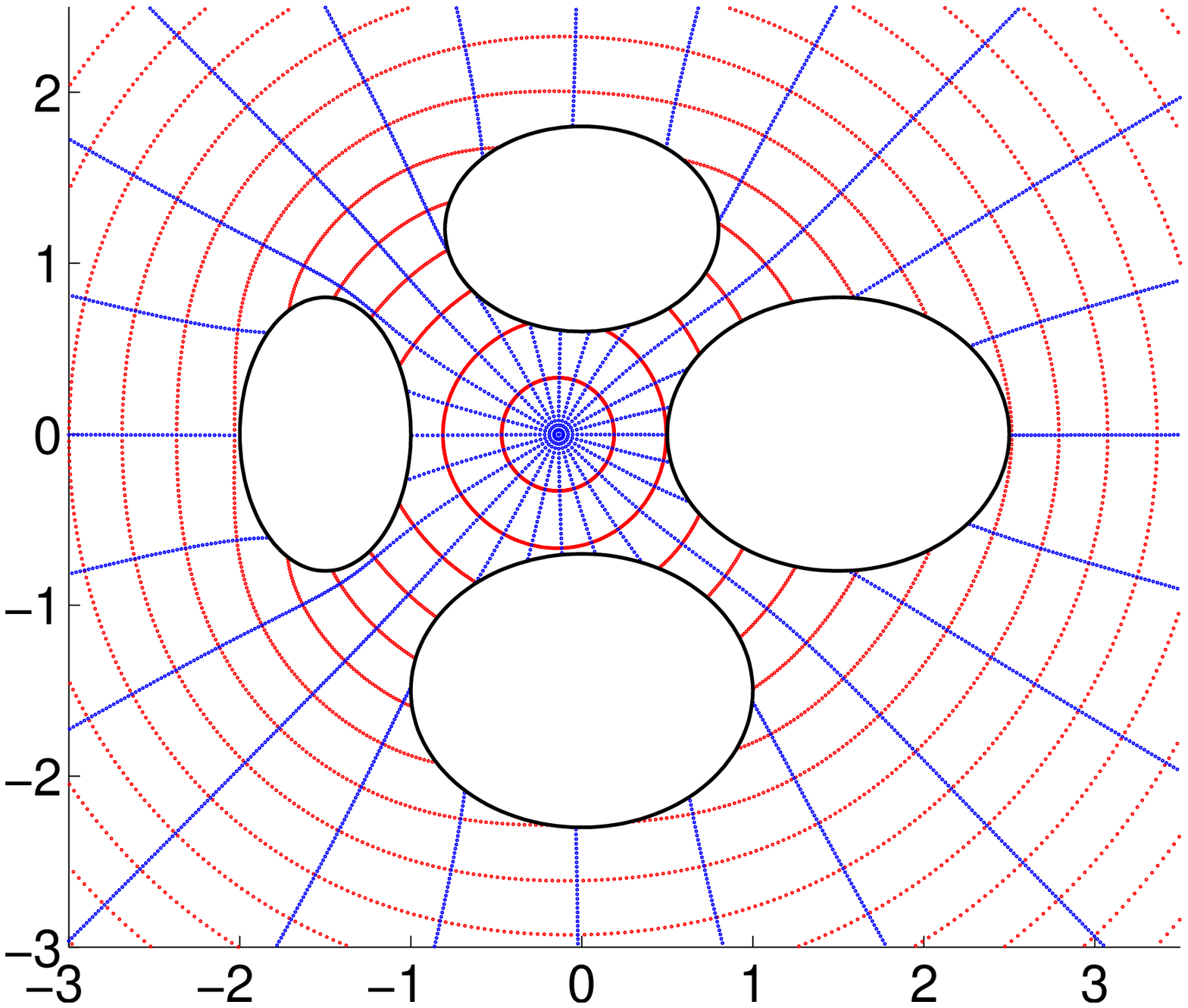}}
}
\caption{\rm The circular region $\Omega$ for Example 4 (left) and its inverse image obtained with $n=128$ (right).} 
\label{f:ex4-inv}
\end{figure}

\begin{figure}%
\centerline{
\scalebox{0.235}{\includegraphics{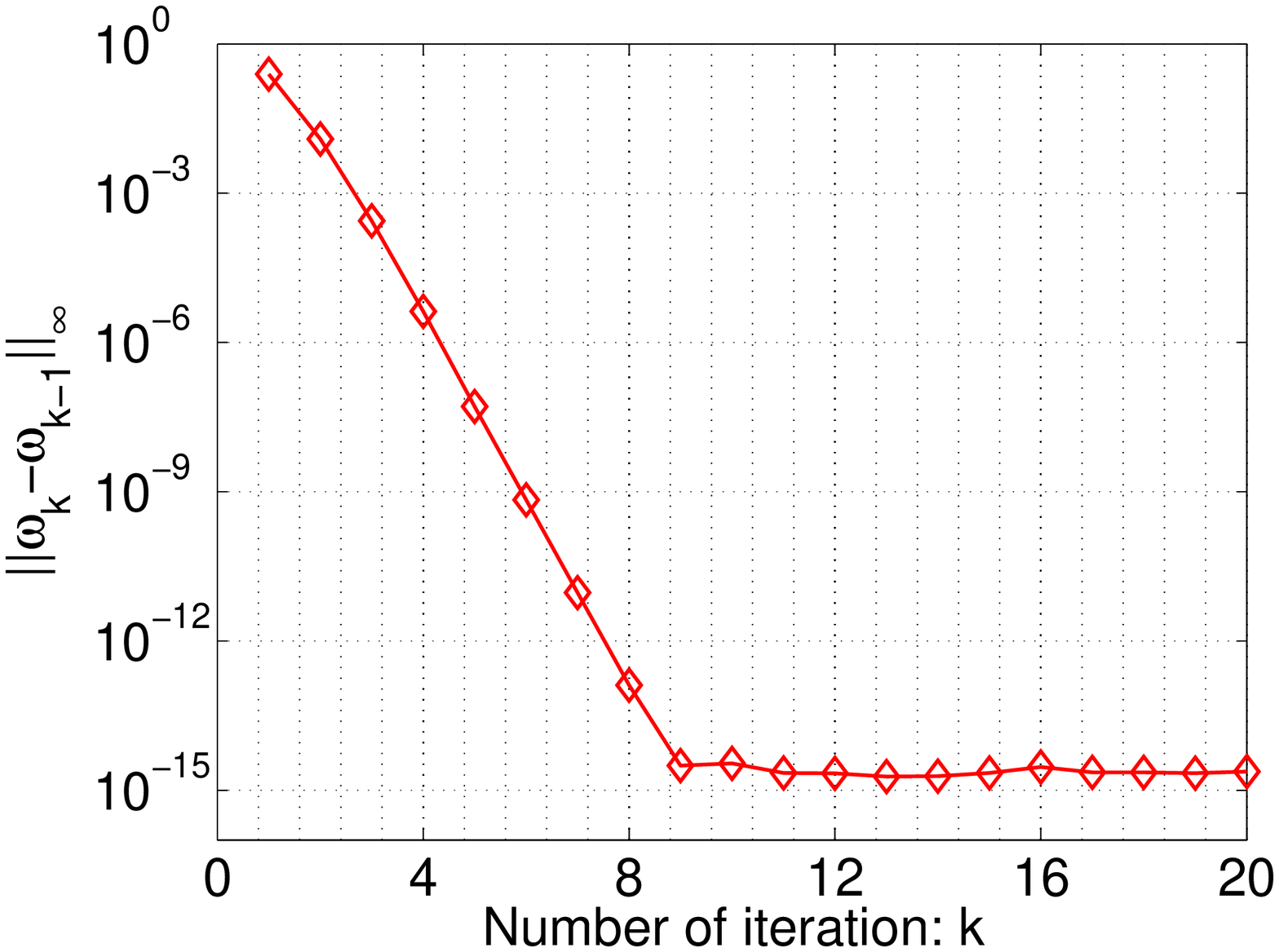}}
\scalebox{0.235}{\includegraphics{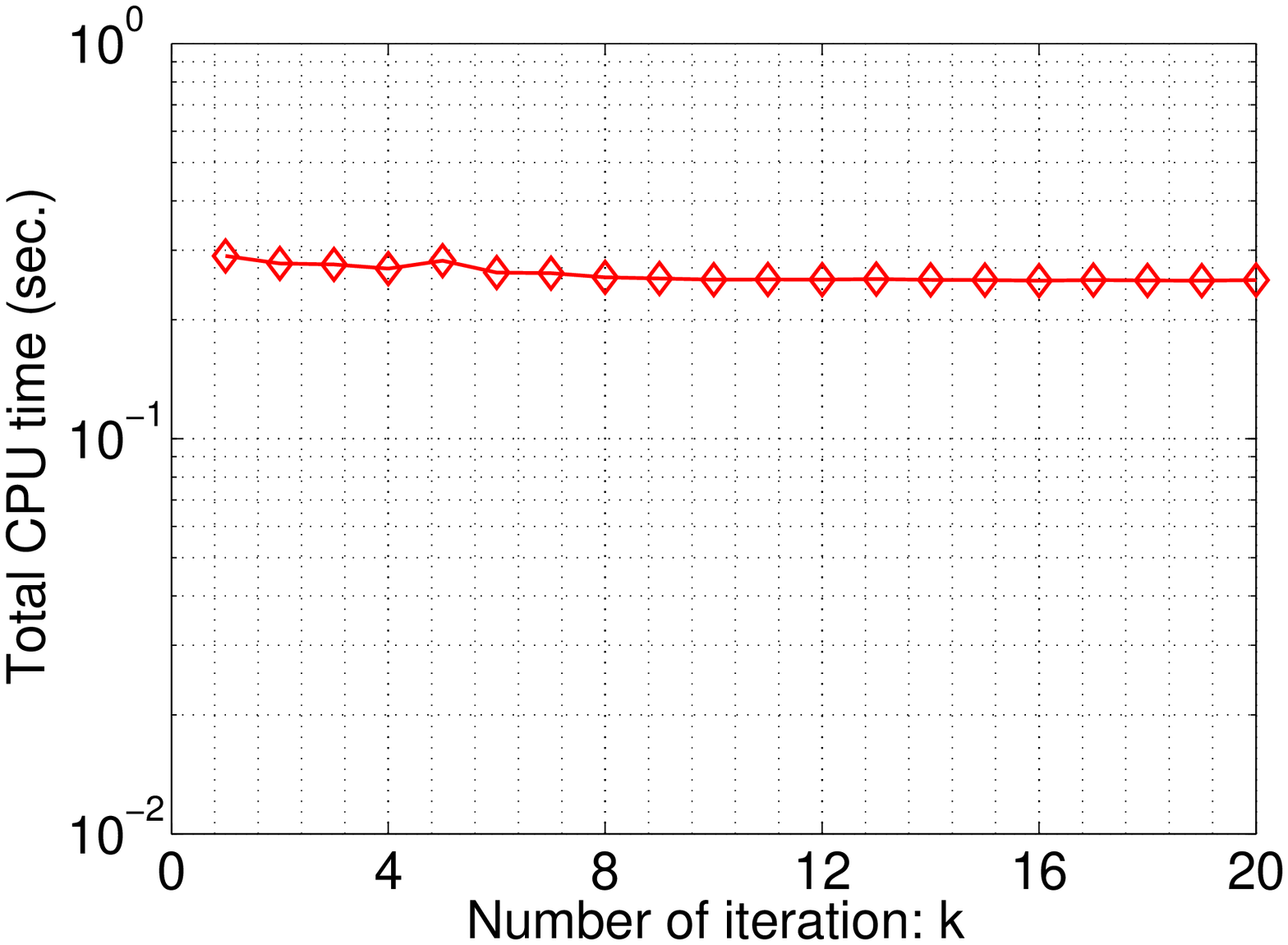}}
}
\caption{\rm The successive iteration errors (left) and the CPU time in seconds (right) obtained for Example 4 with $n=128$.} 
\label{f:ex4-err}
\end{figure}

\begin{example}\label{ex:4}{\rm
In this example, we calculate the mapping function $w=\omega(z)$ maps the unbounded region $G$ exterior to four ellipse in $z$-plane onto an unbounded multiply~connected circular region $\Omega$ in the $w$-plane. This example has been considered in~\cite[Example~2]{Weg01c} for computing $\omega^{-1}$. The boundaries are parametrized by
\begin{eqnarray*}
\Gamma_1 &:& \eta_1(t)= 1.5  +   \cos t-0.8\i\sin t,\\
\Gamma_2 &:& \eta_2(t)= 1.2\i+0.8\cos t-0.6\i\sin t,\\ 
\Gamma_3 &:& \eta_3(t)=-1.5  +0.5\cos t-0.8\i\sin t,\\ 
\Gamma_4 &:& \eta_4(t)=-1.5\i+   \cos t-0.8\i\sin t, \quad 0\le t\le 2\pi.
\end{eqnarray*}
\nopagebreak[0]The numerical results are shown in Figures~\ref{f:ex4-im}--\ref{f:ex4-err}. 
}\end{example}

\begin{figure}%
\centerline{
\scalebox{0.235}{\includegraphics{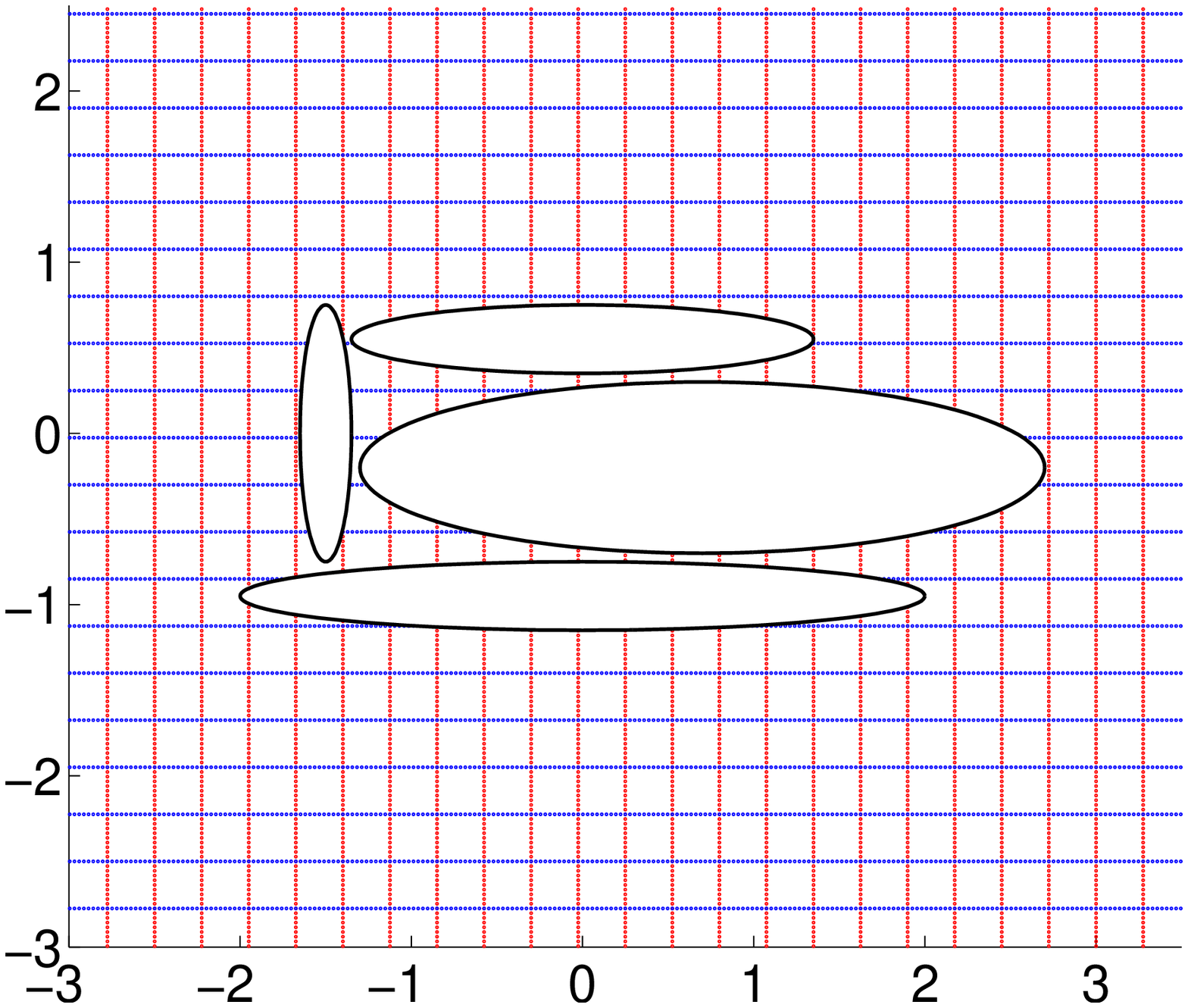}}
\scalebox{0.235}{\includegraphics{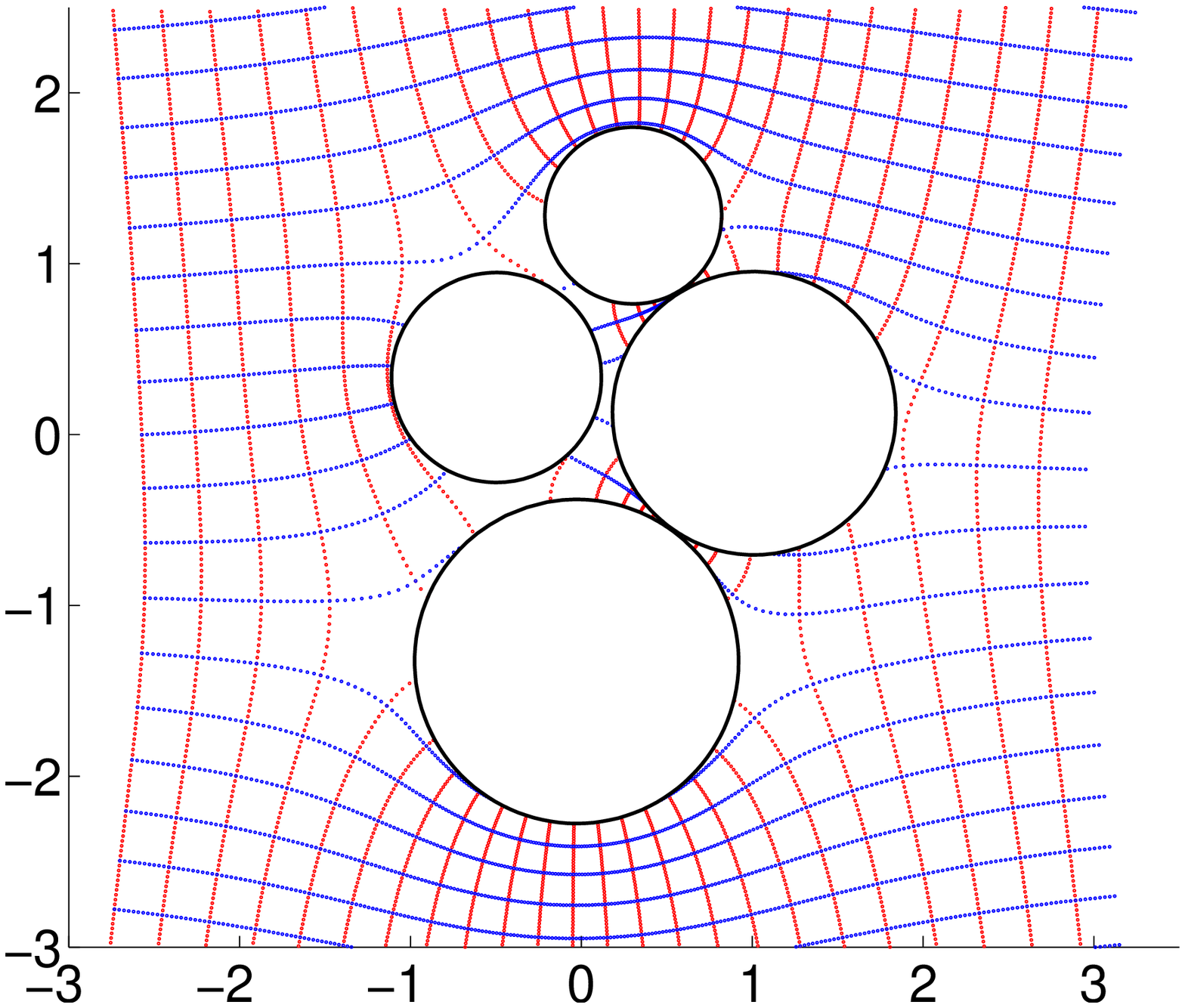}}
}
\caption{\rm The original region $G$ for Example 5 (left) and its image obtained with $n=128$ (right).} 
\label{f:ex5-im}
\end{figure}

\begin{figure}%
\centerline{
\scalebox{0.235}{\includegraphics{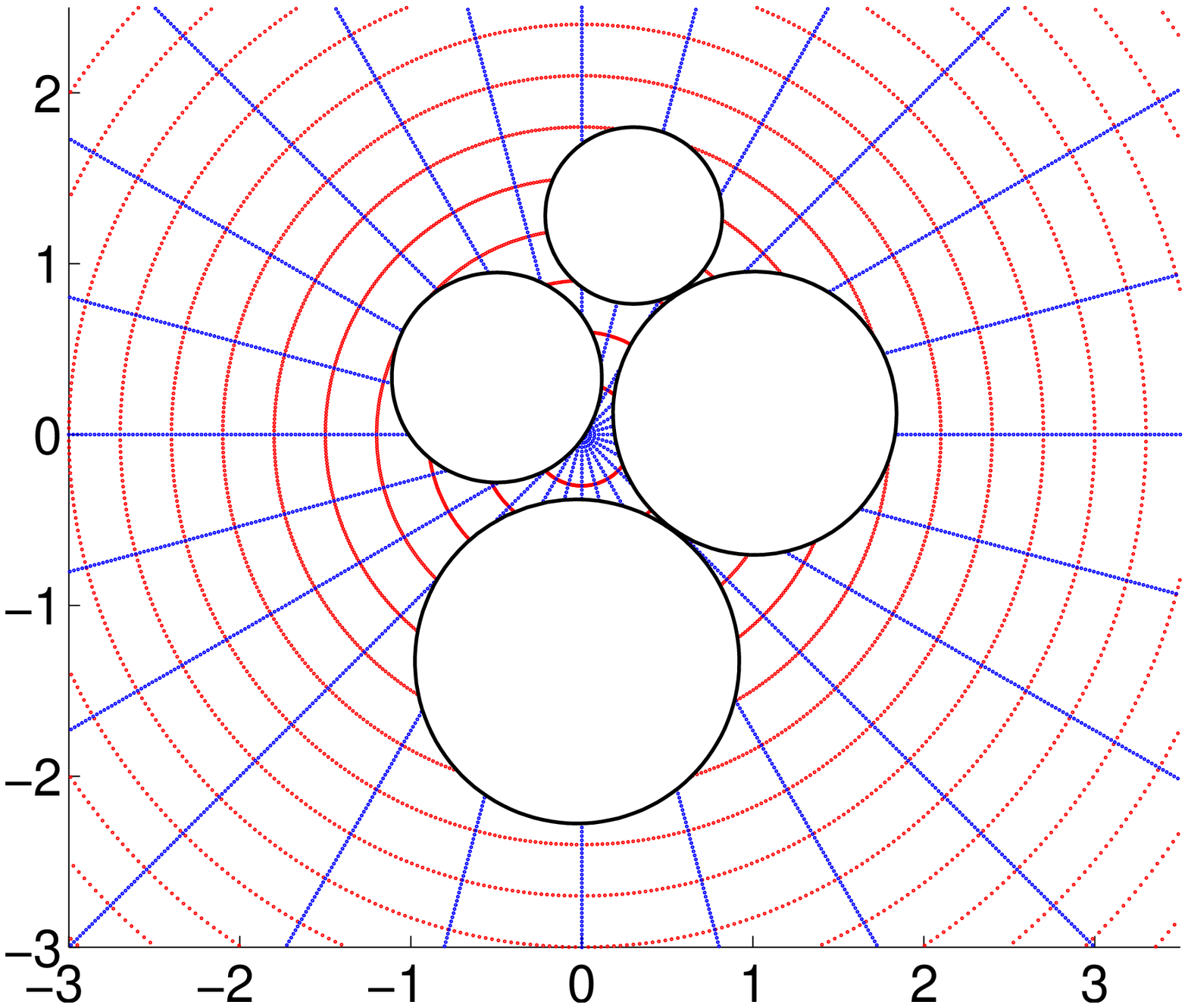}}
\scalebox{0.235}{\includegraphics{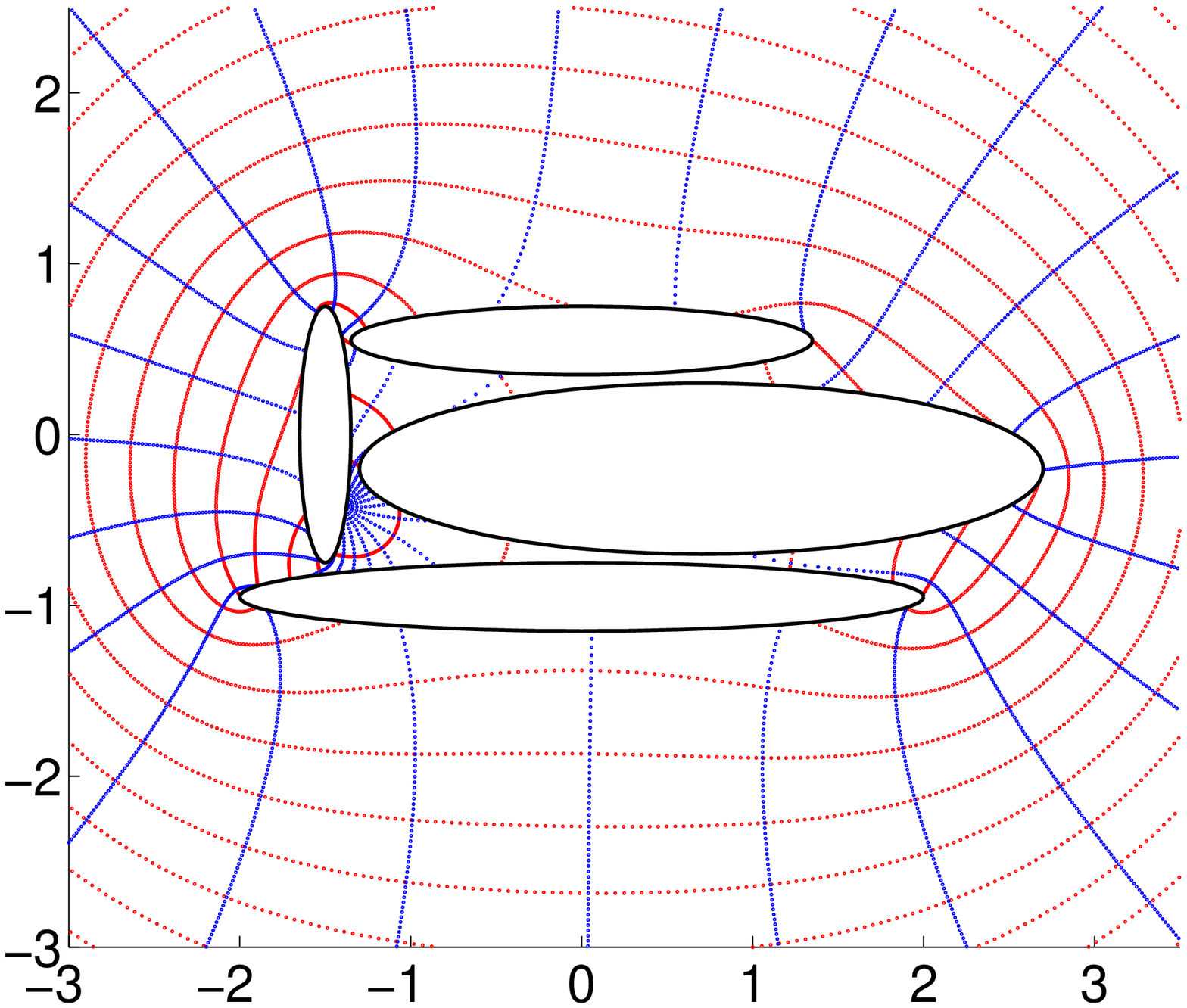}}
}
\caption{\rm The circular region $\Omega$ for Example 5 (left) and its inverse image obtained with $n=128$ (right).} 
\label{f:ex5-inv}
\end{figure}

\begin{figure}%
\centerline{
\scalebox{0.235}{\includegraphics{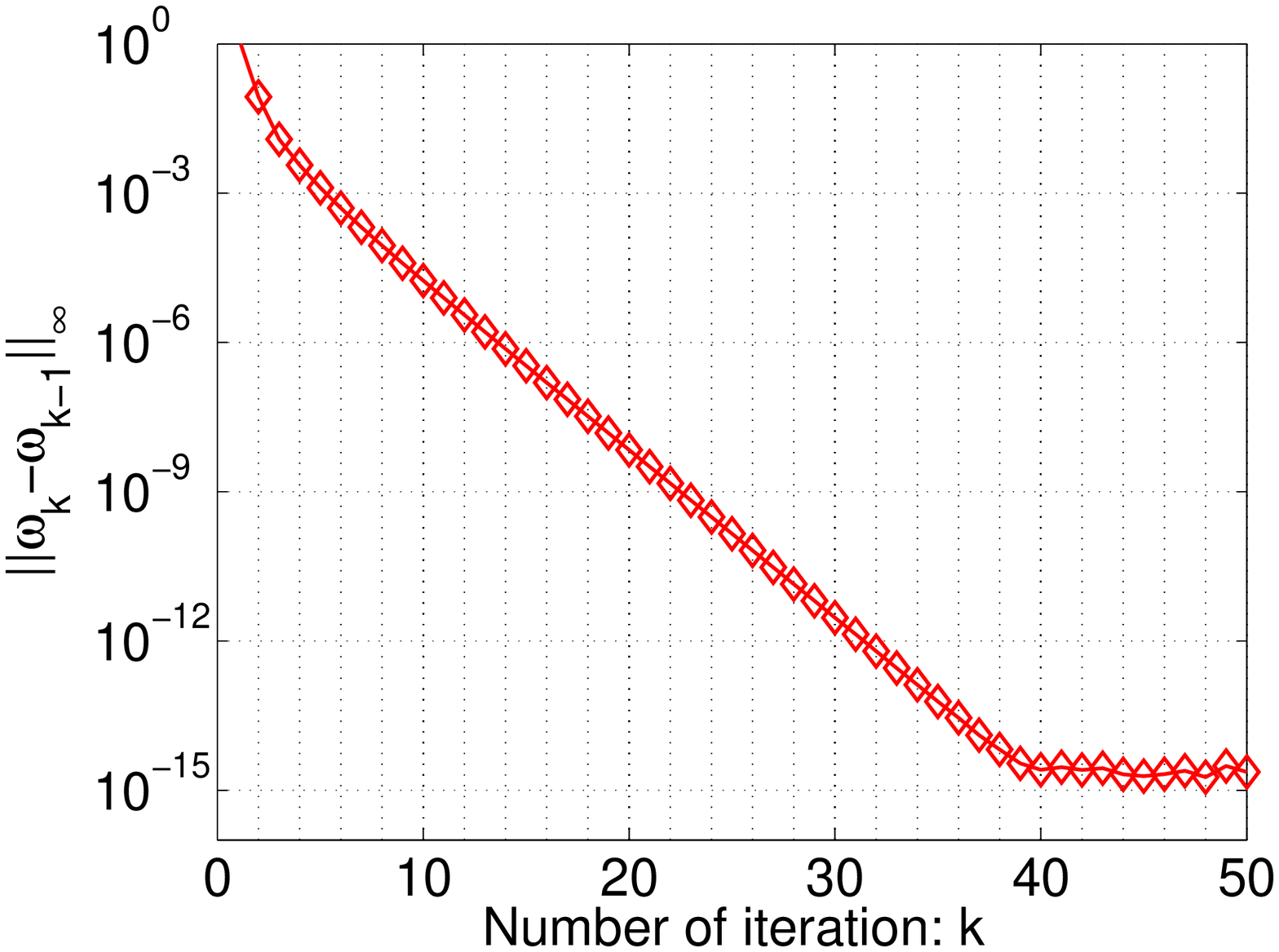}}
\scalebox{0.235}{\includegraphics{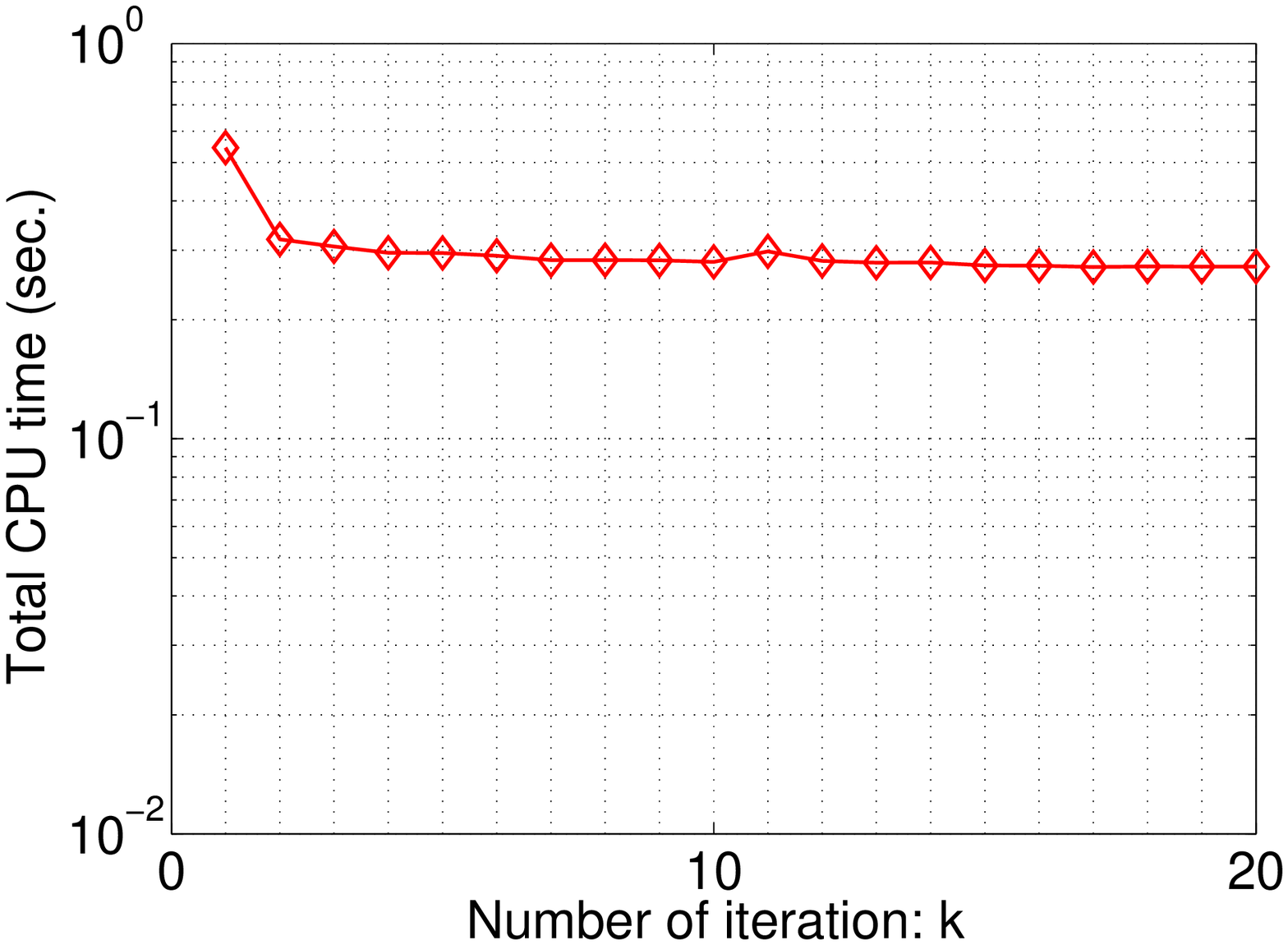}}
}
\caption{\rm The successive iteration errors (left) and the CPU time in seconds (right) obtained for Example 5 with $n=128$.} 
\label{f:ex5-err}
\end{figure}

\begin{example}\label{ex:5}{\rm
In this example, we make the ellipses in the Example~4 much thinner and closer together. The boundaries are parametrized by
\begin{eqnarray*}
\Gamma_1 &:& \eta_1(t)= 0.7-0.20\i+   2\cos t-0.5 \i\sin t,\\
\Gamma_2 &:& \eta_2(t)= 0.55\i    +1.35\cos t-0.2 \i\sin t,\\ 
\Gamma_3 &:& \eta_3(t)=-1.5       +0.15\cos t-0.75\i\sin t,\\ 
\Gamma_4 &:& \eta_4(t)=-0.95\i    +   2\cos t-0.2 \i\sin t, 
\end{eqnarray*}
for $0\le t\le 2\pi$. The numerical results are shown in Figures~\ref{f:ex5-im}--\ref{f:ex5-err}. 
}\end{example}

\begin{figure}%
\centerline{
\scalebox{0.26}{\includegraphics{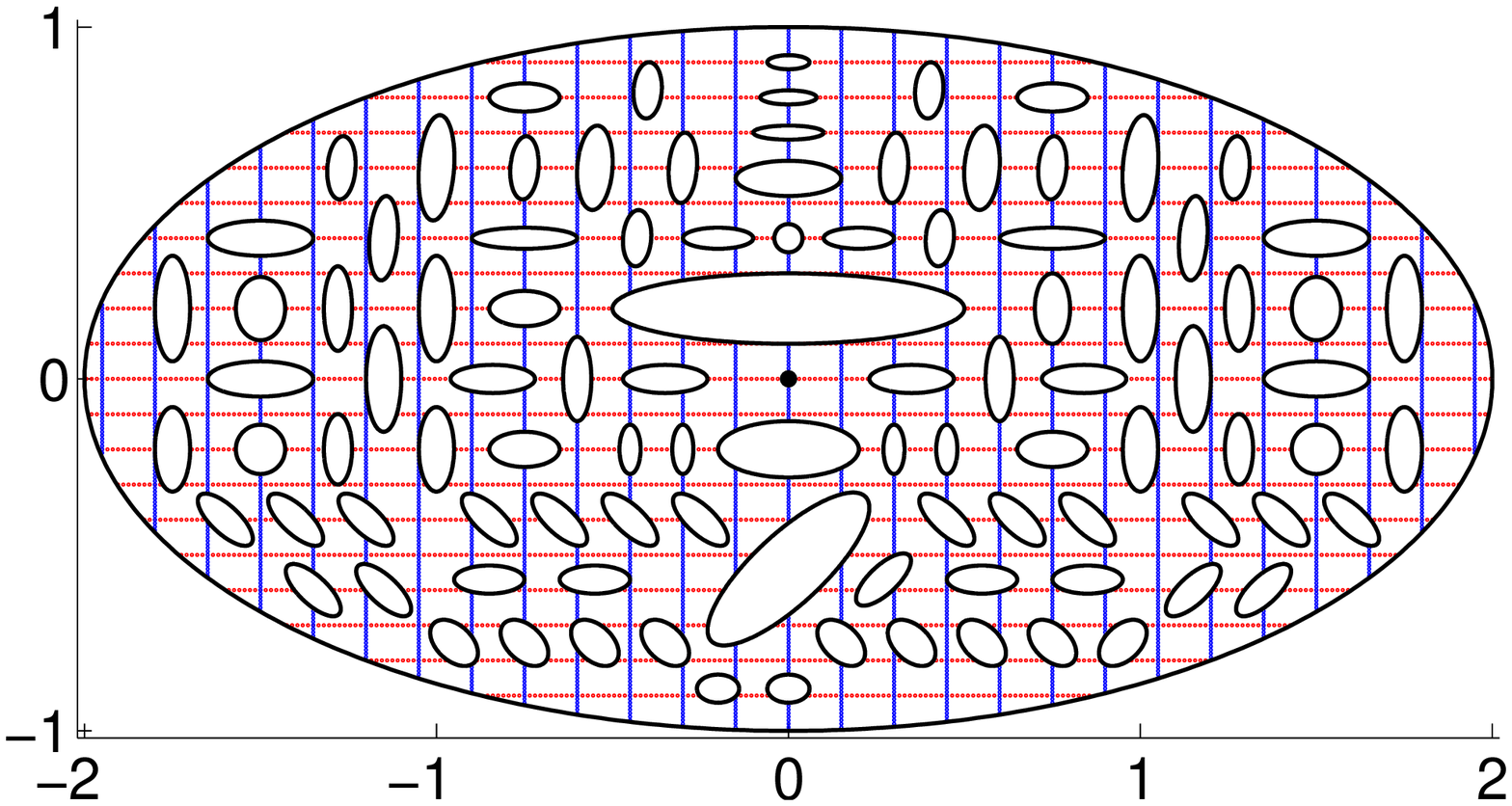}}
\scalebox{0.23}{\includegraphics{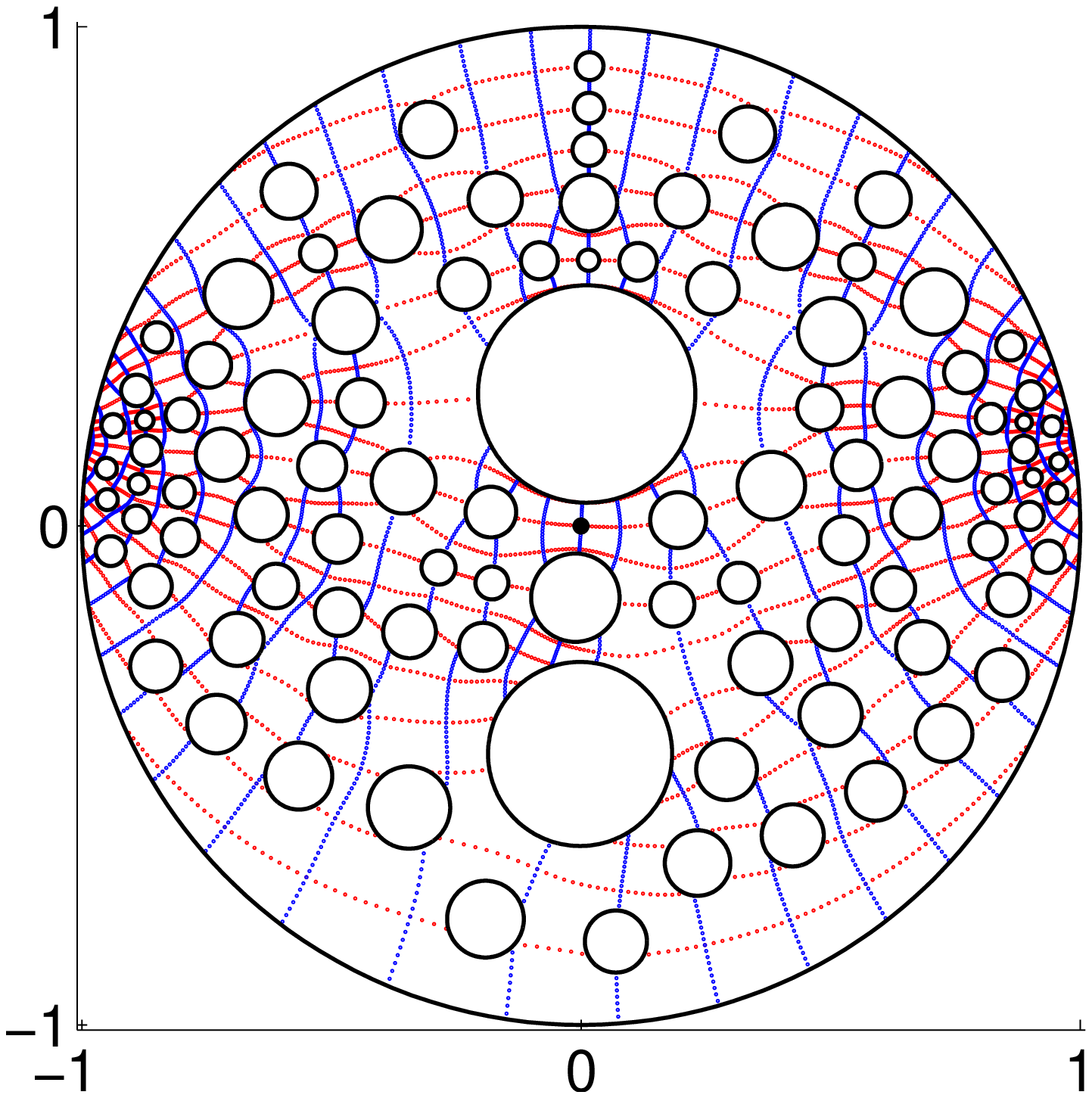}}
}
\caption{\rm The original region $G$ for Example 6 (left) and its image obtained with $n=512$ (right).} 
\label{f:ex6-im}
\end{figure}

\begin{figure}%
\centerline{
\scalebox{0.23}{\includegraphics{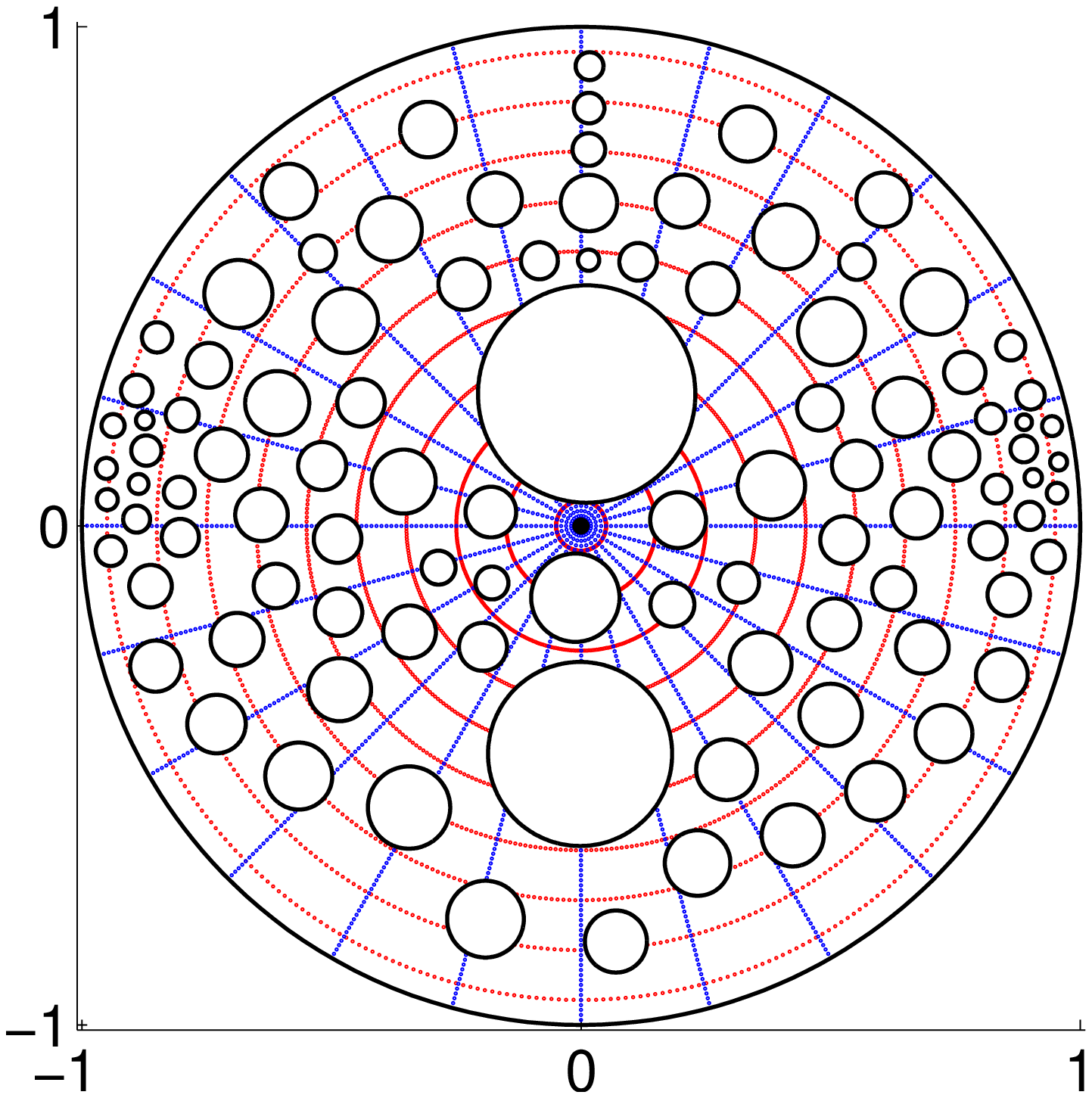}}
\scalebox{0.26}{\includegraphics{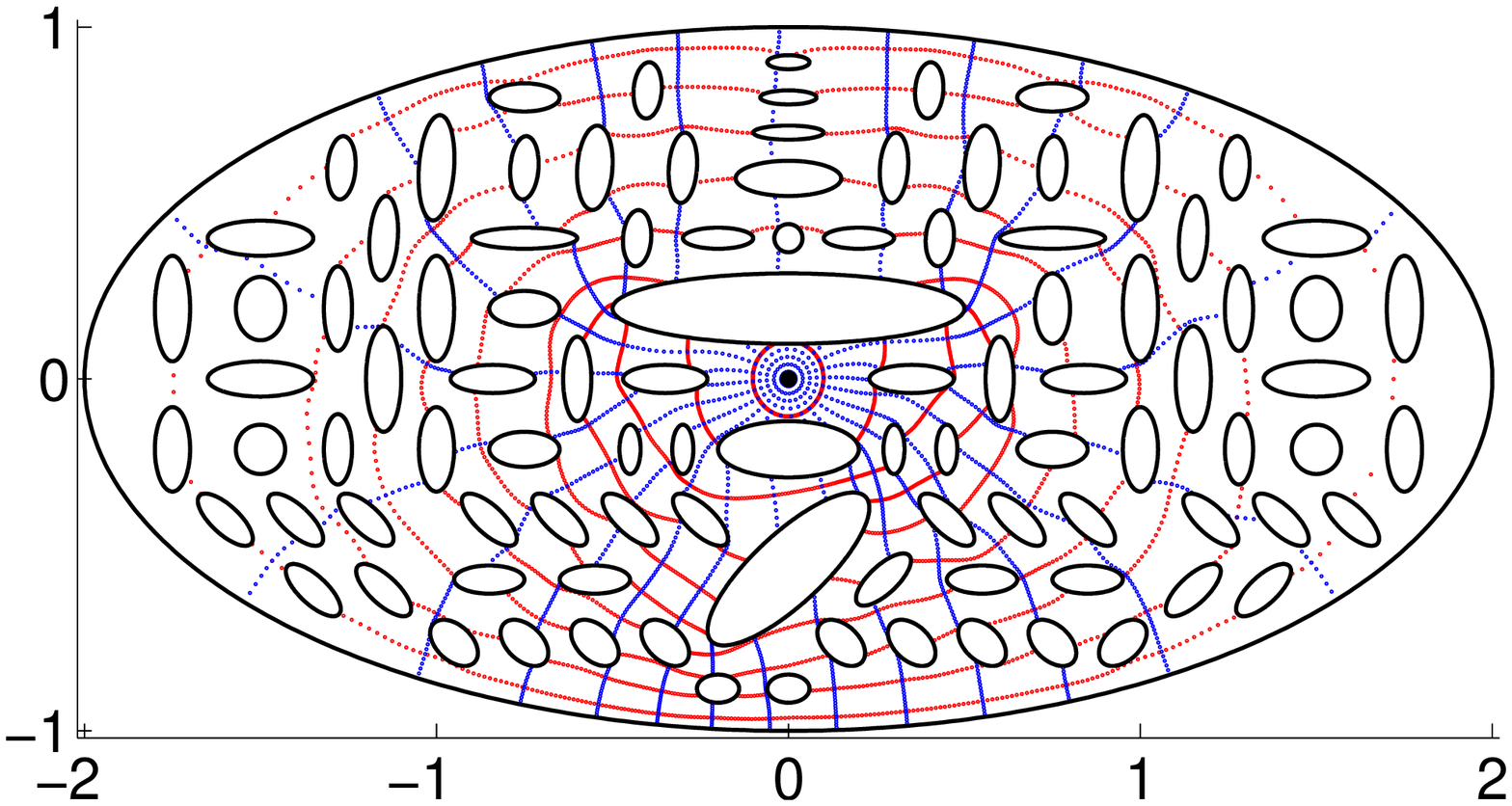}}
}
\caption{\rm The circular region $\Omega$ for Example 6 (left) and its inverse image obtained with $n=512$ (right).} 
\label{f:ex6-inv}
\end{figure}

\begin{figure}%
\centerline{
\scalebox{0.235}{\includegraphics{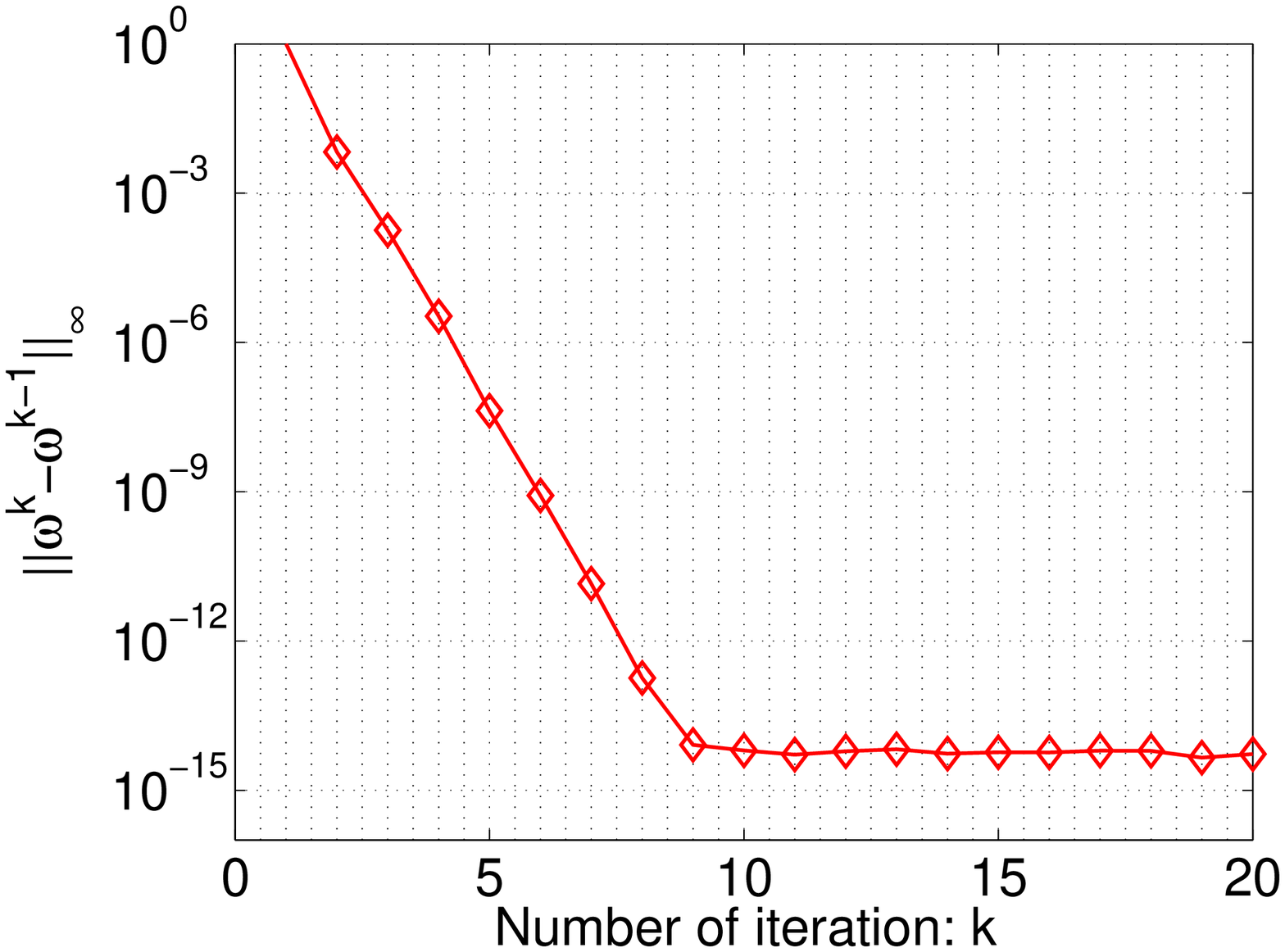}}
\scalebox{0.235}{\includegraphics{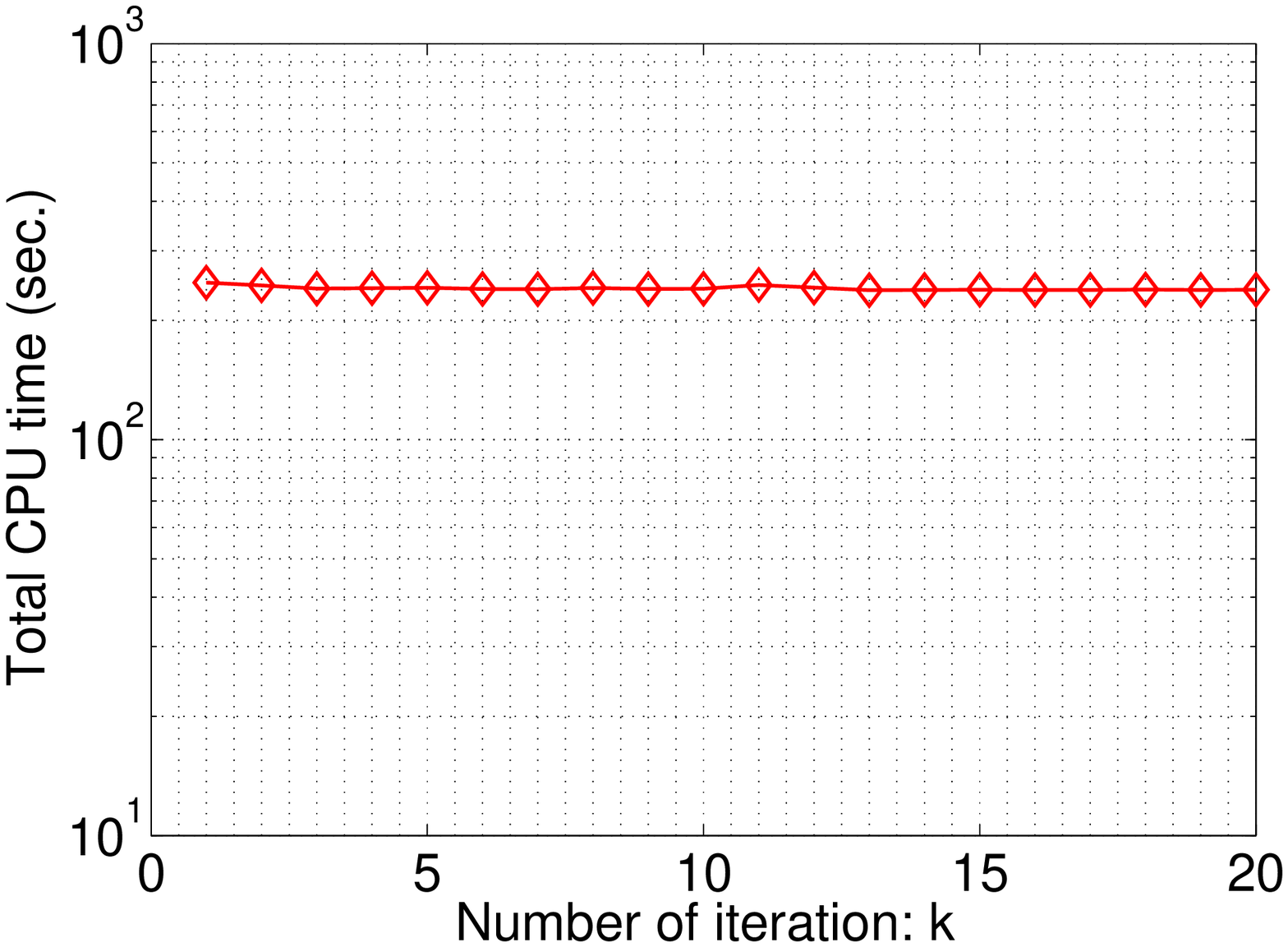}}
\scalebox{0.235}{\includegraphics{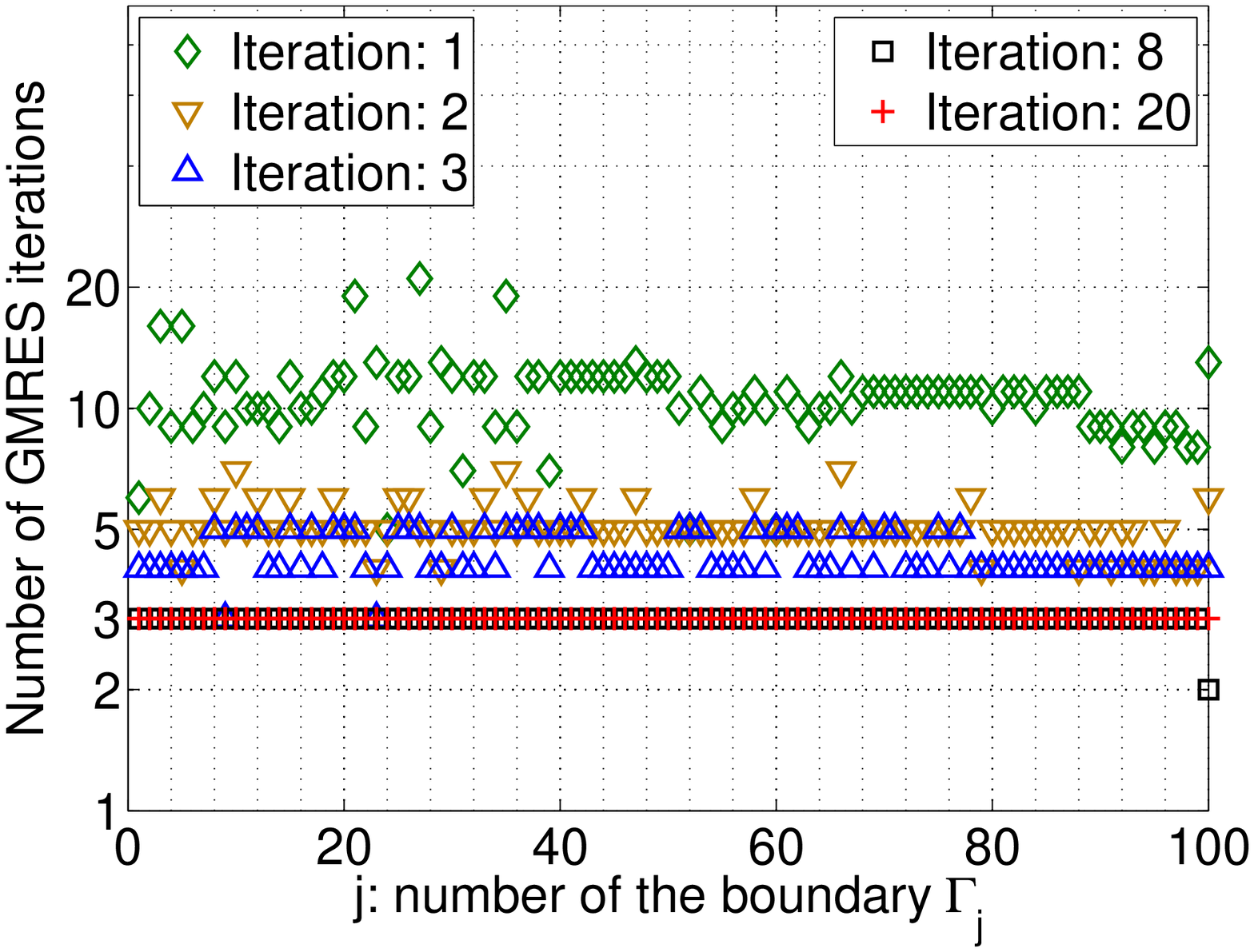}}
}
\caption{\rm Numerical results of Example 6 obtained with $n=512$. The successive iteration errors (left), the CPU time in seconds (middle), and the number of GMRES iterations for each boundary $\Gamma_j$, $j=1,2,\ldots,100$, for several iterations of Koebe's iterative method (right).} 
\label{f:ex6-err}
\end{figure}

\begin{example}\label{ex:6}{\rm
In this example, we compute the conformal mapping from bounded multiply connected region of connectivity $100$ bounded by $100$ ellipses. The numerical results are shown in Figures~\ref{f:ex6-im}--\ref{f:ex6-err}. 
}\end{example}

\begin{figure}%
\centerline{
\scalebox{0.3}{\includegraphics{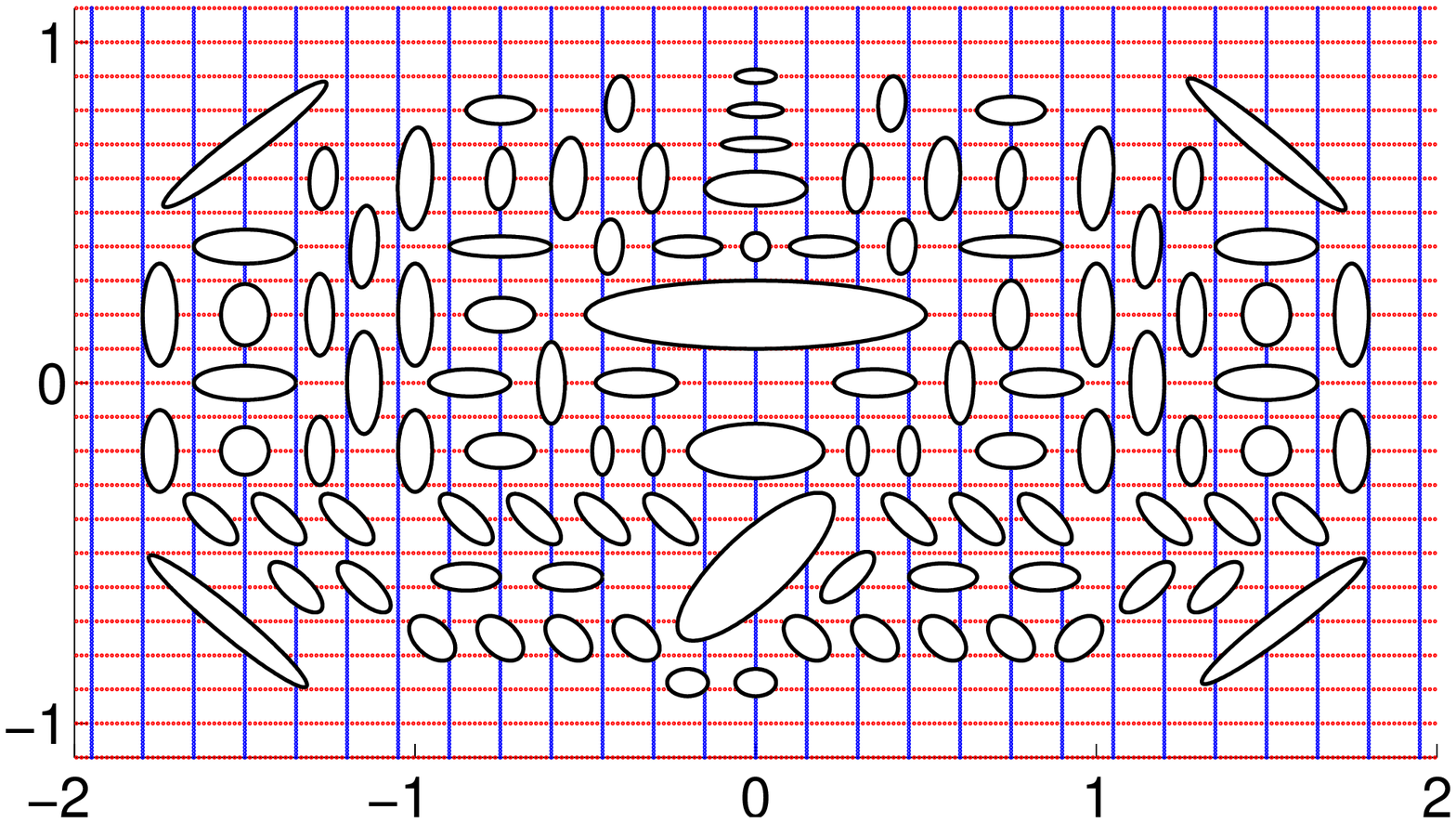}}
\scalebox{0.3}{\includegraphics{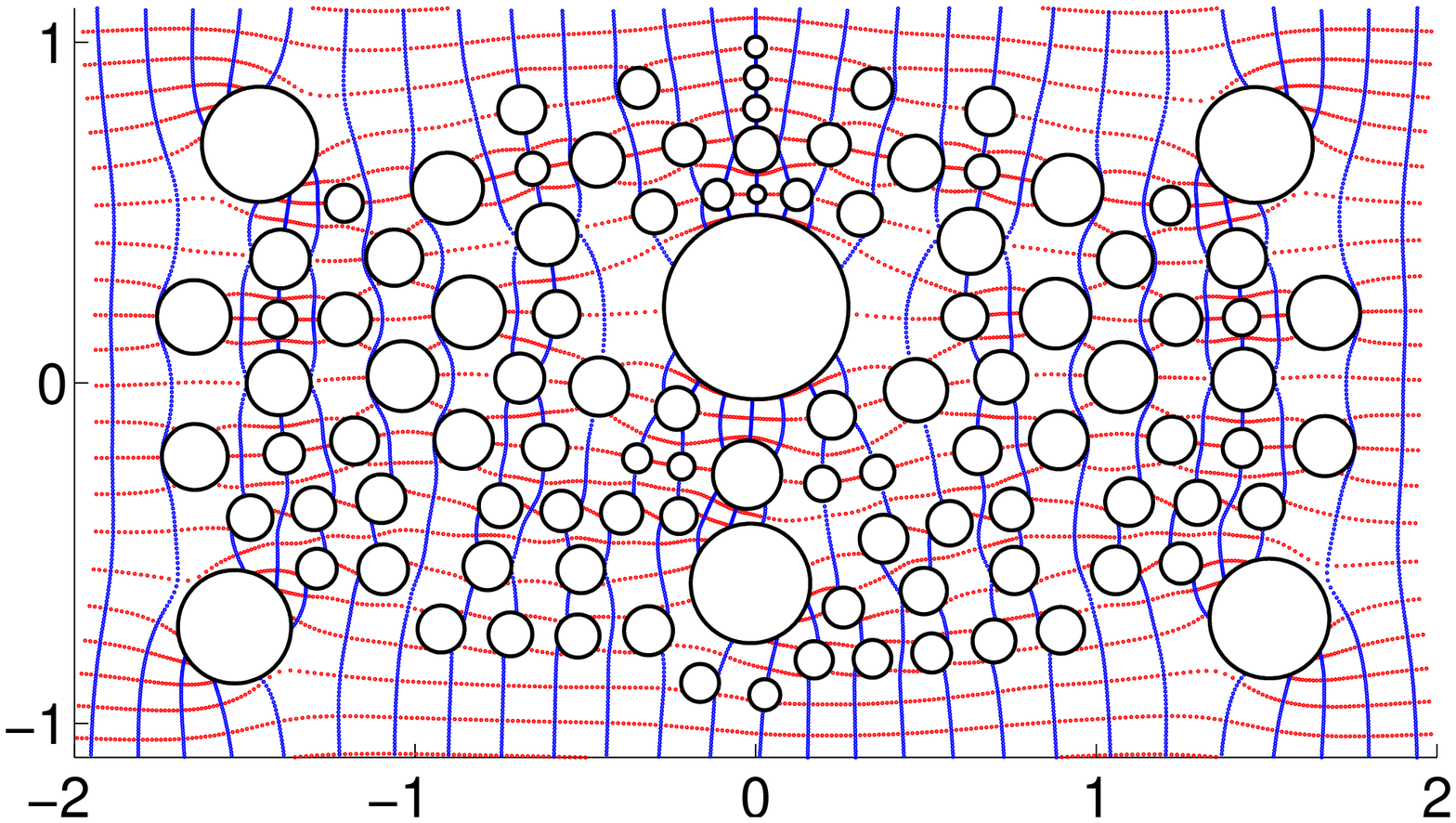}}
}
\caption{\rm The original region $G$ for Example 7 (left) and its image obtained with $n=256$ (right).} 
\label{f:ex7-im}
\end{figure}

\begin{figure}%
\centerline{
\scalebox{0.3}{\includegraphics{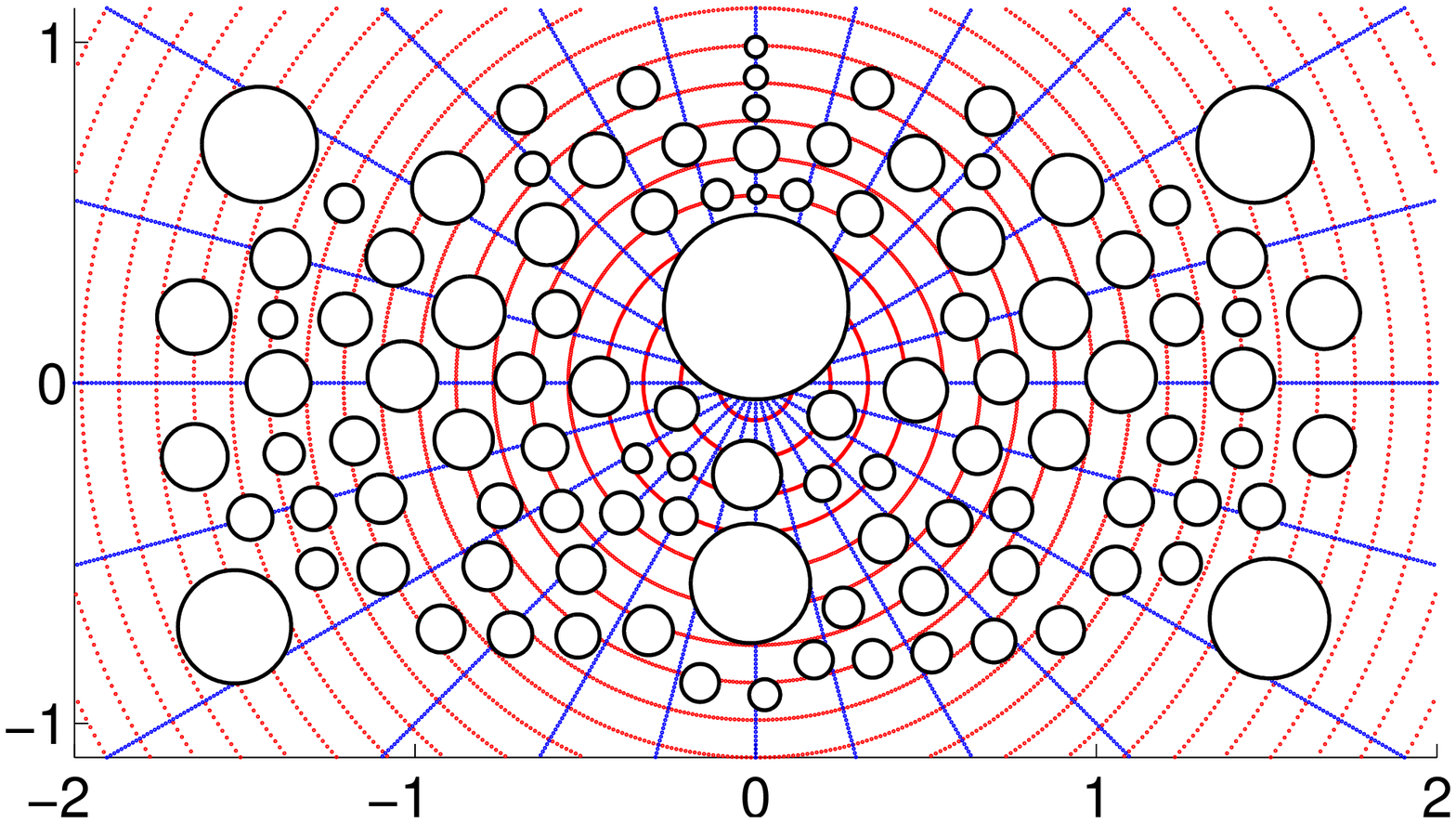}}
\scalebox{0.3}{\includegraphics{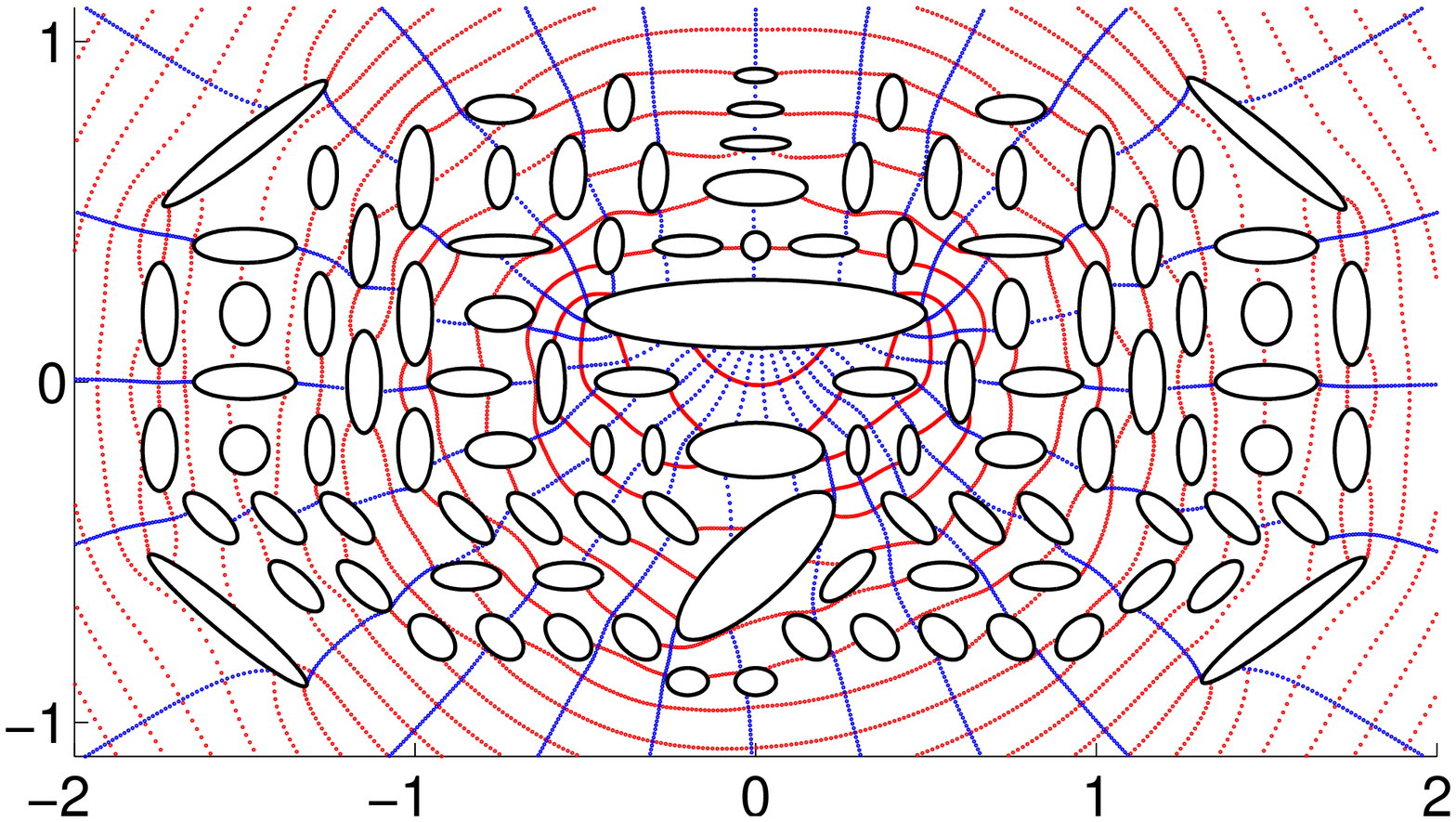}}
}
\caption{\rm The circular region $\Omega$ for Example 7 (left) and its inverse image obtained with $n=256$ (right).} 
\label{f:ex7-inv}
\end{figure}

\begin{figure}%
\centerline{
\scalebox{0.235}{\includegraphics{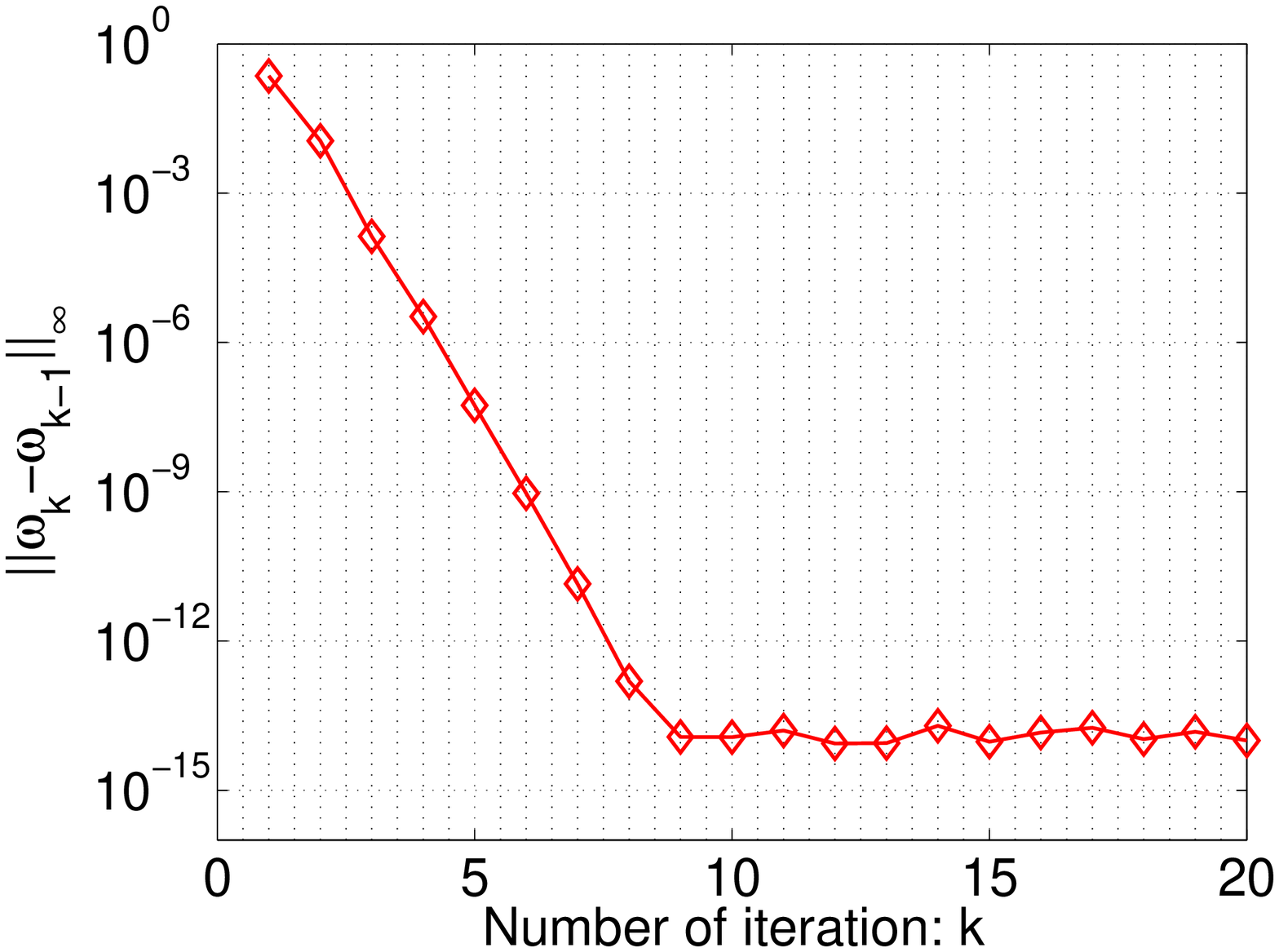}}
\scalebox{0.235}{\includegraphics{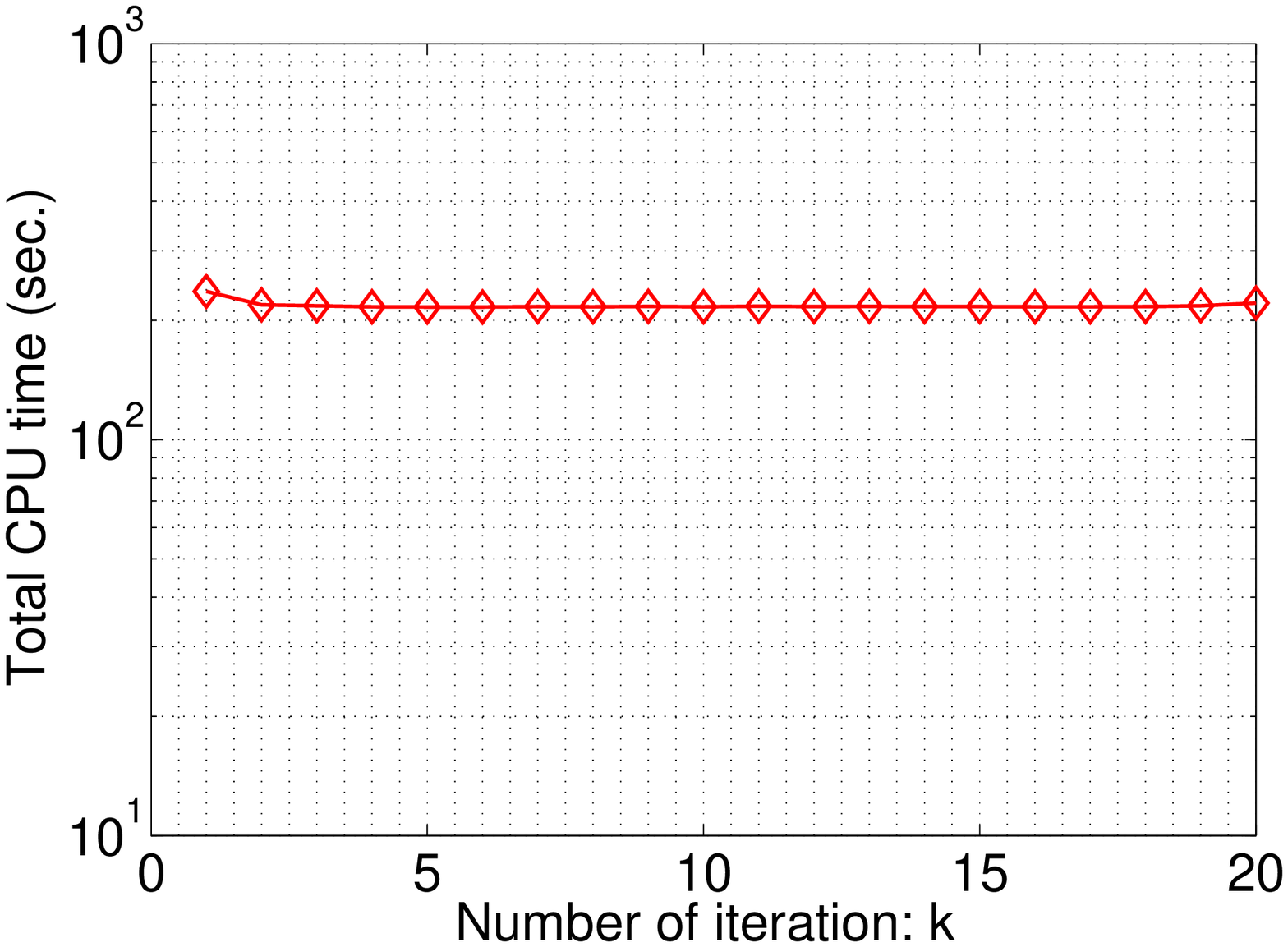}}
\scalebox{0.235}{\includegraphics{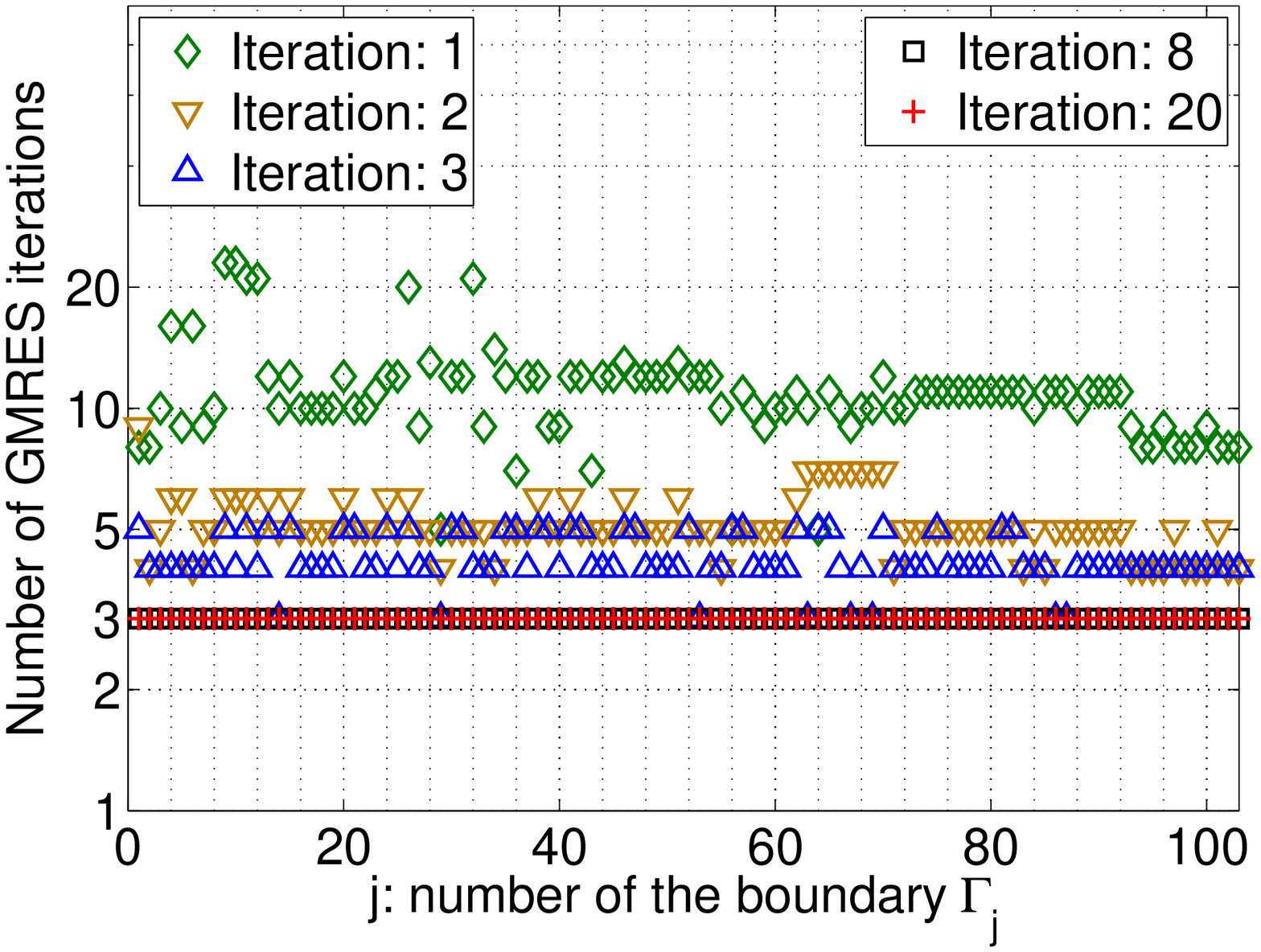}}
}
\caption{\rm Numerical results of Example 7 obtained with $n=256$. The successive iteration errors (left), the CPU time in seconds (middle), and the number of GMRES iterations for each boundary $\Gamma_j$, $j=1,2,\ldots,103$, for several iterations of Koebe's iterative method (right).} 
\label{f:ex7-err}
\end{figure}

\begin{example}\label{ex:7}{\rm
In this example, we compute the conformal mapping from an unbounded multiply connected region of connectivity $103$ bounded by $103$ ellipses. The numerical results are shown in Figures~\ref{f:ex7-im}--\ref{f:ex7-err}. 
}\end{example}

\begin{figure}%
\centerline{
\scalebox{0.24}{\includegraphics{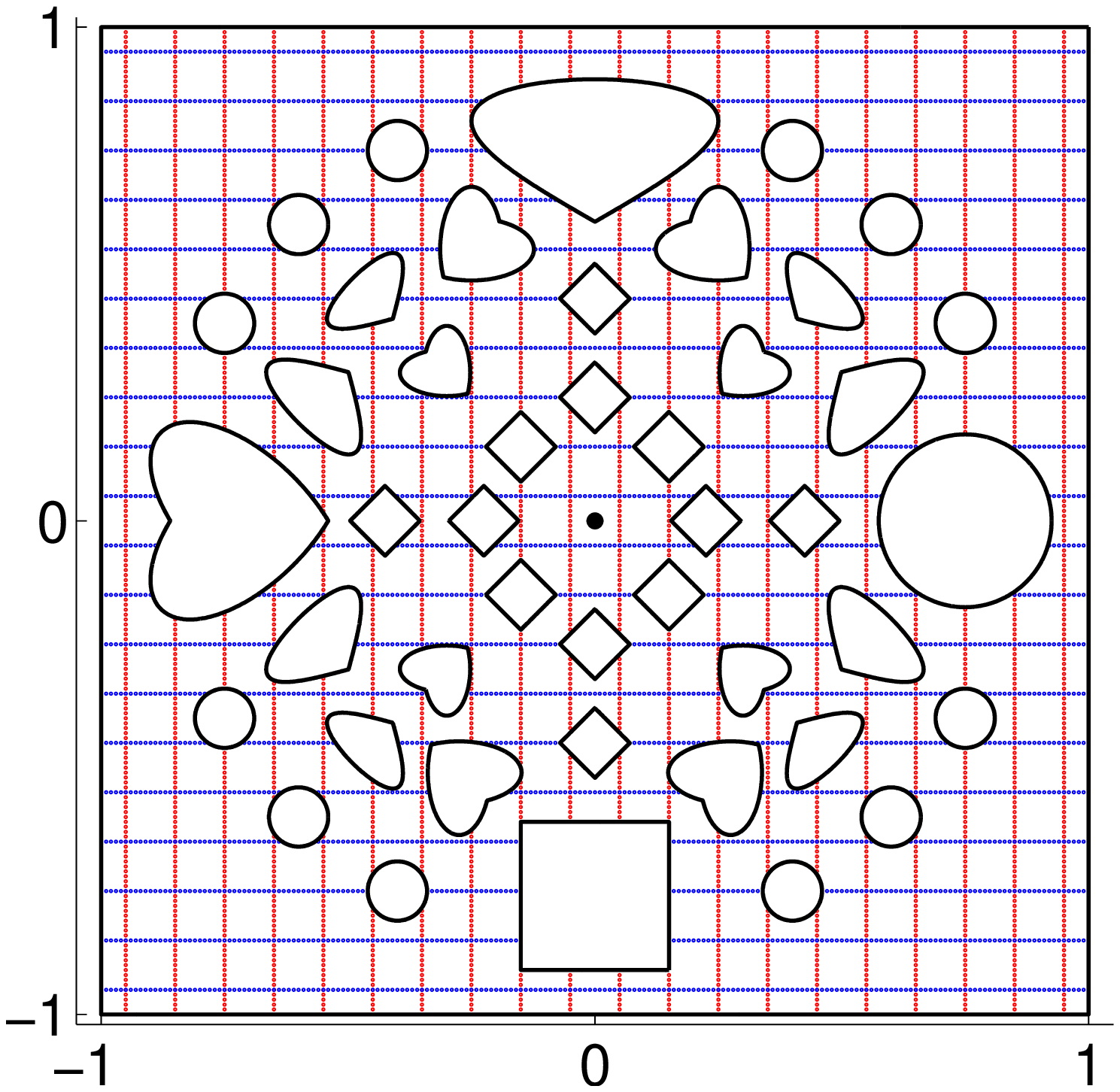}}
\scalebox{0.24}{\includegraphics{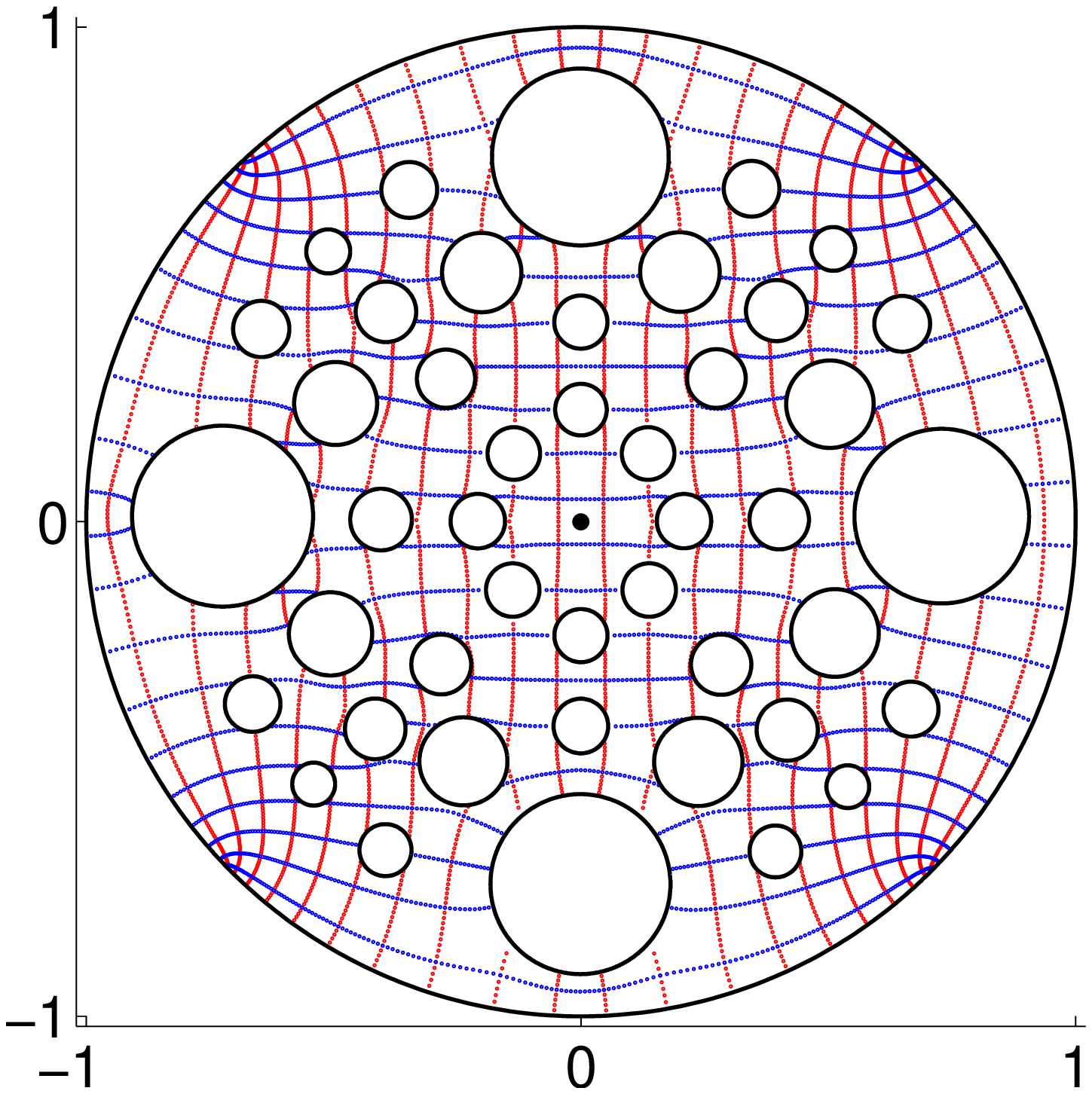}}
}
\caption{\rm The original region $G$ for Example 8 (left) and its image obtained with $n=1024$ (right).} 
\label{f:ex8-im}
\end{figure}

\begin{figure}%
\centerline{
\scalebox{0.24}{\includegraphics{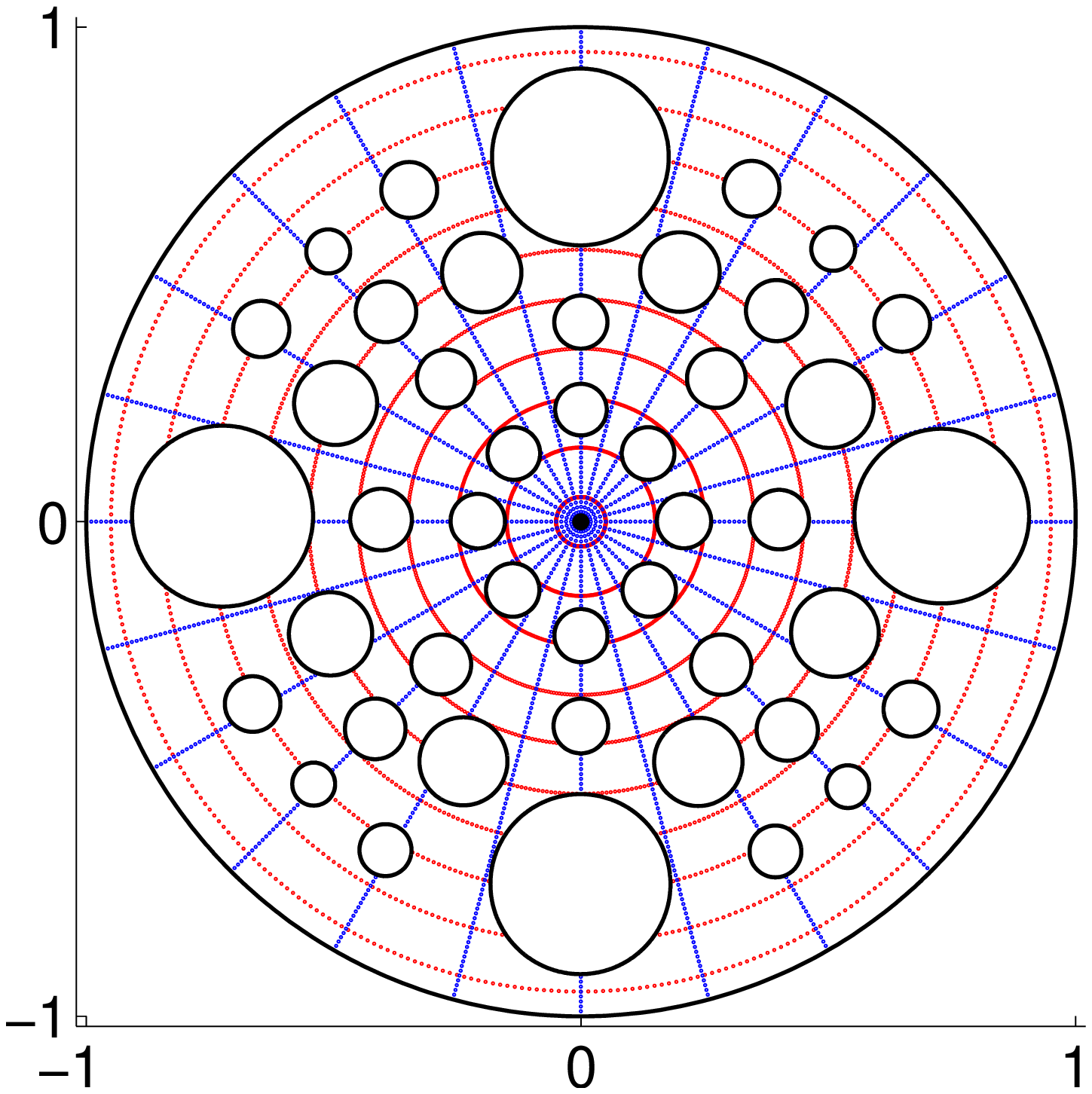}}
\scalebox{0.24}{\includegraphics{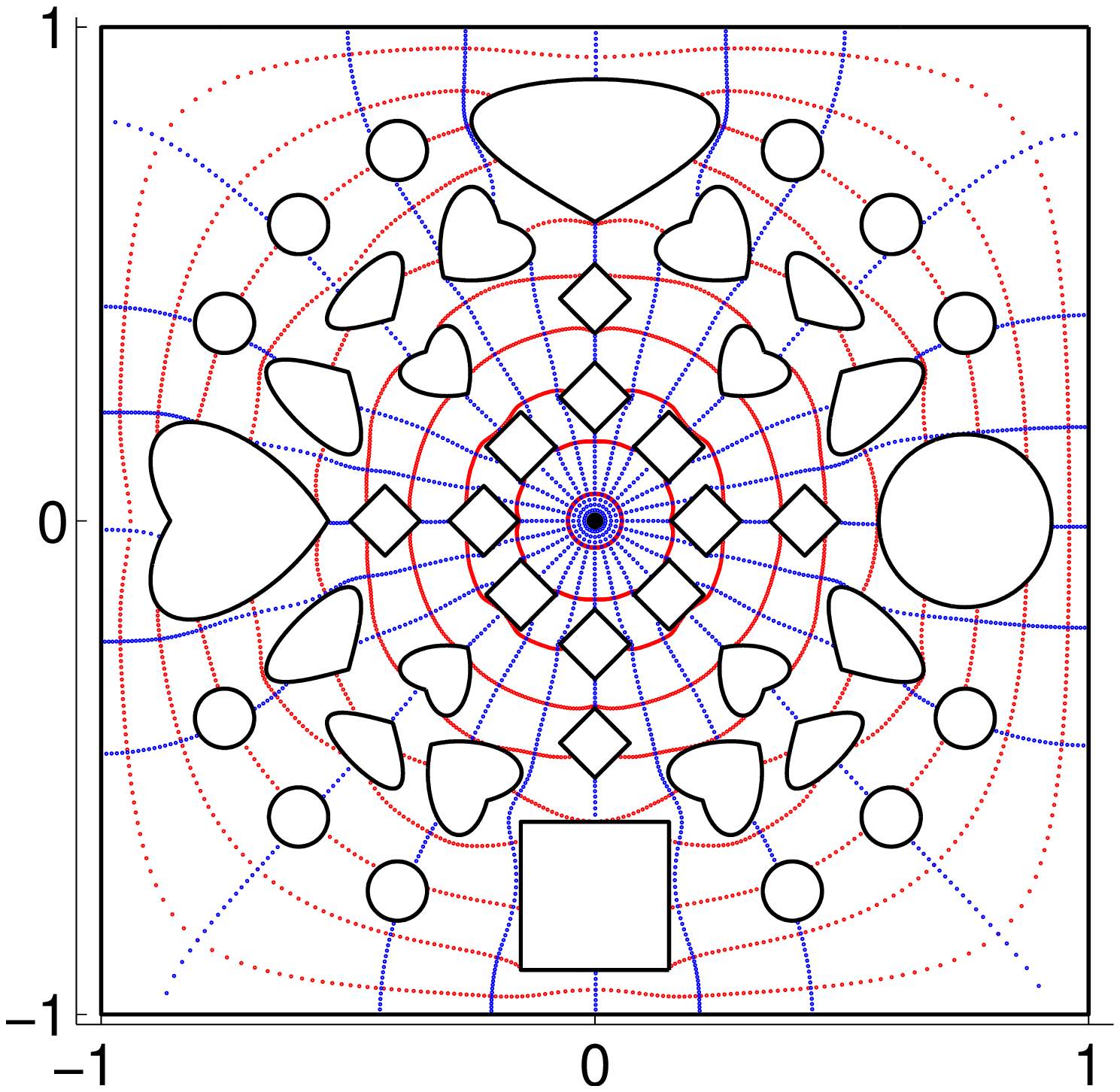}}
}
\caption{\rm The circular region $\Omega$ for Example 8 (left) and its inverse image obtained with $n=1024$ (right).} 
\label{f:ex8-inv}
\end{figure}

\begin{figure}%
\centerline{
\scalebox{0.235}{\includegraphics{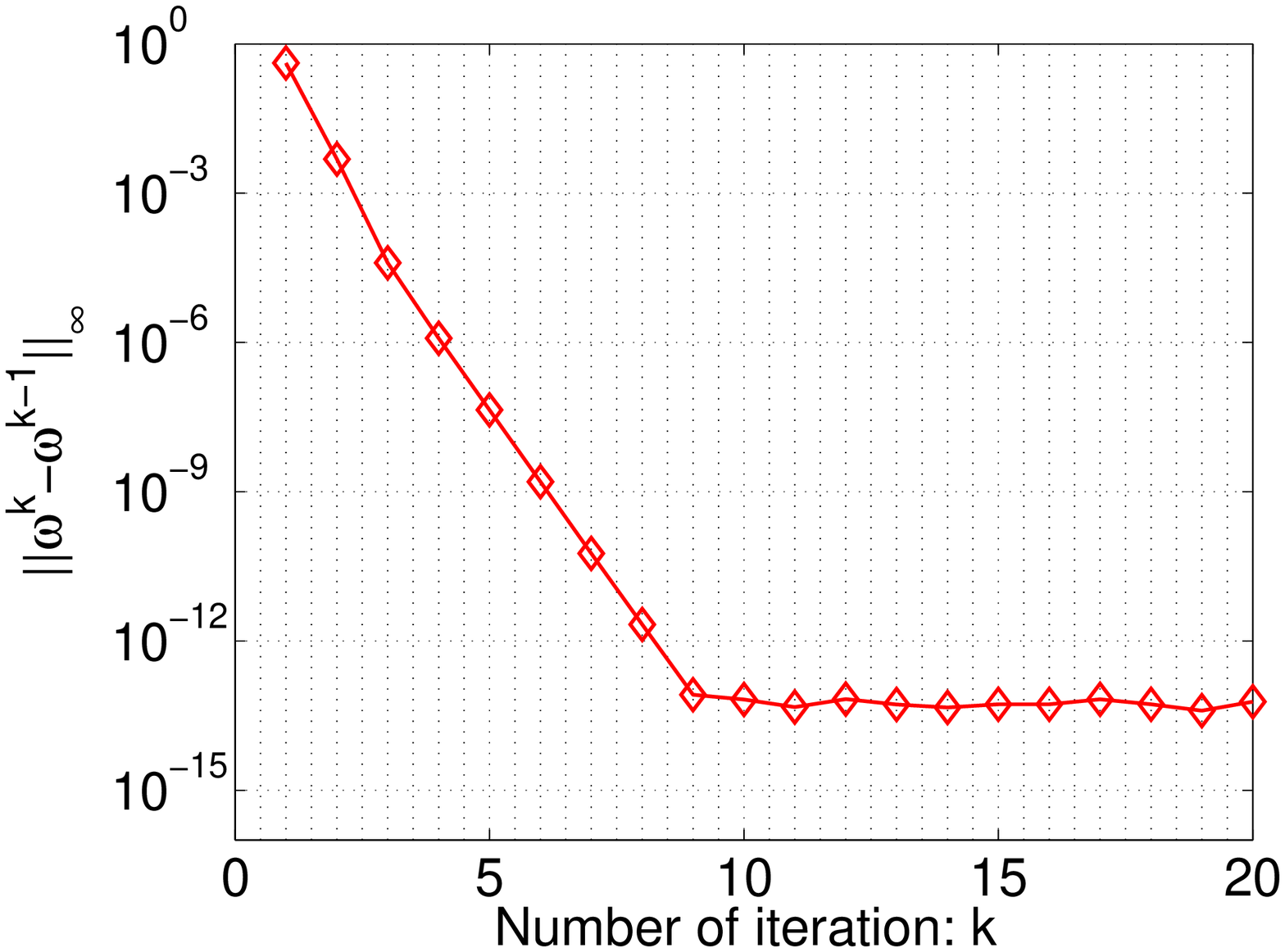}}
\scalebox{0.235}{\includegraphics{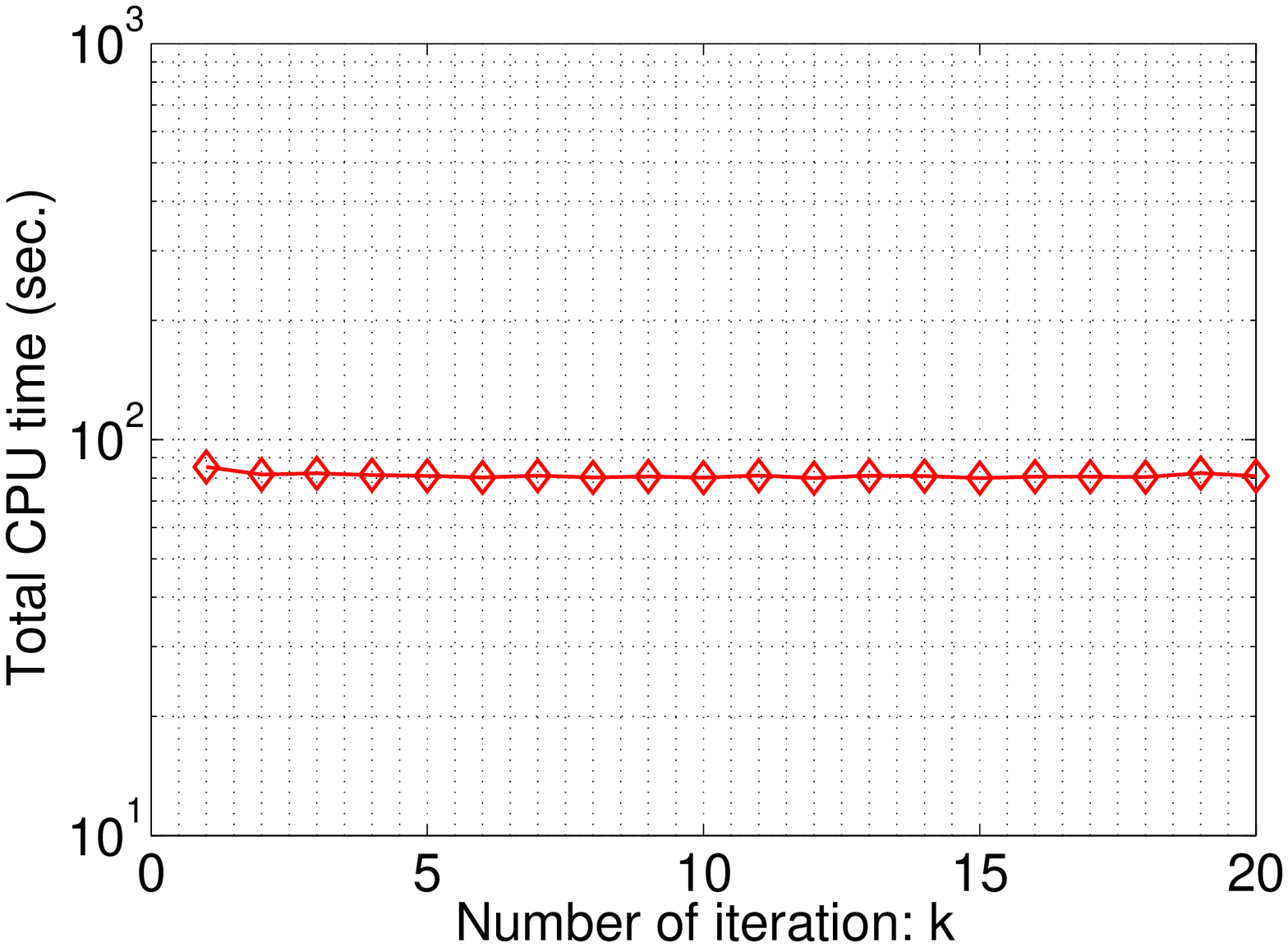}}
\scalebox{0.235}{\includegraphics{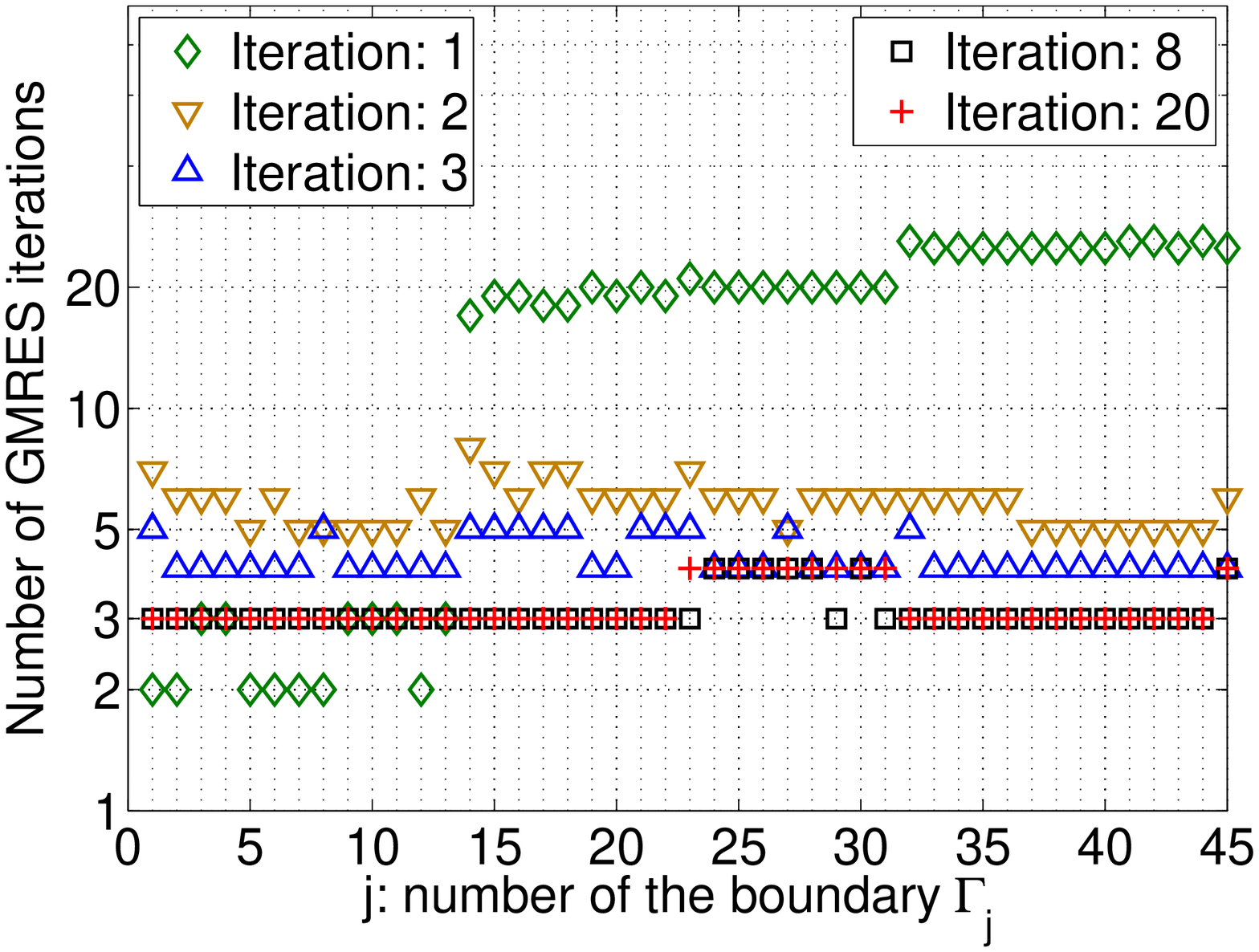}}
}
\caption{\rm Numerical results of Example 8 obtained with $n=1024$. The successive iteration errors (left), the CPU time in seconds (middle), and the number of GMRES iterations for each boundary $\Gamma_j$, $j=1,2,\ldots,45$, for several iterations of Koebe's iterative method (right).} 
\label{f:ex8-err}
\end{figure}

\begin{example}\label{ex:8}{\rm
In this example, we compute the conformal mapping from a bounded multiply connected region of connectivity $45$. The boundaries $\Gamma_1,\ldots,\Gamma_{14}$ are circles, the boundaries $\Gamma_{15},\ldots,\Gamma_{23}$ are piecewise smooth curves with one corner, the boundaries $\Gamma_{24},\ldots,\Gamma_{32}$ are piecewise smooth curves with two corners, and the boundaries $\Gamma_{33},\ldots,\Gamma_{45}$ are piecewise smooth curves with four corners. For this example, we discretize the integral equation by the trapezoidal rule with a graded mesh with grading parameter $p=3$ (see~\cite{Nas-fast}). The numerical results are shown in Figures~\ref{f:ex8-im}--\ref{f:ex8-err}. 
}\end{example}

\section{Conclusions}

Koebe's iterative method is a classical method for computing the conformal mapping of multiply connected regions onto circular regions. The method goes back to 1910~\cite{Koe10}. However, the implementation of the method was not simple. It was stated in~\cite{Kro14} that ``\emph{The coding complexity and the running inefficiency prevent it from broad practical applications}''. This paper presented a fast, an efficient, and easy to program numerical implementation of Koebe's iterative method to compute the circular map from bounded and bounded multiply connected regions of finite connectivity $m$. The computational cost of the presented method is $O(mn\ln n)$ where $n$ is the number of nodes in the discretization of each boundary component. However, the constant in the computational cost $O(mn\ln n)$ of the presented method is large compare to the constant in the computational cost $O(mn\ln n)$ for the method presented in~\cite{Nas-cmft09,Nas-siam09,Nas-jmaa11,Nas-jmaa13,Nas-siam13} for computing the conformal mapping onto canonical slit regions. This is because the computational cost of each iteration of Koebe's method is $O(mn\ln n)$.

The presented method can be used to compute the conformal mapping, its derivative, and its inverse. Thus, the presented method will be useful particularly for fluid problems which requires determining the conformal mapping and its derivative (see e.g.,~\cite{Cro-cyl,Cro-Anz,Cro-Kir}).

%

\end{document}